\documentclass[english]{amsart}
\usepackage{amsxtra}
\usepackage{mathtools}
\usepackage{amsmath}
\usepackage{amsthm}
\usepackage{stmaryrd}
\usepackage{amssymb}
\usepackage{amscd}
\usepackage{epsfig}
\usepackage{xypic}
\usepackage[pdftex,bookmarks,colorlinks,breaklinks=true]{hyperref}
\usepackage{xcolor}
\usepackage{manfnt}
\usepackage{verbatim}
\usepackage{graphicx}
\usepackage{xspace}
\usepackage{ifthen}
\usepackage{tabularx}
\usepackage[all,cmtip,line]{xy}
\usepackage{epigraph}
\usepackage{tikz-cd}           
\usepackage{mathabx}        
\usepackage{charter}

\newtheorem{thm}{Theorem}[section]
\newtheorem*{thm*}{Theorem}%[section]
\newtheorem{lem}[thm]{Lemma}
\newtheorem{prop}[thm]{Proposition}
\newtheorem{cor}[thm]{Corollary}
\newtheorem{defn}[thm]{Definition}

\theoremstyle{definition} 
\newtheorem{rem}[thm]{Remark}

\newtheorem{notation}[thm]{Notation}

\def\B{\mathcal B}
\def\D{\mathcal D}
\def\F{\mathcal F}
\def\G{\mathcal G}%WARNING:THIS IS NEW, PLEASE DO NOT DELETE
\def\L{\mathcal L}
\def\O{\mathcal O}
\def\P{\mathcal P}
\def\Q{\mathcal Q}
\def\U{\mathcal U}
\def\V{\mathcal V}
\def\W{\mathcal W}

\def\phs{\operatorname{PHS}}
\def\rad{\operatorname{rad}}
\def\Isom{\operatorname{Isom}}

\def\Ad{\rm{Ad}}
\def\Der{\operatorname{Der}}
\def\Tri{\operatorname{Tri}}

\def\cry{\mathit{Isoc}}
\def\Ob{\mathtt{Ob}}
\def\Rep{\mathit{Rep}}
\def\ton{\mathit{Tor}}

\def\ect{\mathit{Vec}}

\newcommand{\BA}{{\mathbb{A}}}
\newcommand{\BF}{{\mathbb{F}}}
\newcommand{\BG}{{\mathbb{G}}}
\newcommand{\BM}{{\mathbb{M}}}
\newcommand{\BN}{{\mathbb{N}}}
\newcommand{\BP}{{\mathbb{P}}}
\newcommand{\BQ}{{\mathbb{Q}}}
\newcommand{\BR}{{\mathbb{R}}}
\newcommand{\BZ}{{\mathbb{Z}}}

\newcommand{\Fa}{{\mathfrak{a}}}
\newcommand{\Fb}{{\mathfrak{b}}}
\newcommand{\Fc}{{\mathfrak{c}}}
\newcommand{\Fg}{{\mathfrak{g}}}
\newcommand{\Fh}{{\mathfrak{h}}}

\newcommand{\Fp}{{\mathfrak{p}}}
\newcommand{\Fu}{{\mathfrak{u}}}
\newcommand{\Fv}{{\mathfrak{v}}}
\newcommand{\FF}{{\mathfrak{F}}}

\newcommand{\ol}{\overline}
\newcommand{\ul}{\underline}

\newcommand{\pbar}{{\bar P}}

\newcommand{\gbar}{{\bar G}}
\newcommand{\xbar}{\bar X}

\def\Grp{\mathop{\rm \mathfrak{Grp}}\nolimits}
\def\Scheme{\mathop{\rm\mathfrak{Sch}}\nolimits}
\def\Flat{\mathop{\rm\mathfrak{Flat}}\nolimits}

\def\Sym{\mathop{\rm Sym}\nolimits}

\def\Aut{\mathop{\rm Aut}\nolimits}
\def\uAut{\mathop{\underline{\rm Aut}}\nolimits}
\def\Frob{\mathop{\rm Frob}\nolimits}
\def\Hom{\mathop{\rm Hom}\nolimits}

\def\Ens{\mathop{\rm Ens}\nolimits}
\def\End{\mathop{\rm End}\nolimits}

\def\Spec{\mathop{\bf Spec}\nolimits}

\def\image{\mathop{\rm im}\nolimits}

\def\Tor{\mathop{\rm Tor}\nolimits}

\def\Lie{\mathop{\rm Lie}\nolimits}

\def\intaut{\mathop{\rm int}\nolimits}

\def\diag{\mathop{\rm diag}\nolimits}

\def\proj{\mathop{\rm pr}\nolimits}
\def\GL{\mathop{\rm GL}\nolimits}
\def\Grass{\mathop{\rm Grass}\nolimits}

\def\red{{\rm red}}

\def\tor{{\rm tor}}

\def\alg{{\rm alg}}
\def\id{{\rm id}}

\def\innHom{\underline{\Hom}}

\def\into{\hookrightarrow}
\def\onto{\twoheadrightarrow}
\def\isoto{\arrover{\cong}}

\newbox\mybox
\def\arrover#1{\mathrel{
       \setbox\mybox=\hbox spread 1em
              {\hfil$\scriptstyle#1\vphantom{g}$\hfil}
       \vbox{\offinterlineskip\copy\mybox
             \hbox to\wd\mybox{\rightarrowfill}}}}             
\def\invlim{\mathop{\vtop{\hbox{\rm lim}\vskip-8pt
        \hbox{\hskip1pt$\scriptstyle\longleftarrow$}\vskip-1pt}}}
\def\dirlim{\mathop{\vtop{\hbox{\rm lim}\vskip-8pt
        \hbox{\hskip1pt$\scriptstyle\longrightarrow$}\vskip-1pt}}}
\def\ontoover#1{\mathrel{
       \setbox\mybox=\hbox spread 1.4em{\hfil$\scriptstyle#1$\hfil}
       \vbox{\offinterlineskip\copy\mybox
             \hbox to\wd\mybox{\rightarrowfill\hskip-2.8mm
                               $\rightarrow$}}}}

\title{$(G,\mu)$-Windows and Deformations of $(G,\mu)$-Displays}

\author[O. B\"ultel]{Oliver B\"ultel}
\address{O. B\"ultel: Department of Mathematics and Feza G\"ursey Center for Physics and Mathematics\\
Bo\u{g}azi\c{c}i University\\Istanbul\\Turkey}
\author[M. H. Hedayatzadeh]{S. Mohammad Hadi Hedayatzadeh}
\address{M. H. Hedayatzadeh: School of Mathematics\\ Institute for Research in Fundamental Sciences\\Tehran\\Iran }

\date{}
\begin{document}
\maketitle
\parindent0pt	
\begin{abstract}
Let $k_0$ be a finite field of characteristic $p$, let $G$ be a smooth affine group scheme over $\BZ_p$, and let $\mu$ be a cocharacter of $G_{W(k_0)}$ 
such that the set of $\mu$-weights of $\Lie G$ is a subset of $\{-1,0,1,2,\dots\}$. We prove that the groupoid of adjoint nilpotent $(G,\mu)$-displays as in 
\cite[Definition 1.0.1]{pappas}, \cite{lau2} and \cite[Definition 3.16]{E7} is equivalent to the groupoid of $(G,\mu)$-windows, which are the generalizations of windows as in \cite[Definition 5.1]{langer2} and \cite[Definition 3.3.1]{E7}.
\end{abstract}
	
\tableofcontents	
	
\section{Introduction}

Let $p$ be a prime number. The study of $p$-divisible groups by means of semi-linear algebra over the ring of $p$-typical Witt vectors has a long history. Its 
undisputed point of departure is a perfect field $k\supset\BF_p$, over which one has a classification in terms of the category of \emph{Dieudonn\'e modules} 
of which an object consists of a finitely generated free $W(k)$-module together with a Frobenius-linear endomorphism whose cokernel is annihilated by $p$. 
The theory owes its ease and popularity to $W(k)$ being a complete discrete valuation ring with uniformizing element $p$ and residue field $k$. Now let $A$ 
be a $p$-adically separated and complete commutative ring without $p$-torsion and suppose that some $\sigma\in\End(A)$ lifts the Frobenius endomorphism:
\begin{align*}
A/pA&\to A/pA\\x&\mapsto x^p
\end{align*}

Furthermore, let $\Fa\subset A$ be a $p$-adically open pd-ideal (in the case at hand this means that $pA\subset\sqrt \Fa$ and $x^p\in p\Fa$ for every $x\in \Fa$). 
One calls $(A,\Fa,\sigma)$ a \emph{frame} over the ring $R=A/\Fa$, in which $p$ is nilpotent. One knows that the following three categories are equivalent:
\begin{itemize}
\item
\emph{windows} for the frame $(A,\Fa,\sigma)$ in the sense of \cite{zink1}
\item
nilpotent \emph{displays} over $R$ in the sense of \cite{zink2}.
\item
$p$-divisible formal groups over $R$
\end{itemize}

Please see loco citato and in particular \cite{messing2}, \cite{lau1} for the proofs. A display over $R$ is a quadruple $(P,Q,F,V^{-1})$ with the following properties:
\begin{itemize}
\item
$P$ is a finitely generated projective $W(R)$-module and $Q\subset P$ is the inverse image of a 
direct summand in $R\otimes_{w_0,W(R)}P$, where $w_0:W(R)\onto R$ is the $0$-th ghost coordinate,
\item
$V^{-1}:Q\rightarrow P$ is an $F$-linear homomorphism whose image generates $P$ as a $W(R)$-module, and $F:P\rightarrow P$ is an $F$-linear homomorphism satisfying
\begin{equation}
\label{axiom01}
V^{-1}({^Va}\cdot x)=aF(x)
\end{equation}
for all $a\in W(R)$ and $x\in P$, where $V:W(R)\rightarrow W(R);\,(a_0,a_1,\dots)\mapsto(0,a_0,\dots)$ is the Verschiebung, so that ${^Va}\cdot x\in Q$.\\
\end{itemize}

A window is a similar semi-linear structure in which the projection $A\onto R$ replaces the $0$-th ghost coordinate $w_0:W(R)\onto R$, while the 
Frobenius lift $\sigma$ is a substitute for the absolute Frobenius on $W(R)$, please see \cite{zink1} for more information e.g. on the nilpotence conditions 
in these categories. The axiom \eqref{axiom01} should be regarded as a sharpening of the classical formula $F\vert_Q=pV^{-1}$, inspired by older theories 
on ``strongly divisible lattices'' and the like, please see \cite{fontaine1} and \cite{fontaine2} for more information. It is known that $W(R)$ is $p$-adically 
separated and complete, but it is hard to unravel its structure as a commutative ring, unless $R$ is a perfect field (or a product of these). The equivalence 
between windows and nilpotent displays can be viewed as a remedy for this, as it gives some freedom to work over rings which are nicer than $W(R)$.\\

So given that each of these three categories has its own distinct features, we gain a powerful combination of their advantages from equivalences 
between them: Within the category of $p$-divisible groups, for instance, it is in general not possible to form a tensor product, but the categories of 
displays (resp. windows) can be regarded as full subcategories of the natural rigid $\otimes$-categories of Langer-Zink displays (resp. windows), which 
showed up for the first time in \cite[Definition 2.4-2.6/Definition 5.1]{langer2}. This phenomenon has proved to be a fertile inspiration for the recent 
works \cite{HF1}, \cite{H2} and \cite{G2}. Spurred by this, along with \cite{pappas} and \cite{lau2}, this paper addresses the comparison 
of various concepts of additional structure in the aforementioned three categories, perhaps with an eye towards applications to the numerous moduli 
spaces that are indigenous to the Langlands program, i.e. Rapoport-Zink spaces, period domains, Shimura varieties and their integral models.\\ 

More specifically, our approach to displays with additional structure runs as follows: We fix a smooth affine $\BZ_p$-group $G$ together with a 
cocharacter $\mu$ which is defined over the ring of integers in an unramified $\BQ_p$-extension of finite degree, i.e. over $W(k_0)$, where $k_0$ is a 
finite field of characteristic $p$. Let $P_\mu\subset G_{W(k_0)}$ be the associated parabolic subgroup, so that $\Lie P_\mu$ is the subspace 
of $\Lie G$ whose $\mu$-weights are non-negative. Assuming that $-1$ is the only negative $\mu$-weight of $\Lie G$, we define a key homomorphism 
$$L^+G_{W(k_0)}\hookleftarrow H_\mu\stackrel{\Phi_\mu}{\rightarrow}L^+G_{W(k_0)}$$
where $L^+G$ stands for the affine group scheme with $L^+G(R)=G(W(R))$, $H_\mu$ for the inverse image of $P_\mu$ via the $0$-th ghost map 
$w_0^G:L^+G\onto G$ and $\Phi_\mu$ for the integral model of the group homomorphism $g\mapsto{^F(\mu(p)g\mu(\frac1p))}$. The existence 
of $\Phi_\mu$ is not a triviality and was proved in \cite{pappas} and \cite{E7} only under additional assumptions, most notably that $G$ be reductive. 
In this paper we give an improved exposition, using nothing but $-1$ being a lower bound for the $\mu$-weights of $\Lie G$. At last we let $\B(G,\mu)$ 
be the fpqc quotient stack of $L^+G_{W(k_0)}$ by a certain $H_\mu$-action defined thereon, namely the so-called $\Phi_\mu$-conjugation:
$$L^+G_{W(k_0)}\times_{W(k_0)}H_\mu\rightarrow L^+G_{W(k_0)};\,(U,h)\mapsto h^{-1}U\Phi_\mu(h)$$

According to \cite{pappas} and \cite{E7} a \emph{$(G,\mu)$-display} $\D$ over some $W(k_0)$-scheme $X$ is simply defined to be a 
$1$-morphism from $X$ to $\B(G,\mu)$. Besides $\D$ is called \emph{banal} if the said $1$-morphism factors through the obvious 
presentation $L^+G_{W(k_0)}\rightarrow\B(G,\mu)$. Under the same requirement on $G$ and the $\mu$-weights of $\Lie G$ we suggest a 
natural definition of the groupoid of \emph{$(G,\mu)$-windows} for the frame $(A,\Fa,\sigma)$ and we prove the following theorem:

\begin{thm*}[\ref{main01}]
The groupoid of $(G,\mu)$-windows over the frame $(A,\Fa,\sigma)$ is equivalent to the groupoid of adjoint nilpotent $(G,\mu)$-displays over $A/\Fa$.
\end{thm*}

The adjoint nilpotence is an analog of Zink's nilpotence condition. Unsurprisingly our proof of theorem \ref{main01} has close ties to the study of lifts of 
$(G,\mu)$-displays from $A/\Fa$ to $A$. At this point we invoke the deformation theory established by \cite[Theorem 3.5.4]{pappas}. In fact we have to reconsider 
it slightly, because the proof of loc. cit. (and of similar applicable lemmas in \cite{lau2} and \cite{E7}) seems to require the $p$-adic topological nilpotency of the 
pd-ideal $\Fa$, which we do not want to impose in theorem \ref{main01}. Apart from this technicality the proof of theorem \ref{main01} runs smoothly in an assembly 
line, which roughly speaking works as follows: 
\begin{eqnarray*}
&&(G,\mu)\mbox{-windows over }(A,\Fa,\sigma)\rightsquigarrow\\
&&(G,\mu)\mbox{-windows over }(W(A),W(\Fa)+I(A),F)\rightsquigarrow\\
&&(G,\mu)\mbox{-triples over }(A,\Fa)\rightsquigarrow\\
&&(G,\mu)\mbox{-displays over }A/\Fa
\end{eqnarray*}

For the definition of the above notion of \emph{$(G,\mu)$-triple} over $(A,\Fa)$ we have to refer to the body of the text; It is modeled on Zink's 
notion for $\GL(h)$. In fact the stack of displays of dimension $d$ and height $h$ is just our $\B(\GL(h)_{\BZ_p},\mu_{h,d})$ for the cocharacter: 
\begin{equation*}
\mu_{h,d}:\BG_m\rightarrow\GL(h);z\mapsto\diag(\overbrace{1,\dots,1}^d,\overbrace{z,\dots,z}^{h-d})
\end{equation*}
This normalization follows \cite{pappas}, which is somewhat opposite to the ones of \cite{lau2} and \cite{E7}, but do notice that these stacks neglect 
non-invertible homomorphisms between displays. The usage of $(W(A),W(\Fa)+I(A),F)$ is suggested by \cite{zink1}, of which the main result can be 
reobtained from theorem \ref{main01} by specializing to the case $(G,\mu)=(\GL(h)_{\BZ_p},\mu_{h,d})$.\\

Our second result concerns the rigidity of quasi-isogenies between $(G,\mu)$-displays, for which it is necessary to work over the localization of the 
Witt ring at the prime $p$. This leads to an interesting $1$-morphism from $\B(G,\mu)$ to its rational analog $\B(G)$, which one can think of as the 
sheafification of the prestack whose fiber over $R$ is given by the equivalence relation of $F$-conjugacy on the set $G(\BQ\otimes W(R))$. Thus 
a quasi-isogeny between $(G,\mu)$-displays $\D_1$ and $\D_2$ can be regarded as an isomorphism $\D_1^0\cong\D_2^0$ in $\B(G)$, where 
$\D^0$ stands for the image of some $(G,\mu)$-display $\D$ under the canonical map $\B(G,\mu)\rightarrow\B(G)$. In the following result we write $\D_x$ for the pull-back of a $(G,\mu)$-display $\D$ along the morphism 
$\Spec\kappa(x)\rightarrow X$, where $\kappa(x)$ denotes the residue field of the stalk $\O_{X,x}$, for any $x\in X$:

\begin{thm*}[\ref{rigid01}]
Let $\D$ be a $(G,\mu)$-display over a $W(k_0)$-scheme $X$ whose restriction to the special fiber $X\times_{W(k_0)}k_0$ is adjoint nilpotent.
\begin{itemize}
\item[(i)]
If $X$ is connected and locally noetherian and its structural morphism to $\Spec W(k_0)$ is assumed to be \emph{closed} in the sense of \cite[(2.2.6)]{egai}, then
a self quasi-isogeny $\alpha\in\Aut(\D^0)$ agrees with the identity if and only if this holds for at least one of its specializations $\alpha_x\in\Aut(\D_x^0)$.
\item[(ii)] 
If $p$ is assumed to be locally nilpotent on $\O_X$, then the natural map $\Aut(\D)\rightarrow\Aut(\D^0)$ is 
injective (i.e. $\D\mapsto\D^0$ is faithful).
\end{itemize}
\end{thm*}

The proof of this result is quite technical and we therefore omit an overview here. We should only mention that we make a careful study 
of the loop group functor $LG(R)=G(\BQ\otimes W(R))$, which turns out to be a sheaf for the fpqc topology, in fact it is an ind-scheme. 
We expect that theorem \ref{rigid01} has consequences for the description of the period map of Rapoport-Zink spaces in the setting of \cite{pappas}.\\

We note in passing that \cite[Definition 3.3]{rapoportrichartz} defines an $F$-isocrystal with $G$-structure 
over a connected $\BF_p$-scheme $S$ to be a faithful $\BQ_p$-linear exact $\otimes$-functor 
$$M:\Rep_{\BQ_p}(G)\rightarrow{F-\cry(S)},$$
where $\Rep_{\BQ_p}(G)$ is the $\otimes$-category of finite-dimensional representations of $G_{\BQ_p}$ while ${F-\cry(S)}$ is the category of 
$F$-isocrystals on $S$ in the sense of \cite[section 2.1]{katz} and \cite[section 3.1]{rapoportrichartz}. Loc. cit. attributes the underlying tannakian 
formalism to \cite{simp}, and it might suggest to call a faithful $\BZ_p$-linear exact $\otimes$-functor from the category of $G$-representations 
to the category of Langer-Zink displays a ``Langer-Zink display with $G$-structure''. Please see \cite{daniels} for more precise versions and more 
information on the relationship between the tannakian approach to ours. Please see \cite[\S 2.1]{H2} for more information on the relationship 
between multi-linear maps of displays and multilinear maps of formal groups.\\

Let us give further evidence for the relevance of $\B(G,\mu)$: 

\begin{itemize}
\item[(i)]
Every $(G,\mu)$-display over a perfect field $k\supset k_0$ gives rise to an $F$-isocrystal with $G$-structure over $k$ in the sense of 
\cite[Definition 3.3]{rapoportrichartz}. In particular one obtains a function from the set of ``points'' of $\B(G,\mu)_{k_0}$ (say in the sense of 
\cite[D\'efinition (5.2)]{champions}) to the important set $B(G)$ introduced and studied by Kottwitz and others (cf. \cite{Kott2}, \cite{Kott1} and references therein). 
\item[(ii)]
Every adjoint nilpotent $(G,\mu)$-display $\D$ over a perfect field $k\supset k_0$ possesses a universal formal deformation 
$\tilde\D$ over the $W(k_0)$-algebra\\ $W(k)[[t_1,\dots,t_r]]$, where $r$ is the codimension of $P_\mu$ in $G_{W(k_0)}$.
\end{itemize}

We round off our results with a few further applications of the deformation theory of adjoint nilpotent $(G,\mu)$-displays; Here is a sample:

\begin{thm*}[\ref{compl}]
Let $A$ be a $W(k_{0})$-algebra and $\D$ a $(G,\mu)$-display over $A/\Fa$, where the ideal $\Fa\subset A$ is assumed to satisfy one of the following:
\begin{itemize}
\item[(i)]
$\Fa$ is $(pA+\Fa)$-adically separated and complete.
\item[(ii)]
$\Fa$ allows a divided power structure and is $p$-adically separated and complete.
\end{itemize}

Assume that $\D\times_{A/\Fa}A/(pA+\Fa)$ is adjoint nilpotent. Then, for each lift $\tilde\D$ of $\D$ the reduction 
of automorphisms gives an injection $\Aut(\tilde\D)\into\Aut(\D)$. Moreover there exists at least one such lift.
\end{thm*}
	
Under the assumptions of this theorem it is meaningful to ask for the set of deformations of some $(G,\mu)$-display $\D$ over $A/\Fa$, and our fourth result 
describes this set rather explicitly, provided that we are in the scenario (ii) of theorem \ref{compl}: Fix a lift $\tilde\D$ of $\D$ and consider the pointed $A$-scheme
$$X_{\tilde\D}:=q_0(\tilde\D)*^{P_\mu}G_{W(k_0)}/P_\mu,$$
where $q_0(\tilde\D)$ is a certain $P_\mu$-torsor which we define in the body of the text, our $X_{\tilde\D}$ is a twist of 
the homogeneous space $G_{W(k_0)}/P_\mu$ over $A$. Our result states that lifts of $\D$ are parametrized by the pointed set
$$X_{\tilde\D}(\Fa)=\{x\in X_{\tilde\D}(A)\mid x\equiv1\mod \Fa\},$$ 
i.e. the kernel of $X_{\tilde\D}(A)\rightarrow X_{\tilde\D}(A/\Fa)$. We prove this using the aforementioned 
notion of $(G,\mu)$-triple over $(A,\Fa)$.\\

Our results require a significant amount of background material on torsors, Greenberg transforms, loop groups, arithmetic jet schemes, cocharacters of smooth group schemes, and fpqc-sheaves, which we explain in 
section \ref{admissible08} and in the appendix. We also give as many references as possible and attempt a rather broad generality, which manifests itself in:\\

\begin{itemize}
\item[(i)]
Inclusion of non-reductive smooth affine group schemes, whenever this is possible.
\item[(ii)]
Preference of arbitrary base schemes over merely affine ones.
\item[(iii)]
Avoidance of noetherian (or similar finiteness) assumptions, whenever this is possible.
\item[(iv)]
Occurrences of completions of topological ground rings only with respect to fine as possible a linear topology.\\
\end{itemize}

\textbf{Acknowledgements.} We would like to thank Sebastian Bartling, Chuangxun Cheng, Claudia Glanemann, Simon H\"aberli, Aliakbar Hosseini, George Pappas and Ali Partofard for 
helpful conversations and the IPM for excellent working conditions. In addition, the first named author would like to thank the audience of an introductory 
course on algebraic groups, which was given there on three days in June 2019. 		

\section{Preliminaries}
\label{admissible08}

Recall that a topology on a commutative and unital ring $A$ is called linear if it allows a subbase of the form $\{x+\Fa_i\mid x\in A,\,i\in I\}$ where $I$ is an index set and
\begin{equation}
\label{admissible04}
\{\Fa_i\}_{i\in I}
\end{equation} 
a family of ideals. The index set $I$ can be taken to be $\BN_0$ if and only if the resulting topology satisfies the first axiom 
of countability, in which case the generating family \eqref{admissible04} can be taken to be a descending chain of ideals:
\begin{equation}
\label{admissible05}
A\supset\Fa_0\supset\Fa_1\supset\Fa_2\supset\dots
\end{equation}

An ideal $J$ in a linearly topologized ring $A$ is called topologically nilpotent if $\Fa+J/\Fa$ is nilpotent for any open ideal $\Fa$. A linearly topologized ring $A$ is 
called \emph{preadmissible} if it possesses a topologically nilpotent open ideal (cf. \cite[Chapitre 0, Definition(7.1.2)]{egai}). By the \emph{completion} of $A$ one 
means the filtered inverse limit $\hat A:=\invlim_\Fa A/\Fa$, where $\Fa$ runs through the partially ordered set of open ideals. For each open ideal $\Fa$ one also calls 
\begin{equation}
\label{admissible02}
\hat\Fa:=\ker(\hat A\rightarrow A/\Fa)
\end{equation}
the completion of $\Fa$. Notice that replacing \eqref{admissible04} by $\{\hat\Fa_i\}_{i\in I}$ 
turns $\hat A$ into a linearly topologized ring and that the natural homomorphism 
\begin{equation}
\label{admissible06}
A\rightarrow\hat A
\end{equation}
is continuous, because the inverse image of \eqref{admissible02} is just $\Fa$. The following is straightforward:

\begin{lem}
\label{LemLinTop}
The map \eqref{admissible06} is injective (resp. surjective) if and only if this holds for the map $\Fa\rightarrow\hat\Fa$, where $\Fa$ is any open ideal of $A$.
\end{lem}

In fact $A$ is called separated (resp. complete) if the map \eqref{admissible06} is injective (resp. surjective) and $A$ is called 
\emph{admissible} if it is separated, complete and preadmissible, again according to \cite[Chapitre 0, Definition(7.1.2)]{egai}. 
A certain sort of preadmissible topology plays a key role in this paper, the following lemma illustrates that:

\begin{lem}
\label{Lem-pA+a-adic-top}
Let $A$ be a ring and let $\Fb$ and $\Fa$ be ideals. 
\begin{itemize}
\item[(i)]
The descending chain of ideals $\{\Fb^i\Fa\}_{i\in\BN_0}$ endows $A$ with the structure of a preadmissible topological ring, satisfying the first axiom of countability.  
\item[(ii)]
The ideal $\Fa$ is $\Fb$-adically separated (resp. complete) as an $A$-module, if and 
only if $A$ is separated (resp. complete) with respect to the topology described in (i).
\item[(iii)]
Suppose that $\Fb=pA+\Fa$, then $\{p^i\Fa+\Fa^{i+1}\}_{i\in\BN_0}$ is also a neighborhood basis for $0_A$ with respect to the topology described in (i).
\end{itemize}
\end{lem}
\begin{proof}
The first statement follows from the inclusion $(\Fb\Fa)^i\subset\Fb^i\Fa$, the second one is a consequence of lemma 
\ref{LemLinTop} and the last one follows from the inclusion $(pA+\Fa)^{2i-1}\Fa\subset p^i\Fa+\Fa^{i+1}\subset (pA+\Fa)^i\Fa$.
\end{proof}

\begin{rem}
\label{admissible07}
All topologies on all rings $A$ that occur in this paper may be assumed to be linear and first-countable, so that there does exist a neighborhood basis for 
$0_A$ consisting of a descending chain $\{\Fa_i\}_{i\in\BN_0}$ as in \eqref{admissible05}. Observe that such a topology is preadmissible if and only if:
\begin{equation}
\label{admissible01}
\exists k\in\BN_0\,\forall m\in\BN_0\,\exists l\in\BN_0:\;\Fa_k^l\subset\Fa_{k+m}
\end{equation}
In fact, some readers might like to think of first-countably topologized preadmissible rings as being topologized by a descending chain 
of ideals $\{\Fb_i\}_{i\in\BN_0}$ with $\Fb_0=A$ and $\Fb_i\Fb_j\subset\Fb_{i+j}$. Indeed let us choose $k_0\in\BN_0$ as indicated in 
\eqref{admissible01} and let us write $\Fb_n$ for the ideal generated by all products $\prod_{i=1}^l\Fa_{n_i-1+k_0}$ where $n$ is the sum of the 
$n_i\in\BN$. Then we do have $\Fb_{(l-1)m+1}\subset\Fa_{m+k_0}\subset\Fb_{m+1}$ whenever $m$ and $l$ are as indicated in \eqref{admissible01}.
\end{rem}
	
\begin{lem}
\label{torsor05}
Let $Y$ be quasi-projective $A$-scheme, where $A$ is a ring.
\begin{itemize}
\item[(i)]
Assume that $A$ is noetherian and $\Fa$-adically separated and complete, where $\Fa$ is an ideal. Then, the canonical map 
\begin{equation}
\label{torsor09}
Y(A)\rightarrow\invlim_nY(A/\Fa^n)
\end{equation}
is bijective.
\item[(ii)]
Assume that $A$ is a linearly topologized complete ring. If $Y$ is smooth, then the canonical map 
$Y(A)\rightarrow Y(A/J)$ is a surjection for all closed ideals $J\subset A$ consisting of topologically nilpotent elements.
\end{itemize}
\end{lem}
\begin{proof}
(i) We begin with the injectivity of \eqref{torsor09}, so let $x,y\in Y(A)$ satisfy $x\equiv y\mod{\Fa^n}$ for every $n$. Thinking of 
$(x,y)$ as a morphism $\Spec A\rightarrow Y\times_AY$ we deduce the existence of an ideal $I\subset A$ such that the diagram

\[\begin{tikzcd}
\Spec A/I\ar[r]\ar[d]&Y\ar[d]\\\Spec A\ar[r]&Y\times_{A}Y\end{tikzcd}
\]

is Cartesian, of course the implication ``quasi-projective $\Rightarrow$ separated'' is relevant in this context. Now $x\equiv y\mod{\Fa^n}$ for every $n$, 
implies $I\subset\bigcap_{n\in\BN}\Fa^n=0$, so that $x=y$. Towards the surjectivity of \eqref{torsor09} we note that $Y$ can be regarded as a locally 
closed subscheme of some $\BP_A^k$, by \cite[Corollaire (5.3.3)]{egaii}. However $\BP^k(A)\isoto\invlim_n\BP^k(A/\Fa^n)$ does hold, due to the moduli 
interpretation of $\BP^k$ given in \cite[Proposition (4.2.3)]{egaii} together with \cite[Lemme (18.3.2.1)]{ega4}, which establishes the projectivity of inverse 
limits of finitely generated projective modules. It is elementary to deduce our assertion for locally closed subschemes of $\BP_A^k$.\\

(ii)  Fix some $y\in Y(A/J)$. According to Grothendieck's ``\'Elimination des hypoth\`eses noeth\'eriennes'' the finitely presented $A$-scheme $Y$ allows a form 
$Y_0$ over a finitely generated $\BZ$-algebra $A_0\subset A$ together with some $y_0\in Y_0(A_0/A_0\cap J)$ inducing $y$ (\cite[Proposition (8.9.1.i)]{ega3}). 
The implications ``smooth $\Rightarrow$ locally finitely presented'' and ``separated $\Rightarrow$ quasi-separated'' are relevant in this context, as the meaning 
of finite presentation consists of being locally finitely presented, quasi-separated and quasi-compact. After enlarging $A_0$ we may assume that $Y_0$ is smooth 
and quasi-projective over $\Spec A_0$, by \cite[Th\'eor\`eme (8.10.5.xiv)]{ega3}. Let us endow $A_0$ with the $J_0$-adic topology, where $J_0:=A_0\cap J$. 
Observe that the natural inclusion $A_0\into A$ is a continuous homomorphism, in view of the finite generation of $J_0$ and the assumption on $J$. Passage 
to completions yields a homomorphism $\hat A_0:=\invlim_nA_0/J_0^n\rightarrow A$. According to the smooth quasi-projectivity of $Y_0$ and part (i) of our 
lemma there exists $\hat y_0\in Y_0(\hat A_0)$ which lifts $y_0$. The image of $\hat y_0$ in $Y_0(A)=Y(A)$ is the requested lift of $y$, because $J$ is closed.
\end{proof}

\subsection{Lifts and inverse limits of torsors}
Occasionally we look at inverse limits of smooth group schemes, for the behavior of whose torsors we offer the following easy lemma:

\begin{lem}
\label{proetale01}
Let $A$ be a ring and let $G\rightarrow\Spec A$ be the inverse limit of a projective system of smooth affine group schemes over $A$
$$ G_1\leftarrow G_2\leftarrow\cdots$$ 
such that all transition homomorphisms are surjective and smooth.		
\begin{itemize}
\item[(i)]
An inverse limit of a projective system of $G_n$-torsors is a $G$-torsor and all $G$-torsors arise in this way. 
\item[(ii)]
If $P$ is a $G$-torsor, there exists a suitable direct limit of rings
$$A\rightarrow R_1\rightarrow R_2\dots$$ with \'etale faithfully flat transition maps, such that $P(\dirlim_nR_n)\neq\emptyset$. 
\end{itemize}
\end{lem}
\begin{proof}
(i): Let $\Spec A\leftarrow P_1\leftarrow P_2\leftarrow\dots$ be a projective system of $G_n$-torsors. First and foremost notice that all torsors in question are affine schemes, by part (i) of corollary \ref{descent04}. In view of this and according to \cite[Proposition (8.2.3)]{ega3} we may consider the inverse limit $P=\invlim_nP_n$, which is again an affine scheme. In fact $P$ is formally smooth and faithfully flat over $A$, as all transition maps $P_{n+1}\rightarrow P_n$ are smooth and surjective, which may be deduced from \cite[Proposition (17.7.1.ii)]{ega4} and \cite[Proposition (2.7.1.iv)]{ega2}, as $G_{n+1}\rightarrow G_n$ is smooth and surjective.\\

In the other direction, consider the projective system of $G_n$-torsors $P_n$ arising from some $G$-torsor $P$ by an extension of structure 
group via $G\rightarrow G_n$. In order to check that $P\rightarrow\invlim_nP_n$ is an isomorphism we can use $G\isoto\invlim_nG_n$ and 
pass to a faithfully flat extension, which is legitimate according to \cite[Proposition (8.2.5)]{ega3} and \cite[Proposition (2.7.1.viii)]{ega2}.\\

(ii): We build up the direct limit together with a compatible system of sections $x_n\in P_n(R_n)$ by induction on $n$. 
Indeed we may apply \cite[Corollaire (17.16.3.ii)]{ega4} to find a lift $x_{n+1}\in P_{n+1}(R_{n+1})$ of $x_n\in P_n(R_n)$ 
for some \'etale faithfully flat $R_n$-algebra $R_{n+1}$. Passing to the limit we get $P(\dirlim_nR_n)\neq\emptyset$.
\end{proof}

\begin{lem}
\label{torsor03}
Fix an affine group scheme $G\rightarrow\Spec\Lambda$, where $\Lambda$ is a discrete valuation 
ring and $G$ is the inverse limit of a projective system of smooth affine group schemes over $\Lambda$
$$G_1\leftarrow G_2\leftarrow\cdots$$
such that all transition homomorphisms are surjective and smooth. Assume that some $\Lambda$-algebra $A$ is complete 
with respect to some linear topology and that $J\subset A$ is a closed ideal consisting of topologically nilpotent elements. 
\begin{itemize}
\item[(i)]
The natural map $\Hom_G(P,Q)\rightarrow\Hom_G(P\times_AA/J,Q\times_AA/J)$ is surjective for any $G$-torsors $P$ and $Q$ over $A$.
\item[(ii)]
Every $G$-torsor over $A/J$ allows a lift to $A$.	
\end{itemize}
\end{lem}
\begin{proof}
If $P$ and $Q$ are as in (i), we let $P_n$ (resp. $Q_n$) be the $G_n$-torsor that arrises from $P$ (resp. $Q$) by extension of the structure group via 
$G\rightarrow G_n$. Let $Y_n=\innHom_{G_n}(P_n,Q_n)$ be the smooth and affine $A$-scheme representing $$R\mapsto\Hom_{G_n}(P_n\times_AR,Q_n\times_AR)$$ 
according to part (ii) of corollary \ref{descent04}. It is easy to see that all transition maps $Y_{n+1}\rightarrow Y_n$ are smooth while the $A$-scheme $Y:=\invlim_nY_n$ represents the functor $R\mapsto\Hom_G(P\times_AR,Q\times_AR)$. Fix a section $y\in Y(A/J)$ and let $y_n$ be its image in $Y_n(A/J)$. By induction we build 
up a sequence $x_n\in Y_n(A)$ such that $x_{n+1}$ maps to $x_n$ while $y_n$ is the $\mod J$ reduction of $x_n$. Indeed consider the commutative square: 
\[
\begin{tikzcd}
\Spec A/J\arrow[r, "y_{n+1}"]\ar[d]&Y_{n+1}\ar[d]\\
\Spec A\arrow[r, "x_n" below]&Y_n
\end{tikzcd}
\]
Since $Y_{n+1}\rightarrow Y_n$ is smooth we may apply part (ii) of lemma \ref{torsor05} to the smooth affine $A$-scheme 
$X:=Y_{n+1}\times_{Y_n,x_n}A$ and thus obtain an oblique morphism in the above diagram, i.e. a $x_{n+1}\in Y_{n+1}(A)$ 
mapping to $x_n$ and having $\mod J$ reduction equal to $y_{n+1}$.\\

(ii): In order to construct lifts of torsors we proceed in three steps. Let us start with verifying the assertion for the special case $G=\GL(h)$. According to 
Hilbert 90 it suffices to prove that any finitely generated projective $A/J$-module $M$, say of constant rank $h$, can be lifted to a projective $A$-module 
of the same rank. Now, it is folklore that the functor associating to a scheme $T$ the set of all projective quotients $\O_T^n\onto\F$ of rank $h$ is representable 
by a smooth projective $\BZ$-scheme, namely the Grassmanian $\Grass(h,n)$. So upon a choice of generators $x_1,\dots,x_n$ of our $A/J$-module $M$ 
we obtain a surjection $(A/J)^n\onto M$, and hence an $A/J$-valued point of $\Grass(h,n)$, to which we can apply part (ii) of lemma \ref{torsor05}.\\

The second step consists of checking the result whenever $G$ is smooth, so that it allows a closed immersion
\begin{equation}
\label{torsor08}
i:G\into\GL(h)_\Lambda
\end{equation}
for some $h\in\BN$, according to \cite[Expos\'e VIB, Proposition 11.11(v)]{sga1} and \cite[Remarque 11.11.1]{sga1}. This allows to mimic the  method of 
\cite[Lemma B.0.17b)]{pappas}, so let us start with a $G$-torsor $P_0$ over $A/J$: By an extension of structure group we get a $\GL(h)_\Lambda$-torsor $P_0*^{G,i}\GL(h)_\Lambda=:Q_0$ over $A/J$. According to the previous step there exists a $\GL(h)_\Lambda$-torsor $Q$ over $A$ which lifts $Q_0$. By lemma \ref{descent06} there exists a smooth and quasi-projective $A$-scheme $Y$ representing the fpqc-quotient $Q/G$. Its functor may be described as $G$-torsors of which the extension of 
structure group via \eqref{torsor08} is isomorphic to $Q$. So our $G$-torsor $P_0$ together with the isomorphism $Q_0=P_0*^{G,i}\GL(h)_\Lambda\cong Q\times_AA/J$ 
constitutes an $A/J$-valued point of $Y$. We are done, because we may apply part (ii) of lemma \ref{torsor05} to this point.\\

The general case can be reduced to this in a third and last step: Choose $G$-torsors $\tilde P_i$ that lift $P*^GG_i$ from $A/J$ to 
$A$. Observe that $\tilde P_i$ and $\tilde P_{i+1}*^{G_{i+1}}G_i$ have like restrictions to $\Spec A/J$, and use part (i) of our lemma 
to choose lifts $\psi_i:\tilde P_i\isoto\tilde P_{i+1}*^{G_{i+1}}G_i$ of the natural isomorphisms between their $\mod J$-reductions. 
The projective system thus defined possesses a limit $\tilde P$, by  part (i) of lemma \ref{proetale01}, and this is the desired lift of $P$.
\end{proof}
	
The following is one of our prime tools:
		
\begin{cor}
\label{torsor06}
Let $\Lambda$ and $G\rightarrow\Spec\Lambda$ be as in lemma \ref{torsor03}. Let $A$ be a $\Lambda$-algebra and assume that a 
descending chain of ideals $\{\Fa_i\}_{i\in\BN_0}$ defines an admissible topology on $A$. Then the natural functor from the category of 
$G$-torsors over $A$ to the category of projective systems of $G$-torsors over $A/\Fa_i$ is an equivalence of categories.
\end{cor}
\begin{proof}
Let $\{P_i\}_{i\in\BN_0}$ be a projective system of $G$-torsors, where $P_i$ is the specimen defined over $A/\Fa_i$. Choose $k_0\in\BN_0$ 
as indicated in \eqref{admissible01}, so that $\Fa_{k_0}$ is a topologically nilpotent open ideal in $A$. Thanks to part (ii) of lemma \ref{torsor03} 
we deduce the existence of a $G$-torsor $P$ over $A$ which lifts $P_{k_0}$. We claim that $P$ serves as an inverse limit of the projective system 
$\{P_i\}_{i\in\BN_0}$: Indeed we have an isomorphism $P_{k_0}\cong P\times_AA/\Fa_{k_0}$. Thanks to part (i) of lemma \ref{torsor03} we can lift 
it to obtain an isomorphism $P_{k_0+1}\cong P\times_AA/\Fa_{k_0+1}$. We keep lifting that and thus obtain a preimage of $\{P_i\}_{i\in\BN_0}$.
\end{proof}

\begin{rem}
\label{admissible03}
In view of the remark \ref{admissible07}, the above is more general than \cite[Lemma B.0.17.b)]{pappas}, which treats reductive group schemes only.
\end{rem}
		
\subsection{Greenberg transforms}
\label{properties}

By means of the usual addition and multiplication one can regard $\BA^1$ as a ring scheme over $\BZ$. Now pick $n\in\BN\cup\{\infty\}$ 
and a prime number $p$ and recall the ring scheme $W_{n}$ of Witt vectors of length $n$ with respect to $p$: Its underlying scheme is 
$\Spec\BZ[x_0,\dots,x_{n-1}]$, on which there exists a unique ring schemes structure, whichs makes the so-called ghost coordinates
$$w_k:\Spec\BZ[x_0,\dots,x_{n-1}]\rightarrow\BA^1;(x_0,\dots,x_{n-1})\mapsto\sum_{i=0}^kp^ix_i^{p^{k-i}}$$
into homomorphisms from the ring scheme $W_n$ to the ring scheme $\BA^1$, for every $k\in\{0,\dots,n-1\}$. If $n=\infty$ we also write
$$W:=W_\infty\cong\invlim_nW_n,$$ 
in which case there is also an additive map $V:W\rightarrow W;(x_0,x_1,\dots)\mapsto(0,x_0,\dots)$, called the Verschiebung, and the 
so-called absolute Frobenius $F:W\rightarrow W$, which is the unique ring endomorphism that satisfies 
$$\forall k\in\BN_{0}:\quad w_k\circ F=w_{k+1}$$ 
Notice that $F(V(x))=px$ and $V(xF(y))=yV(x)$ hold for any two Witt-vectors $x$ and $y$. We denote Cartier's diagonal by
$$\Delta:W\rightarrow W\circ W,$$ 
this is the unique transformation satisfying $w_k\circ\Delta=F^k$, and for its truncation of length $n$ we write $\Delta_n:W\rightarrow W_n\circ W$. At last we 
let $I\subset W$ be the closed ideal scheme defined by the equation $x_0=0$. One also writes $[x]$ for the Teichm\"uller lift, which is the Witt-vector $(x,0,\dots)$.\\
	
We write $B$ for a fixed commutative unital ring, $B-\alg$ for the category of commutative $B$-algebras with unit and let $\Ens$ be the category of sets. 
A covariant functor from $B-\alg$ to $\Ens$ will simply be called a $B$-functor. Notice that every (not necessarily affine) $B$-scheme naturally gives 
rise to a $B$-functor, and that the set of $B$-morphisms between any two $B$-schemes can be recovered from the functorial transformations between their 
$B$-functors. In this optic every $B$-scheme ``is'' (i.e. represents) a $B$-functor. In this subsection we are particularly interested in \emph{arithmetic jet space of $X$ of length $n$}, which is the $B$-functor: 
\begin{equation}
\label{loop01}
M_{n,B}X:B-\alg\rightarrow\Ens:\,R\mapsto X(W_n(R)),
\end{equation}
where $n\in\BN\cup\{\infty\}$ is fixed and $X$ is a $W_n(B)$-scheme. Notice that $W_n(R)$ can be regarded as a 
$W_n(B)$-algebra, so that the notation $X(W_n(R))$ is meaningful. Also, note that the construction of $M_{n,B}$ is compatible 
with base change, more precisely, for any $B$-algebra $R$, and any $W_{n}(B)$-scheme $X$ there is a canonical isomorphism 
\begin{align}
\label{MnBBaseChange}
\big(M_{n,B}X\big)_{R}\cong M_{n,R}(X_{W_{n}(R)})
\end{align}
Observe that for every $k\in\{0,\dots,n-1\}$ there is a transformation of $B$-functors
\begin{equation}
\label{loop02}
M_{n,B}X\stackrel{w_k}{\rightarrow}X\times_{W_n(B),w_k}B
\end{equation}
that is canonically induced by the ghost maps $w_k:W_n(R)\rightarrow R$. It seems to be 
well-known (e.g. see \cite[Proposition 29]{kreidl} and references therein), that the functor 
\begin{align*}
B-\alg&\rightarrow W_n(B)-\alg\\
R&\mapsto W_n(R)
\end{align*}
allows a left-adjoint, say
$$\Lambda_{n,B}:W_n(B)-\alg\rightarrow B-\alg,$$
so that 
\begin{align}
\label{LeftRightAdjMLambda}
M_{n,B}\Spec A\cong\Spec\Lambda_{n,B}(A)
\end{align} 
where $A$ is an arbitrary $W_n(B)$-algebra. Please see \cite{borger} for much more on this issue, 
e.g. for generalizations to algebraic spaces. In particular, the affine case of \eqref{loop02} reads:
$$m_k^\sharp:A\otimes_{w_k,W_n(B)}B\rightarrow\Lambda_{n,B}(A)$$ 
We refer to these as the co-ghost maps (and may occasionally confuse them with the $W_n(B)$-linear pendants 
$m_k:A\rightarrow\Lambda_{n,B}(A)_{[w_k]}$ for any $k<n$). Before we move on, let us list some properties of the functor $\Lambda_{n,B}$

\begin{lem}
\label{LemPropLambda}
Let $B$ be a ring. We have the following facts about the functor $\Lambda_{n,B}$:

\begin{enumerate}
\item[(i)]
it preserves arbitrary direct limits,
\item[(ii)]
it turns surjections into surjections
\item[(iii)]
there is a canonical isomorphism $$\Lambda_{n,B}(W_n(B)[t_{1},\dots,t_{r}]) \cong B[t_{ij};1\leq i\leq r,1\leq j\leq n]$$ where $t_{i}$ and $t_{ij}$ are independent variables.
\item[(iv)]
it turns finitely generated (resp. finitely presented) $W_n(B)$-algebras into finitely generated (resp. finitely presented) $B$-algebras, provided that $n$ is finite.
\item[(v)]
it turns formally smooth $W_n(B)$-algebras into formally smooth $B$-algebras.
\item[(vi)]
more generally, there is a canonical and functorial isomorphism  
\begin{align}
\label{LambdaOmegaId}
\Omega_{\Lambda_{n,B}(A)/B}^1\cong\prod_{k=0}^{n-1}\big(\Lambda_{n,B}(A)\otimes_{m_k,A}\Omega_{A/W_n(B)}^1\big)
\end{align}
where $A$ is a $W_{n}(B)$-algebra
\end{enumerate}
\end{lem}
\begin{proof}
The first two assertions follow from the fact that $\Lambda_{n,B}$ has a right adjoint. Statement (iii) follow from the adjunction formula \eqref{LeftRightAdjMLambda}, 
and the canonical set-theoretical bijection $W_{n}(R)\cong R^{\times n}$. Statement (iv) follows from (ii) and (iii). Assertion (v) is a consequence of the formula 
\begin{equation}
\label{NO}
W_n(C\oplus M)\cong W_n(C)\oplus\prod_{k=0}^{n-1}M_{[w_k]}
\end{equation} 
where the multiplication on $C\oplus M$ is $(c'+m')(c+m)=c'c+c'm+cm'$ and $M$ is any $C$-module, the dual number 
case of this appeared for the first time in \cite{norman}. Finally, the last assertion would follow by a similar argument.
\end{proof}

Now suppose that $X$ is an affine $W(B)$-scheme, so that the $B$-scheme $M_BX:=M_{\infty,B}X$ is well-defined. It will not 
cause confusion to write $M_{n,B}X$ for 
$M_{n,B}(X\times_{W(B)}W_n(B))$ (for some $n\in\BN$). Moreover, when using this convention one has $M_BX\cong\invlim_nM_{n,B}X$. We need the following 
easy variant of \cite[Proposition 2.2.1c),d)]{pappas}:

\begin{lem}
\label{answer04}
Fix $0\leq n'<n\leq\infty$ and a smooth affine $W_n(B)$-scheme $X$. 
\begin{itemize}
\item[(i)]
The truncation map
\begin{equation}
\label{loop03}
M_{n,B}X\rightarrow M_{n',B}X
\end{equation}
is smooth unless $n=\infty$, in which case it is still flat and formally smooth.
\item[(ii)]
The morphism \eqref{loop03} is faithfully flat if and only if $X\times_{W_n(B),w_0}B\rightarrow\Spec B$ is surjective or $n'>0$.
\end{itemize}
In particular the zeroth ghost map $M_{n,B}X\rightarrow X\times_{W_n(B),w_0}B$ (cf. \eqref{loop02}) is faithfully flat and formally smooth.
\end{lem}
\begin{proof}
Without loss of generality we assume that $n$ is finite, as the case $n=\infty$ follows by a limit process, regardless of whether we are discussing the 
assertion (i) or (ii). Observe also that \cite{GreenbergII} or \cite[Table 3]{borger} settle the smoothness of each $M_{n,B}X\rightarrow\Spec B$, which 
is the case $n'=0$. To deduce the smoothness of \eqref{loop03}. We observe that $\O_{M_{n,B}X}\otimes_{\O_{M_{n',B}X}}\Omega_{M_{n',B}X/B}^1$ 
is a direct factor of $\Omega_{M_{n,B}X/B}^1$, which one can read off from the formula \eqref{LambdaOmegaId}. We deduce the smoothness of 
\eqref{loop03} according to \cite[Theoreme 17.11.1]{ega4}, more specifically the implication ``c)$=>$b)'' therein.\\

Towards the surjectivity of \eqref{loop03} we let $x$ be a point of $M_{n',B}$ over an algebraically closed field $k$, i.e. $x\in X(W_{n'}(k))$. If the characteristic 
of $k$ is $p$ we know that $W_n(k)\twoheadrightarrow W_{n'}(k)$ is a nilpotent thickening, provided of course, that $n'>0$! So in this case the existence of a lift of 
$x$ to $M_{n,B}X$ follows easily from formal smoothness of $X$. If the characteristic of $k$ is different from $p$, we have $W_n(k)\cong W_{n'}(k)\times k^{n-n'}$. 
It follows that there exists a $W_n(k)$-linear section $W_{n'}(k)\rightarrow W_n(k)$, provided of course, that $n'>0$! Again we are able to lift $x$ to $M_{n,B}X$. 
At last note that the scheme $M_{1,B}X$ is the base change of $X$ along $W_n(B)\to B$ and $M_{0,B}X\cong\Spec B$, thus settling the fringe case $n'=0$.
\end{proof}

\begin{lem}
\label{torsor01}
Let $B$ be a ring and $n\in\BN\cup\{\infty\}$. Let $G$ be a smooth affine group scheme over $W_n(B)$. 
\begin{itemize}
\item[(i)]
If $P$ is a $G$-torsor, then $M_{n,B}P$ is a $M_{n,B}G$-torsor. 
\item[(ii)]
If $Q$ is another $G$-torsor, then the natural map from $M_{n,B}\innHom_G(P,Q)$ to 
$\innHom_{M_{n,B}G}(M_{n,B}P,M_{n,B}Q)$ is an isomorphism of $B$-schemes.
\end{itemize}
\end{lem}
\begin{proof}
(i): Observe that $M_{n,B}P$ is a well-defined affine $B$-scheme because $P$ is affine as a consequence of part (i) of 
corollary \ref{descent04}. It is clear that the natural action of $M_{n,B}G$ on $M_{n,B}P$ from the right gives the latter 
the structure of a pseudo-torsor, as $M_{n,B}$ preserves fiber products (cf. \cite[Proposition 2.2.1a)]{pappas}). Since 
$P\rightarrow\Spec W_n(B)$ is surjective we are able to deduce the faithful flatness of $M_{n,B}P\rightarrow M_{0,B}X\cong\Spec B$, 
in view of part (ii) of lemma \ref{answer04}. So that we are done according to lemma \ref{descent05}.\\

(ii): According to \cite[Proposition (2.7.1.viii)]{ega2} we can pass to a faithfully flat extension, and use \eqref{MnBBaseChange}. 
\end{proof}
		
Now fix a finite field $k_0$ of characteristic $p$ and let $G$ be an affine group scheme over $W(k_0)$. Let $R$ be a $W(k_0)$-algebra. Using the Cartier map 
$W(k_0)\to W(W(k_0))$, we can endow the ring $W_{n}(R)$ with the structure of a $W(k_0)$-algebra. We write $L_n^+G=M_{n,W(k_0)}(G_{W_n(W(k_0))})$, 
resp. $L^+G=M_{W(k_0)}(G_{W(W(k_0))})$, for the affine group schemes over $W(k_0)$ whose group of $R$-valued points is given by $G(W_n(R))$, resp. 
$G(W(R))$. Notice that there are homomorphisms of group schemes
\begin{eqnarray}
\label{frob01}
&&L^+G\stackrel{F^G}{\rightarrow}{^FL^+G}=L^+{^FG}\\
\label{ghost01}
&&w_k^G:L^+G\rightarrow{^{F^k}}G,
\end{eqnarray}
which are canonically induced by the absolute Frobenius $F:W(R)\rightarrow W(R)$ and the ghost maps:
$$w_k:W(R)\rightarrow R_{[F^k]};(x_0,x_1,\dots)\mapsto\sum_{i=0}^kp^ix_i^{p^{k-i}}$$
Furthermore, we write $L^{>0}G\subset L^+G$ for the kernel of $w_0^G:L^+G\rightarrow G$. Also notice that $M_{n,k_0}(G_{W_n(k_0)})$, resp. 
$M_{k_0}(G)$ are the special fibers of $L_n^+G$, resp. $L^+G$ (see \eqref{MnBBaseChange}). One more piece of notation will prove useful in the 
sequel: If $\Fa$ is a commutative but not necessarily unital $W(k_0)$-algebra, so that $W(k_0)\oplus\Fa$ is a unital $W(k_0)$-algebra with augmentation 
ideal $\Fa$, we write $G(\Fa)$ for the kernel of the map from $G(W(k_0)\oplus\Fa)$ to $G(W(k_0))$. Observe that we have an exact sequence
\begin{equation*}
1\rightarrow G(\Fb)\rightarrow G(\Fa)\rightarrow G(\Fa/\Fb),
\end{equation*}
whenever $\Fb\subset\Fa$ is a $W(k_0)$-submodule with $\Fb\Fa\subset\Fb$, which are the so-called ideals of $\Fa$. In the same vein we have an exact sequence 
\begin{equation*}
1\rightarrow G(\Fa)\rightarrow G(R)\rightarrow G(R/\Fa),
\end{equation*}
whenever $\Fa$ happens to be an ideal in a $W(k_0)$-algebra $R$. In this optic one can think of $L^{>0}G(R)$ as the group of $I(R)$-valued points of $G$, 
where $I(R)$ is regarded as a non-unital $W(k_0)$-algebra. Recall that $G$ is called a vector group scheme provided that it can be written in the form
\begin{equation}
\label{frob04}
\Spec\Big(\bigoplus_{l=0}^\infty\Sym_{W(k_0)}^l\check\Fv\Big)=:\ul\Fv
\end{equation}
for some free $W(k_0)$-module $\Fv$ of finite rank. In this case the map \eqref{frob01} has an important variant, 
observe that $\ul\Fv$, $L^+\ul\Fv$ and $L^{>0}\ul\Fv$ represent respectively the group functors on $W(k_{0})$-algebras:
$$R\mapsto\begin{cases}
\Fv(R)=R\otimes_{W(k_0)}\Fv\\
L^+\ul\Fv(R)=W(R)\otimes_{W(k_0)}\Fv\\
L^{>0}\ul\Fv(R)=I(R)\otimes_{W(k_0)}\Fv,
\end{cases}$$
allowing us to define the mutually inverse maps:
\begin{eqnarray}
\label{frob03}
&&V_\Fv^{-1}:L^{>0}\ul\Fv\rightarrow{^FL^+\ul\Fv};\,a\otimes
X\mapsto V^{-1}(a)\otimes X\\
\label{frob02}
&&V_\Fv:{^FL^+\ul\Fv}\rightarrow L^{>0}\ul\Fv;\,a\otimes X\mapsto V(a)\otimes X,
\end{eqnarray}
where
$${^FL^+\ul\Fv}= L^+\ul\Fv\times_{W(k_0),F}W(k_0)\cong L^+\ul{W(k_0)\otimes_{F,W(k_0)}\Fv}.$$
	
\begin{lem}
\label{torsor02}
Let $B$ be a $p$-adically separated and complete $W(k_0)$-algebra and $n\in\BN\cup\{\infty\}$. Let $G$ be a smooth affine group scheme over $W(k_0)$. The functor:
$$P\mapsto L_{n,B}^+P$$ 
provides an equivalence from the category of $G$-torsors over $W_n(B)$ to the category of $L_n^+G$-torsors over $B$.
\end{lem}
\begin{proof}
Corollary \ref{torsor06} and part (i) of lemma \ref{proetale01} are reducing the problem to the 
discrete finite-length case, i.e. we are allowed to assume $pB\subset\sqrt{0_B}$ and $n\neq\infty$.\\

The full faithfulness of our functor is clear form part (ii) of lemma \ref{torsor01}. It remains to show the essential 
surjectivity. Let $Q$ be a $L_n^+G$-torsor over $B$ and let $\tilde B$ be an \'etale faithfully flat extension that trivializes 
$Q$. Observe that $W_n(B)\rightarrow W_n(\tilde B)$ is \'etale faithfully flat too and that we have canonical isomorphisms 
$$G(W_n(\tilde B)^{\otimes_{W_n(B)}m})\cong L_n^+G(\tilde B^{\otimes_Bm})$$
by \cite[Proposition A.8]{langer1} and \cite[Proposition A.12]{langer1}. Hence the result follows from descent theory. 
\end{proof}
	
\begin{rem}
The case $G=\GL(h)_{W(k_0)}$ is nothing but \cite[Corollary 34]{zink2} in disguise. Again, please see \cite[Proposition B.0.15b)]{pappas} for the reductive case.
\end{rem}
	
\subsection{Basics on algebraic groups and their cocharacters}
\label{simplerversion}
Let $k_0$ be a finite field. For any smooth affine group scheme $G$ over $W(k_0)$ we let $\Rep_{W(k_0)}(G)$ be the category 
of representations of $G$ on free $W(k_0)$-modules of finite rank. Furthermore we write $\omega^G$ for the natural fiber 
functor from $\Rep_{W(k_0)}(G)$ to the category $W(k_0)-\ect$ of free $W(k_0)$-modules of finite rank. Let us fix a cocharacter 
$\mu:\BG_{m,W(k_0)}\rightarrow G$ and an integer $m$. For a representation $\rho\in\Ob_{\Rep_{W(k_0)}(G)}$ we write
$$\omega_\mu^G(m,\rho)\subset\omega^G(\rho)$$
for the subspace of $\mu$-weight $m$ and
$${\rm Fil}_\mu^m\rho=\bigoplus_{l=m}^\infty\omega_\mu^G(l,\rho)$$
for the filtration which is induced by $\mu$. For every non-negative $m\in\BN_0$ we consider the associated subgroup functor of $G$:
\begin{equation}
\label{filtrategroup}
R\mapsto\{g\in G(R)\,|\,\forall\rho,l_0:\\
(\rho(g)-1)(\omega_\mu^G(l_0,\rho))\subset R\otimes_{W(k_0)}{\rm Fil}_\mu^{l_0+m}\rho\}
\end{equation}
We note the following facts:
\begin{itemize}
\item[(i)]
According to \cite[Chapitre IV, 2.1.3]{rivano} there exists a closed subgroup scheme $U_\mu^m\into G$ which represents the subgroup 
functor \eqref{filtrategroup}, moreover $U_\mu^m$ is connected, unipotent and smooth over $W(k_0)$, and its Lie algebra can be computed from
\begin{equation*}
\Lie U_\mu^m\cong{\rm Fil}_\mu^m\Ad^G,
\end{equation*}
by \cite[Chapitre IV, Proposition 2.1.4.1]{rivano}. Please see \cite[2.1]{gabber} or \cite[4th part of Theorem 4.1.7]{conradnote} for alternative treatments.
\item[(ii)]
If $g\in U_\mu^m(R)$ and $h\in U_\mu^n(R)$ then $ghg^{-1}h^{-1}\in U_\mu^{m+n}(R)$, which can be deduced from:
$$1-\rho(ghg^{-1}h^{-1})=\rho(g)[1-\rho(g)^{-1},1-\rho(h)]\rho(h)^{-1}$$
where $\rho\in\Ob_{\Rep_{W(k_0)}(G)}$ is a representations of $G$ on some free $W(k_0)$-module of finite rank. 
\end{itemize}
	
We will use the notations:
\begin{eqnarray*}
&&P_\mu:=U_\mu^0\\
&&U_\mu:=U_\mu^1\\
&&\Fp_\mu:=\Lie P_\mu\cong{\rm Fil}_\mu^0\Ad^G\\
&&\Fu_\mu:=\Lie U_\mu\cong{\rm Fil}_\mu^1\Ad^G\\
&&\Fg:=\Lie G=\omega^G(\Ad^G)
\end{eqnarray*}
	
The following lemma can easily be checked by means of an embedding into a suitable general linear group (cf. \cite[A.0.9]{pappas}):
	
\begin{lem}
\label{IntAction}
There exists an action
\begin{equation}
\label{movie}
\intaut_\mu:\BA_{W(k_0)}^1\times_{W(k_0)}P_\mu\rightarrow P_\mu
\end{equation}
of the monoid $\BA_{W(k_0)}^1$ on the group scheme $P_\mu$ by group scheme endomorphisms, whose restriction to $\BG_{m,W(k_0)}$ is given by the formula:
$$\intaut_\mu(z,g)=\mu(z)g\mu(\frac1z)$$
\end{lem}

Proofs of the following result can be found in \cite[4th part of Theorem 4.1.1]{conradnote} and \cite[3rd part of Proposition 2.1.8]{gabber}:
	
\begin{lem}
\label{facta}
Let $k_0$ and $G$ be as above. Let $\mu:\BG_{m,W(k_0)}\rightarrow G$ be a cocharacter. The group multiplication
$m_G: P_\mu\times_{W(k_0)}U_{\mu^{-1}}\rightarrow G$ identifies the product $P_\mu\times_{W(k_0)}U_{\mu^{-1}}$ with an open subscheme $G^*\subset G$
\end{lem}
	
\begin{defn}
Let $k_0$, $G$ and $\mu$ be as above. Let $A$ be a $W(k_0)$-algebra and $I$ an ideal in $A$. We denote 
by $\Gamma_\mu^{A,I}$ the preimage of $P_\mu(A/I)$ in $G(A)$ under the projection $G(A)\to G(A/I)$. We call $\Gamma_\mu^{A,I}$ the \emph{parahoric subgroup} of $G$ attached to $A, I$ and $\mu$.
\end{defn}

We have the following remarkable corollary:
	
\begin{cor}
\label{productIII}
Let $k_0$ and $G$ be as above. Let $\mu:\BG_{m,W(k_0)}\rightarrow G$ be a cocharacter and let $A$ be a $W(k_0)$-algebra. 
Consider an ideal $I$ which is contained in the Jacobson radical $\rad(A)$ of $A$. Then, the multiplication in $G(A)$ gives identifications
\begin{eqnarray*}
&&P_\mu(A)\times U_{\mu^{-1}}(I)\cong\Gamma_\mu^{A,I}\\
&&P_\mu(I)\times U_{\mu^{-1}}(I)\cong G(I)
\end{eqnarray*}
\end{cor}
\begin{proof}
The second identification follows from the first one. The idea of 
the proof is as in \cite[proposition 3.1.3]{pappas} and for convenience of the reader we spell out the details: Consider $h\in\Gamma_\mu^{A,I}\subset G(A)$. We would like to show that the corresponding $h: \Spec(A)\rightarrow G$ factors through $G^*$. 
Since $G^*$ is an open subscheme of $G$ by the previous lemma, 
it suffices to show that the set theoretic image of $h$ is 
contained in $G^*$, which is equivalent to the emptiness of 
$h^{-1}(G\backslash G^*)$. Notice that every non-empty closed 
subset of $\Spec A$ meets $\Spec A/I$ in at least one point, given that every maximal ideal of $A$ contains $I$. Whence it suffices to check the disjointness of $h^{-1}(G\backslash G^*)$ and $\Spec A/I$. However, from $h\in\Gamma_\mu^{A,I}$ we deduce 
$h(\Spec A/I)\subset P_\mu\subset G^*$.
\end{proof}
	
Recall the Tannakian formalism involving the forgetful fiber functor $\omega^G$ on the representation category 
$\Rep_{W(k_0)}(G)$ of our smooth affine group $G$, which is flat over the Dedekind ring $W(k_0)$: For every $W(k_0)$-algebra 
$R$ one recovers $G(R)$ as the group of $\otimes$-preserving functorial transformations from the $R-\ect$-valued 
$\otimes$-functor $\rho\mapsto R\otimes_{W(k_0)}\omega^G(\rho)$ to itself, i.e.: 
$$G(R)=\Aut_R^\otimes(R\otimes_{W(k_0)}\omega^G)=\End_R^\otimes(R\otimes_{W(k_0)}\omega^G)\subset\End_R(R\otimes_{W(k_0)}\omega^G)$$ 
We note a useful Lie-theoretic analog: One recovers $R\otimes_{W(k_0)}\Lie G$ as the 
additive group of $R$-linear transformations $D:R\otimes_{W(k_0)}\omega^G\rightarrow R\otimes_{W(k_0)}\omega^G$ such that
$$D(\rho_1\otimes\rho_2)=D(\rho_1)\otimes_{W(k_0)}\id_{\omega^G(\rho_2)}+\id_{\omega^G(\rho_1)}\otimes_{W(k_0)}D(\rho_2)$$
holds for any two objects $\rho_1$ and $\rho_2$ of $\Rep_{W(k_0)}(G)$. In the sequel we refer to this as the group $\Der_R^\otimes(R\otimes_{W(k_0)}\omega^G)$ of 
$R$-linear $\otimes$-derivations of $R\otimes_{W(k_0)}\omega^G$. Finally, observe that both $G(R)$ and $R\otimes_{W(k_0)}\Lie G$ can be viewed as subfunctors of: 
$$R\mapsto\End_R(R\otimes_{W(k_0)}\omega^G)=
R\otimes_{W(k_0)}\End_{W(k_0)}(\omega^G)$$
	
For an arbitrary connected unipotent algebraic group $U$ 
over a field $K$ of characteristic $0$ there exists a 
canonical isomorphism of $K$-varieties
$$\exp_U:\ul{\Lie U}\isoto U,$$ 
which can be described as follows: Let $R$ be a $K$-algebra 
and consider a $R$-valued point of $\ul{\Lie U}$, i.e. some 
$D\in R\otimes_K\Lie U$, which can be viewed as a $\otimes$-derivation of $R\otimes_K\omega^U$ in the above sense (note 
that since the characteristic of $K$ is zero, $U$ is smooth). 
Then $\exp_U(D)$ is the $\otimes$-automorphism of 
$R\otimes_K\omega^U$ given by the classical series
$$\exp_U(D)(\rho)=\sum_{n=0}^{\infty}\frac{D(\rho)^n}{n!}$$
Note that $D(\rho)^n$ vanishes on an $n$-dimensional representation $\rho$ of $U$, and that all terms on the right are elements 
in the ring $R\otimes_K\omega^U(\check\rho\otimes_K\rho)$, see \cite[Chapter IV.2, Proposition 4.1]{tome I} for more details.

\begin{lem}
\label{factb}
As before we fix a $W(k_0)$-rational cocharacter $\mu$ of some 
smooth affine $W(k_0)$-group $G$, so that the generic fiber of 
$U_\mu$ is a connected unipotent algebraic group over the fraction field $K:=W(k_0)[\frac1p]$. If the $\mu$-weights of $\Fg$ are strictly less than $p$, then $\exp_{{U_\mu}_K}$ extends uniquely 
to an integral morphism $e:\ul{\Fu_\mu}\rightarrow U_\mu$.
\end{lem}
\begin{proof}
We have to show that for all $n\geq 1$, all 
$W(k_{0})$-algebras $R$, all representations $\rho$ and all 
$D\in R\otimes_{W(k_{0})}\Lie U_{\mu}$, $n!$ divides $D(\rho)^{n}$.  Thanks to the formula
$$\frac{\sum_{j=1}^ka_j\otimes D_j(\rho)^n}{n!}=\sum_{\sum_{j=1}^ki_j=n}\prod_{j=1}^ka_j^{i_j}\otimes\prod_{j=1}^k\frac{D_j(\rho)^{i_j}}{i_j!}$$ 
it is enough to check that $n!$ is a divisor of $D(\rho)^n$ for 
every $D\in\Lie U_\mu$. So we prove $n!\mid D(\rho)^n$. Assume 
that the statement is true when $n$ is a power of $p$. If $n$ is 
not a power of $p$, then we can write $n$ is a sum $\sum_lm_lp^l$ with $0\leq m\leq p-1$. Since the $p$-adic valuations of $(\sum_lm_lp^l)!$ and $\prod_l(p^l)!^{m_l}$ are equal, $n!$ 
divides $D(\rho)^{n}$ if and only if $\prod_l(p^l)!^{m_l}$ 
does. We have
\[\frac{D(\rho)^{n}}{\prod_l(p^l)!^{m_l}}=\prod\frac{D(\rho)^{m_{l}p^{l}}}{(p^{l})!^{m_{l}}}=\prod\big(\frac{D(\rho)^{p^{l}}}{(p^{l})!}\big)^{m_{l}}\]
By the assumption, the quotient $\frac{D(\rho)^{p^{l}}}{(p^{l})!}$ 
is integral and thus so is $\frac{D(\rho)^{n}}{n!}$. Therefore, 
it remains to show the statement when $n$ is a power of $p$, say $p^{l}$. We prove by induction on $l$. When $l=0$ there is nothing 
to prove. We have
\[\frac{D^{p^{l}}}{(p^{l})!}=\frac{(p^{l-1}!)^{p}p}{(p^{l})!}\cdot\frac{1}{p}\cdot\big(\frac{D^{p^{l-1}}}{p^{l-1}!}\big)^{p}\]
Note that $\frac{(p^{l-1}!)^{p}p}{(p^{l})!}$ is a unit of $\BZ_{p}$, and so, since by assumption $\frac{D^{p^{l-1}}}{p^{l-1}!}$ 
is integral, we have to show that 
$\big(\frac{D^{p^{l-1}}}{p^{l-1}!}\big)^{p}$ is divisible by $p$. Consider the elements
\begin{eqnarray}
&&\sigma=\sum_{j=1}^{p^l}\epsilon_j\\
&&\pi=\frac{(p^l)!}p\prod_{j=1}^{p^l}\epsilon_j
\end{eqnarray}
in the ring $R=W(k_0)[\epsilon_1,\dots\epsilon_{p^l}]/(\epsilon_1^2,\dots,\epsilon_{p^l}^2)$. Let $R'$ be the sub-$W(k_{0})$-algebra 
generated by these two elements. In $R'$ the ideal $(p\pi-\sigma^{p^l},\pi^2,\sigma\pi)$ is trivial. Notice that $\sigma^{p^l+1}$ vanishes 
and that $\{\frac{\sigma^k}{k!}\mid\,k=1,\dots,p^l\}$ are given by the elementary symmetric polynomials in $\epsilon_1,\dots,\epsilon_{p^l}$. It follows 
$\BQ\otimes W(k_0)[\sigma]/(\sigma^{p^l+1})\cong\BQ\otimes R'$, so that the $W(k_0)$-rank of the free $W(k_0)$-module $R'$ is equal to $p^l+1$. We deduce 
that the canonical projection $W(k_0)[\sigma,\pi]/(p\pi-\sigma^{p^l},\pi^2,\sigma\pi)\twoheadrightarrow R'$ must be isomorphism, as both sides have the same 
$W(k_0)$-rank. Next consider the element
$$g:=\prod_{i=1}^{p^l}(1+\epsilon_iD)\in U_\mu(R),$$
and observe that the formula
$$\frac{(p^{l-1}!)^pp}{p^l!}\cdot\pi\cdot(\frac{D^{p^{l-1}}}{p^{l-1}!})^p+\sum_{k=0}^{p^l-1}\sigma^k\frac{D^k}{k!}=g$$
holds in $\End_R(R\otimes_{W(k_0)}\omega^{U_\mu})$. Consider the ring homomorphism 
$R'\rightarrow k_0[\epsilon]/\epsilon^2$, defined by $\sigma\mapsto0$ and $\pi\mapsto\epsilon$. 
We deduce that the $\pmod p$-reduction of $(\frac{D^{p^{l-1}}}{p^{l-1}!})^p$ is a $k_0$-linear $\otimes$-derivation of $k_0\otimes_{W(k_0)}\omega^{U_\mu}$, 
but due to $\mu$-weight considerations this can only be the trivial derivation, so that we have $p\mid(\frac{D^{p^{l-1}}}{p^{l-1}!})^p$.
\end{proof}

\begin{cor}
\label{factc}
As before we fix a $W(k_0)$-rational cocharacter $\mu$ of some 
smooth affine $W(k_0)$-group $G$, so that the generic fiber of 
$U_\mu$ is a connected unipotent algebraic group over the fraction field $K:=W(k_0)[\frac1p]$. If the $\mu$-weights of $\Fg$ are 
bounded by $1$, then $\exp_{{U_\mu}_K}$ extends uniquely to an integral isomorphism of group schemes $e:\ul{\Fu_\mu}\isoto U_\mu$.
\end{cor}
\begin{proof}
In view of $[\Fu_\mu,\Fu_\mu]\subset{\rm Fil}_\mu^2\Ad^G=0$, the Baker-Campbell-Hausdorff formula implies that $\exp_{{U_\mu}_K}$ 
is a homomorphisms of group schemes, hence 
$e:\ul{\Fu_\mu}\rightarrow U_\mu$ is one too. Since 
$\exp_{{U_\mu}_K}$ induces the identity on the tangent spaces at 
the neutral section, the same holds for $e$. In view
\cite[Th\'eor\`eme (17.11.1)]{ega4} we may deduce that $e$ is 
\'etale, so that its kernel is an \'etale group scheme. In fact, 
we obtain the triviality of $\ker(e)$ from looking at the generic fiber. So $e$ is an \'etale monomorphism, i.e. an open immersion, 
due to \cite[Th\'eor\`eme (17.9.1)]{ega4}. The connectedness of $U_\mu$ finishes the proof.
\end{proof}

\begin{rem}
The above is more general than \cite[Lemma 6.3.2]{lau2} and \cite[A.0.13]{pappas}, which treat reductive group schemes only.
\end{rem}

\begin{rem}
In the more general situation of lemma \ref{factb} the morphism $e$ fails to respect the group law, 
but it remains to be an isomorphism of schemes, as can be seen as follows: With a little bit of extra work, 
along the same lines of the proof of corollary \ref{factc}, one can show that $e$ induces isomorphisms 
$$\ul{\omega_\mu^G(m,\Ad^G)}\isoto U_\mu^m/U_\mu^{m+1}$$ 
(for all $m$). Moreover, observe that $U_\mu/U_\mu^{m+1}$ is a $U_\mu^m/U_\mu^{m+1}$-torsor 
over $U_\mu/U_\mu^m$. Arguing inductively, one sees that the induced maps of quotients 
$$\ul{\bigoplus_{l=1}^m\omega_\mu^G(l,\Ad^G)}\rightarrow U_\mu/U_\mu^{m+1}$$ 
are isomorphisms, because they can be viewed as homomorphisms between\\ $\ul{\omega_\mu^G(m,\Ad^G)}$-torsors over $U_\mu/U_\mu^m$.
\end{rem}

\begin{cor}
\label{factd}
Let the assumptions on $k_{0}$, $G/W(k_{0})$ and $\mu:\BG_{m,W(k_{0})}\rightarrow G$ be as in corollary \ref{factc}. Then 
$$\Gamma_{\mu^{-1}}^{A,pA}\subset\mu(p)G(A)\mu(\frac1p)$$ 
holds for any torsion-free $W(k_0)$-algebra $A$.
\end{cor}
\begin{proof}
Without loss of generality we can assume $pA\subset\rad A$, as $(1+pA)^{-1}A\times A_p$ is faithfully flat over $A$. 
By corollary \ref{productIII} we can consider the groups 
$P_{\mu^{-1}}(A)$ and $U_\mu(pA)$ separately. However,
$$\mu(\frac1p)P_{\mu^{-1}}(A)\mu(p)\subset P_{\mu^{-1}}(A)\subset G(A)$$
follows immediately from lemma \ref{IntAction}, and
$$\mu(\frac1p)U_\mu(pA)\mu(p)\subset U_\mu(A)\subset G(A)$$ 
follows from corollary \ref{factc}. Observe that for $X\in A\otimes_{W(k_0)}\Lie U_\mu$ we have: 
$$\mu(\frac1p)\exp_{U_\mu}(X)\mu(p)=\exp_{U_\mu}(\Ad^G(\mu(\frac1p))X)$$
\end{proof}

%\subsection{Definition of $H_\mu$ and $\Phi_\mu$}
\subsection{On the display group $H_\mu$ 
and the conjugated Frobenius $\Phi_\mu$}
\label{productII}

\begin{defn}
    The \emph{display group} attached to $\mu:\BG_{m,W(k_{0})}\rightarrow G$ (as in corollary 
    \ref{factc}) is the subgroup scheme $H_\mu:=L^+G\times_{w_0^G,G}P_\mu$ of $L^+G$. This is a flat 
    and formally smooth $W(k_0)$-group scheme characterized by $H_\mu(R)=\Gamma_\mu^{W(R),I(R)}$, for any $W(k_0)$-algebra $R$. Please note that we have the closed immersions:
    $$L^{>0}G\into H_\mu\into L^+G$$
\end{defn}

Consider $H_\mu=L^+G\times_{w_0^G,G}P_\mu$, where $w_0^G$ is the zeroth ghost component (cf. \eqref{ghost01}). This is a flat and 
formally smooth $W(k_0)$-group scheme characterized by $H_\mu(R)=\Gamma_\mu^{W(R),I(R)}$, please note that we have the closed immersions:
$$L^{>0}G\into H_\mu\into L^+G$$

\begin{prop}
\label{productI}
Assume that $pR$ is contained in the Jacobson radical $\rad(R)$ of $R$. Then multiplication in $L^+G(R) $ gives an identification
$$L^+P_\mu(R)\times L^{>0}U_{\mu^{-1}}(R)\cong H_\mu(R).$$
\end{prop}
\begin{proof}
By part (i) of \cite[Lemma 2.1.1]{pappas} it follows that $I(R)\subset\rad(W(R))$. The result now follows from the 
definition of $H_\mu$ and corollary \ref{productIII}.
\end{proof}

\begin{prop}
\label{coFrob}
There is a group scheme homomorphism
$$\Phi_\mu: H_\mu\rightarrow{^FL^+G}=L^+{^FG}$$
characterized by the following property: If $g$ is an $R$-valued point of $H_\mu$, where $R$ is a $W(k_0)[\frac1p]$-algebra, then
\begin{equation}
\label{characterization}
\Phi_\mu(g)={^F(\mu(p)g\mu(\frac1p))}\in{^FG(W(R))}.
\end{equation}
\end{prop}
\begin{proof}
We begin with the introduction of maps
\begin{eqnarray}
\label{nonpositive}
&\Phi_0:L^+P_\mu\rightarrow{^FL^+P_\mu}&\\
\label{positive}
&\Phi_1:L^{>0}U_{\mu^{-1}}\rightarrow{^FL^+U_{\mu^{-1}}}&
\end{eqnarray}
satisfying \eqref{characterization}. We define $\Phi_0$ 
to be the composition of the relative Frobenius $F:L^+P_\mu\rightarrow{^FL^+P_\mu}$ 
with $L^+(\intaut_\mu(p))$, i.e. the endomorphism on $L^+P_\mu$ 
that arises from applying the functor $L^+$ to the extension of 
$g\mapsto\mu(p)g\mu(\frac1p)$, in the sense of \eqref{movie}, 
as given by lemma \ref{IntAction}. In order to define $\Phi_1$ 
we can use the map $V_{\Fu_{\mu^{-1}}}^{-1}:
L^{>0}\ul{\Fu_{\mu^{-1}}}\rightarrow{^FL^+\ul{\Fu_{\mu^{-1}}}}$, 
as in \eqref{frob02}, and we appeal to corollary \ref{factc} for 
our transport of structure, namely by the isomorphism $\ul{\Fu_{\mu^{-1}}}\isoto U_{\mu^{-1}}$ which is induced 
by $\exp_{{U_{\mu^{-1}}}_K}$ where $K=W(k_0)[\frac1p]$.\\

Now let $h\in H_\mu(A_0)$ be the universal point of $H_\mu$ with values in its coordinate ring $A_0$, which is flat over $W(k_0)$. 
It is enough to show $^F(\mu(p)h\mu(\frac1p))\in{^FG(W(A_0))}$. 
We observe that we can apply proposition \ref{productI} over the localization $(1+pA_0)^{-1}A_0=R$. So there exists a factorization $h=h_0h_1$ into elements $h_0\in L^+P_\mu(R)$ and 
$h_1\in L^{>0}U_{\mu^{-1}}(R)$. 
It follows that $^F(\mu(p)h\mu(\frac1p))=\Phi_0(h_0)\Phi_1(h_1)\in{^FG(W(R))}$, 
which suffices, as we have $^F(\mu(p)h\mu(\frac1p))\in{^FG(W(A_0[\frac1p]))}$ 
too. Here note that localizations are flat and that 
$\{\Spec R,\Spec A_0[\frac1p]\}$ is an fpqc covering of 
$\Spec A_0$, as any prime ideal $\Fp\subset A_0$ with $\Fp\cap(1+pA_0)\neq\emptyset$ cannot have $pA_0\subset\Fp$.
\end{proof}

\begin{rem}
\label{Lieanalog}
In the sequel (in particular in the proof of proposition 
\ref{GMZCF}) we are using Lie-theoretic analogs of 
$H_\mu$ and $\Phi_\mu:H_\mu\rightarrow L^+{^FG}$, namely
\begin{eqnarray*}
&&\Fh_\mu(R):=
L^{>0}\ul{\Fu_{\mu^{-1}}}\times_{W(k_0)}L^+\ul{\Fp_\mu}\\
&&\phi_\mu:\Fh_\mu\rightarrow{^F{L^+\ul{\Fg}}}
\end{eqnarray*}
which we define by $V_{\Fu_{\mu^{-1}}}^{-1}$ on $L^{>0}\ul{\Fu_{\mu^{-1}}}$ and simply by $p^mF$ on the summand $L^+\ul{\Fg_m}$ of $\mu$-weight 
$m\geq0$. Actually, one recovers the map $$\Fh_\mu\rtimes H_\mu\to ^F{L^+\ul{\Fg}}\rtimes ^FL^+G;\quad(X,k)\mapsto(\phi_\mu(X),\Phi_\mu(k))$$ 
when working with $\ul{\Fg}\rtimes G$ in the proposition \ref{coFrob}. Please see \cite[subsubsection 3.5.3]{pappas} for a similar construction.
\end{rem}

\section{$(G,\mu)$-Displays}
For any affine group scheme $G$ over a ring $B$ we denote by $\Tor_S(G)$ the groupoid of $G$-torsors over a $B$-scheme $S$, and 
we write $\ton^G$ for the fpqc-stack over $\Spec B$ whose fibers over $S$ are given by $\ton^G(S):=\Tor_S(G)$. By pull-back we obtain 
$\ton_S^G:=S\times_{\Spec B}\ton^G$, which is an fpqc-stack over $S$. We need to introduce a couple of definitions:

\begin{defn}
\label{concept02}
Fix a finite field $k_0$. By a display datum over $W(k_0)$ we 
mean a pair $(G,\mu)$, where $G$ is a smooth affine group scheme 
over $\BZ_p$ and 
$$\mu:\BG_{m,W(k_0)}\rightarrow G_{W(k_0)}$$ 
is a cocharacter all of whose weights on the adjoint representation $\Fg=\omega^G(\Ad^G)$ are contained in the set $\{-1,0,1,2,\dots\}$, i.e.
\begin{equation}
\Fu_{\mu^{-1}}=\omega_\mu^G(-1,\Ad^G)
\end{equation}
\end{defn}

\begin{defn}
\label{concept01}
Fix a finite field $k_0$ and let $(G,\mu)$ be a display datum, as above.
\begin{itemize}
\item[(1)]
By the conjugated Frobenius map of $(G,\mu)$ we mean the morphism 
$$H_\mu\stackrel{\Phi_\mu}{\rightarrow}{L^+G_{W(k_0)}},$$ 
as exhibited in proposition \ref{coFrob}.
\item[(2)]
The \emph{stack of $(G,\mu)$-displays}, denoted by $\B(G,\mu)$, is the $W(k_0)$-stack rendering the diagram

\begin{equation}
\label{poshII}
\begin{tikzcd}
\B(G,\mu)\arrow[rrr]\ar[dd]&&&\ton^{L^+G}_{W(k_0)}\arrow[dd,"\Delta_{\ton^{L^+G}_{W(k_0)}}"]\\\\
\ton^{H_\mu}\arrow[rrr,"(\ton^{\Phi_\mu}\times\id)\circ\Delta_{\ton^{H_\mu}}" below]&&&\ton^{L^+G}_{W(k_0)}\times_{W(k_0)}\ton^{L^+G}_{W(k_0)}
\end{tikzcd}
\end{equation}

$2$-Cartesian.
\item[(3)]
A \emph{$(G,\mu)$-display} $\D$ over a $W(k_0)$-scheme $S$, is an $S$-point of the stack $\B(G,\mu)$, 
its image under the vertical projection is called the \emph{underlying $H_\mu$-torsor} and will be denoted by $q(\D)$. 
\item[(4)]
The precomposition of the natural $1$-morphism $\ton^{H_\mu}\rightarrow\ton^{P_\mu}$, induced by the projection ${H_\mu}\to P_{\mu}$ 
(cf. subsection \ref{productII}), with the vertical projection is called the \emph{lowest truncation} and will be denoted by 
\begin{equation}
\label{poshIV}
q_0:\B(G,\mu)\rightarrow\ton^{P_\mu}
\end{equation}
\end{itemize}
\end{defn}

Suppose that $S$ is a $W(k_0)$-scheme and that $\D$ is an $S$-valued section of $\B(G,\mu)$. We write 
$\P_{\D}\in\Ob_{\ton^{L^+G}(S)}$ (resp. $\Q_{\D}\in\Ob_{\ton^{H_\mu}(S)}$) for the $L^+G$-torsor (resp. $H_\mu$-torsor) 
over $S$, that arises from $\D$ by means of the horizontal (resp. vertical) projection. Moreover, we denote by
$$u_{\D}:\Q_{\D}*^{\Phi_\mu}L^+G\isoto\P_{\D}$$
the morphism making the diagram commute. Sometimes, when there is little risk of confusion, we drop the subscript 
$\D$ from the notation. Upon identifying $\P_{\D}$ with the structure group extension of $\Q_{\D}$ along the inclusion 
$H_{\mu}\into L^{+}G$, we see that the triple $\D=(\P_{\D}, \Q_{\D}, u_{\D})$ is determined by $(\Q_{\D}, u_{\D})$.\\

The generic fiber of $\B(G,\mu)$ has a rather ``pathological'' behavior:

\begin{lem}
\label{banana04}
The generic fiber of the lowest truncation \eqref{poshIV} is an isomorphism, and thus yields an equivalence:
$$\B(G,\mu)_{W(k_0)[\frac1p]}\cong\ton_{W(k_0)[\frac1p]}^{P_\mu}$$
\end{lem}
\begin{proof}
The diagram \eqref{poshII} being 2-Cartesian suggests immediately to consider the natural $1$-morphism
\begin{equation}
\label{banana03}
\ton^{H_\mu^{\Phi_\mu}}\rightarrow\B(G,\mu)
\end{equation}
where $H_\mu^{\Phi_\mu}\into H_\mu$ is defined to be the equalizer of $\Phi_{\mu}:H_{\mu}\to L^+G$ and the 
inclusion of $H_{\mu}$ into $L^+G$. In general the $1$-morphism \eqref{banana03} fails to be an isomorphism, as the 
criterion (i) of lemma \ref{FuncChangGroup5} is not satisfied. Nevertheless our proof is finished by the following remarks:
\begin{itemize}
\item
Upon a change of base to $W(k_0)[\frac1p]$ the criterion (i) of lemma \ref{FuncChangGroup5} is satisfied.
\item
Over $W(k_0)[\frac1p]$ the composition
$$H_\mu^{\Phi_\mu}\into H_\mu\onto P_\mu$$
becomes an isomorphism.
\end{itemize}
One easily checks the above by using the ghost coordinates of the Witt vectors functor. Indeed, if $R$ is a 
$\BZ[\frac1p]$-algebra, then the zeroth ghost component induces an isomorphism  \[w_0:W(R)^{F=\id}\xrightarrow{\cong} R\]
\end{proof}

\begin{defn}
\label{BanDisp}
A $(G,\mu)$-display $(\P,\Q,u)$ is called \emph{banal} if its underlying $H_\mu$-torsor $\Q$ is trivial.
\end{defn}

Banal $(G,\mu)$-displays are important for explicit computations, here are some facts:

\begin{lem}
\label{banana01}
Let $(G,\mu)$ and $k_0$ be as before, and let $\D$ be a $(G,\mu)$-display over an affine scheme $\Spec A$, where $A$ is a $W(k_0)$-algebra.
\begin{itemize}
\item[(i)]
Assume that $A$ is a linearly topologized complete ring and that $J$ is a closed ideal consisting of 
topologically nilpotent elements. Then $\D$ is banal if and only if this holds for its pull-back to $\Spec A/J$.
\item[(ii)]
Suppose that $pA\subset\sqrt{0_A}$ holds. Then $\D$ is banal if and only if $q_0(\D)$ has a global section.
\item[(iii)]
Suppose that $pA=A$ holds. Then $\D$ is banal if and only if $q_0(\D)$ has a global section.
\item[(iv)] 
The groupoid of banal $(G,\mu)$-displays over a $W(k_0)$-scheme $S$ is equivalent to the quotient groupoid $[L^+G(S)/ _{\Phi_\mu} H_\mu(S)]$ 
(cf. appendix \ref{app1st} for the notation). In particular the groupoids of banal $(G,\mu)$-displays over $S$ and over $\Spec\Gamma(S,\O_S)$ agree.
\end{itemize}
\end{lem}
\begin{proof}
The first statement follows from part (i) of lemma \ref{torsor03} while the second one can be handled using the method of \cite[Proposition B.0.15a)]{pappas} (see also \cite[Corollary B.0.16]{pappas}). The third statement is a consequence of lemma \ref{banana04}. For the fourth statement consult \cite[Subsubsection 3.2.7]{pappas}.
\end{proof}

\begin{rem}
Observe that lemma \ref{banana01} implies that a $(G,\mu)$-display $\D$ over a locally noetherian $W(k_0)$-scheme $S$ allows a banalization over an fpqc 
covering $\{U_i\}_{i\in I}$ where each $U_i$ is noetherian. More specifically, if $\{U_i\}_{i\in I}$ is an \'etale covering trivializing $q_0(\D)$ and such that each 
$U_i$ is affine, say equal to $\Spec A_i$, then $\D$ is banal over the localization $A_i[\frac1p]$ and the $p$-adic completion $\hat A_i$, and the collection 
$\{\Spec A_i[\frac1p]\}_{i\in I} \dot \cup \{\Spec\hat A_i\}_{i\in I}$ is indeed an fpqc covering of $S$ banalizing $\D$ and consisting of noetherian schemes.
\end{rem}

\begin{rem}
\label{RemBanalDisp}
Let $k_0\subset\Omega$ be a field extension with $\Omega$ algebraically closed, and let $\D$ be a $(G,\mu)$-display over 
$\Omega$. Then $\D$ is banal and so, by lemma \ref{banana01} (iv), it is represented by an element of $G(W(\Omega))$.
\end{rem}

\section{Deformation of Banal $(G,\mu)$-Displays}
\label{deform04}
Throughout this section, $k_{0}$ is a finite field and $(G,\mu)$ is a display datum over $W(k_{0})$ (cf. definition \ref{concept02}).\\

The following definition is adopted from \cite[Definition 3.4.2]{pappas} in conjunction with \cite[Lemma 3.4.4]{pappas}. 
We denote the absolute Frobenius endomorphism of a $\BF_p$-algebra $R_0$ by $\Frob_{R_0}$.

\begin{defn}
\label{stressful02}
Let $R$ be a $W(k_0)$-algebra, and let $\pi$ 
be the projector on $W(k_0)\otimes_{\BZ_p}\Fg$ with $\image(\pi)=\Fu_{\mu^{-1}}$ and $\ker(\pi)=\Fp_\mu$. An element $U\in G(R)$ 
is said to satisfy the \emph{adjoint nilpotence condition} if and only if the $p$-linear endomorphism on $R/pR\otimes_{\BZ_p}\Fg$ defined by
\begin{equation}
\label{stressful04}
{\Ad}^G(U)\circ(\Frob_{R/pR}\otimes{\id}_\Fg)\circ({\id}_{R/pR}\otimes\pi)
\end{equation}
is nilpotent. We will denote by $C_R^{nil}(G, \mu)$ the subset of $G(R)$ consisting of elements satisfying 
the adjoint nilpotence condition. Note that the construction of $C_R^{nil}(G, \mu)$ is functorial in $R$.
\end{defn}

\begin{lem}
\label{stressful01}
Let $R$ be a $W(k_0)$-algebra, and let $R_0$ be the perfection of $R/\sqrt{pR}$. 
\begin{itemize}
\item[(i)]
The subset $C_R^{nil}(G,\mu)\subset G(R)$ is the preimage of $C_{R_0}^{nil}(G,\mu)\subset G(R_0)$ under the canonical map $G(R)\to G(R_0)$.
\item[(ii)]
The subset $C_{W(R)}^{nil}(G,\mu)\subset L^+G(R)=G(W(R))$ is the preimage of $C_R^{nil}(G,\mu)$ in $G(R)$ via the 
zeroth ghost map $G(W(R))\to G(R)$, and it is invariant under $\Phi_\mu$-conjugation, i.e. for all $h\in H_\mu(R)$, we have 
$$h^{-1}C_{W(R)}^{nil}(G,\mu)\Phi_\mu(h)=C_{W(R)}^{nil}(G,\mu)$$
\end{itemize}
\end{lem}
\begin{proof}
That $C_R^{nil}(G,\mu)$ depends only on $R/pR$ was built into the very definition while passage to $(R/pR)_\red$ 
and $R_0$ does not alter the nilpotence of an $F$-linear endomorphism either. The next claim follows from the formula 
$I(R)^2=pI(R)$, which implies the following (note that $pW(R)+I(R)$ is the kernel of the projection $W(R)\to R/pR$): 
$$(pW(R)+I(R))^2=p(pW(R)+I(R))\subset pW(R)\subset pW(R)+I(R)$$
For the last assertion we may assume that $R$ is perfect (of characteristic $p$) and so consider the $F$-linear endomorphism
\begin{align*}
\phi_U^{\Ad}:W(R)[\frac1p]\otimes_{\BZ_p}\Fg&\rightarrow W(R)[\frac1p]\otimes_{\BZ_p}\Fg\\
X&\mapsto\Ad^G(U{^F\mu(p)}){^FX}.
\end{align*}
Observe that $p\phi_U^{\Ad}$ preserves the integral structure $\FF=W(R)\otimes_{\BZ_p}\Fg$ and that its reduction to $R\otimes_{\BZ_p}\Fg$ 
agrees with the endomorphism defined by the formula \eqref{stressful04}. It follows that $(\phi_U^{\Ad})^n(\FF)\subset p^{1-n}\FF$ 
holds for $n>>0$ if and only if $w_0(U)\in C_{R}^{nil}(G,\mu)$, which in turn is equivalent to $U\in C_{W(R)}^{nil}(G,\mu)$. The 
proof is completed by $\Ad^G(h)^{-1}\circ\phi_U^{\Ad}\circ\Ad^G(h)=\phi_{h^{-1}U\Phi_\mu(h)}^{\Ad}$ for any $h\in H_\mu(R)$.
\end{proof}

\begin{defn}
\label{stressful03}
Let $\D$ be a $(G,\mu)$-display over a $W(k_0)$-scheme $X$, on which $p$ is Zariski-locally nilpotent. Then $\D$ is called \emph{adjoint nilpotent} 
if each geometric fiber $\D\times_X\Spec\Omega$ is represented by an element in $C_{W(\Omega)}^{nil}(G,\mu)$, where $\Omega$ stands for 
an algebraically closed field. We denote by $\B^{nil}(G,\mu)$ the (full) substack of $\B(G,\mu)$ given by adjoint nilpotent $(G, \mu)$-displays.
\end{defn}

Recall that a pd-structure on an ideal $\Fa\subset R$ consists of a family of maps $\{\gamma_k:\Fa\rightarrow\Fa\}_{k\in\BN}$ such that the following holds:

\begin{itemize}
\item[(i)]
$\forall n\in\BN,\,x,y\in\Fa:\,\gamma_n(x+y)-\gamma_n(x)-\gamma_n(y)=\sum_{k=1}^{n-1}\gamma_k(x)\gamma_{n-k}(y)$
\item[(ii)]
$\forall n\in\BN,\,x\in R,\,y\in\Fa:\gamma_n(xy)=x^n\gamma_n(y)$
\item[(iii)]
$\forall n\in\BN,\,x\in \Fa,\,n!\gamma_n(x)=x^n$ 
\item[(iv)]
$\forall n,m\in\BN,\,x\in\Fa:\,\gamma_m(\gamma_n(x))=\frac{(mn)!}{m!n!^m}\gamma_{mn}(x)$
\end{itemize} 

Notice that a pd-structure $\{\gamma_k:\Fa\rightarrow\Fa\}_{k\in\BN}$ on an ideal $\Fa\subset R$ induces a $W(R)$-linear isomorphism
\begin{align}
\label{GMZCF03}
w':W(\Fa)&\isoto\prod_{n=0}^\infty\Fa_{[w_n]};\,(x_0,x_1,\dots)\mapsto(w'_0,w'_1,\dots),
\end{align}
which is given by Zink's \emph{logarithmic ghost coordinates}
$$w'_n(x_0,x_1,\dots,x_n)=\sum_{i=0}^n(p^i-1)!\gamma_{p^i}(x_{n-i})$$
(the usual ones being $w_n(x_0,x_1,\dots,x_n)=\sum_{i=0}^np^{n-i}x_{n-i}^{p^i}$). Note that the isomorphism above induces a canonical decomposition 
$$W(\Fa)\cong\Fa_{[w_0]}\oplus I(\Fa),$$ 
since $I(\Fa)\cong\prod_{n=1}^\infty\Fa_{[w_n]}$ (as $w'_{0}=w_{0}$). We write $W(\Fa)\onto I(\Fa);x\mapsto x^{>0}$ for 
the associated $W(R)$-linear projector whose kernel is  $\bigcap_{n=1}^\infty\ker w'_n\cong\Fa_{[w_0]}$, so
$$w'_n(x^{>0})=\begin{cases}w'_n(x)&n>0\\0&\mbox{$n=0$}\end{cases}.$$
Notice that the assignment $x\mapsto x^{>0}$ is additive and multiplicative, so that we obtain an induced group homomorphism from $L^+G(\Fa)$ to $L^{>0}G(\Fa)$.\\

Before we prove a proposition that shows the significance of adjoint-nilpotence, we need two auxiliary lemmas:

\begin{lem}
\label{GMZCF06}
Let $R$ be a $\BZ_p$-algebra, let $\Fa\subset R$ be a pd-ideal, and let $U$ be a $W(R)$-valued point of $G$. Then
\begin{itemize}
\item[(i)] for all $h\in G(W(\Fa))$, we have  $(UhU^{-1})^{>0}=Uh^{>0}U^{-1}$ 
\item[(ii)] for all $X\in W(\Fa)\otimes\Fg$, we have $(\Ad^G(U)X)^{>0}=\Ad^G(U)X^{>0}$.
\end{itemize}
\end{lem}
\begin{proof}
Observe that conjugation with $U\in G(W(R))$ (resp. $\Ad^G(U)$) defines an automorphism of $G\times_{\BZ_p}W(R)$ (resp. 
$\ul\Fg\times_{\BZ_p}W(R)$). The assertions follow as $W(\Fa)\rightarrow I(\Fa);x\mapsto x^{>0}$ is a homomorphism of 
non-unital $W(R)$-algebras (please see subsection \ref{properties} for our conventions on groups of points over non-unital algebras).
\end{proof}

\begin{lem}
\label{GMZCF05}
Let $R$ be a $W(k_0)$-algebra and let $\Fa\subset R$ be a pd-ideal with $p\Fa=0$. If $X\in I(\Fa)\otimes\Fg$ 
corresponds to $h\in G(I(\Fa))$ (via $I(\Fa)^2=0$), then $\phi_\mu(X)$ corresponds to $\Phi_\mu(h)$, where 
$\phi_\mu:\Fh_\mu\rightarrow L^+\ul\Fg\times_{\BZ_p}W(k_0)$ is the transformation constructed in remark \ref{Lieanalog}.
\end{lem}
\begin{proof}
Observe that we have a decomposition $$G(I(\Fa))\cong P_\mu(I(\Fa))\times U_{\mu^{-1}}(I(\Fa))$$ and its Lie analog 
$$I(\Fa)\otimes\Fg\cong I(\Fa)\otimes_{W(k_0)}\Fp_\mu\oplus I(\Fa)\otimes_{W(k_0)}\Fu_{\mu^{-1}}.$$ It follows that we may verify 
the assertions separately for $h\in P_\mu(I(\Fa))$ and $h\in U_{\mu^{-1}}(I(\Fa))$. In either case we deduce the lemma from 
the explicit descriptions of $\Phi_\mu$ on $L^+P_\mu$ and $L^{>0}U_{\mu^{-1}}$ (\eqref{nonpositive} and \eqref{positive}).
\end{proof}

\begin{prop}
\label{GMZCF}
Fix a $p$-adically separated and complete pd-ideal $\Fa$ in a $W(k_0)$-algebra $R$ (i.e. such that $R$ is separated and complete with respect to the 
topology defined by the descending chain of pd-ideals $\{p^i\Fa\}_{i\in\BN}$ by part (ii) of lemma \ref{Lem-pA+a-adic-top}). Let $O$ and $U$ be elements of 
$C^{nil}_{W(R)}(G,\mu)\subset G(W(R))$, satisfying  $UO^{-1}\in G(W(\Fa))$. Then there exists a unique $h\in G(W(\Fa))$ with $O=h^{-1}U\Phi_\mu(h^{>0})$.
\end{prop}
\begin{proof}
By a limit process followed by a straightforward induction we may assume $p\Fa=0$. 
According to the isomorphism \eqref{GMZCF03} we deduce $I(R)W(\Fa)=0$, given that
$$w_n({^Vx})=\begin{cases}pw_{n-1}(x)&n>0\\0&\mbox{$n=0$}\end{cases}$$
implies $w_n(I(R))\subset pR$. We deduce $pW(\Fa)=0$ by the same token, so that:
$$(pW(R)+I(R))W(\Fa)=0=(pW(R)+W(\Fa)+I(R))I(\Fa)$$
In particular we may and will regard 
$I(\Fa)$ (resp. $W(\Fa)$) as a module over $R/(pR+\Fa)$ (resp. over $R/pR$). Also, notice that the solution sets
\begin{eqnarray*}
&A=\{h\in G(W(\Fa))\,|h=U\Phi_\mu(h^{>0})O^{-1}\}&\\
&B=\{h\in G(I(\Fa))\,|h=(U\Phi_\mu(h)O^{-1})^{>0}\}&
\end{eqnarray*}
correspond to each other under the mutually inverse bijections: 
\begin{eqnarray*}
&&B\rightarrow A;\quad h\mapsto U\Phi_\mu(h)O^{-1}\\
&&A\rightarrow B;\quad h\mapsto h^{>0}
\end{eqnarray*} 
Although the set $A$ seems to be more natural to work with, we nevertheless prefer to work with the set $B$, because of the isomorphism
$$G(I(\Fa))\cong I(\Fa)\otimes_{\BZ_p}\Fg,$$
as $I(\Fa)^2=0$. In order to show that the set $B$ has exactly one element, we note that
$$(U\Phi_\mu(h)O^{-1})^{>0}=(U\Phi_\mu(h)^{>0}U^{-1})(UO^{-1})^{>0}$$
holds for all $h\in G(I(\Fa))$, according to lemma \ref{GMZCF06}. Let $C\in I(\Fa)\otimes\Fg$ correspond to the element 
$(UO^{-1})^{>0}\in G(I(\Fa))$. Consider the endomorphism of the $R/(pR+\Fa)$-module $I(\Fa)\otimes_{\BZ_p}\Fg$, sending an element $X$ to 
$$\Ad^G(U)\big(\phi_\mu(X)^{>0}\big)+C$$
 In view of lemma \ref{GMZCF05}, we ought to show that this endomorphism
has a unique fixed point. Our claim follows, once we prove that the $F$-linear endomorphism
$$X\mapsto M(X):=\Ad^G(U)\big(\phi_\mu(X)^{>0}\big)\overset{\eqref{GMZCF06}}{=}\big(\Ad^G(U)\phi_\mu(X)\big)^{>0}$$
is nilpotent on $I(\Fa)\otimes_{\BZ_p}\Fg$. In order to accomplish this it suffices to prove the nilpotence of
$$X\mapsto L(X):=\Ad^G(U)\phi_\mu(X^{>0})$$
on $W(\Fa)\otimes_{\BZ_p}\Fg$ (if $L^{m}=0$, then $M^{m+1}=0$). We can rewrite $W(\Fa)\otimes_{\BZ_p}\Fg$ as $W(\Fa)\otimes_{R/pR}(R/pR\otimes_{\BZ_p}\Fg)$ and the endomorphism $L$ 
thereon can be written as the tensor product of the endomorphisms $x\mapsto V^{-1}(x^{>0})$ on $W(\Fa)$ and the endomorphism \eqref{stressful04} 
on $R/pR\otimes_{\BZ_p}\Fg$, where the $R$-valued element of $G$ which appears in it is just the zeroth ghost component $w_0(U)$.
\end{proof}

\begin{rem}
Over artinian local rings with perfect residue field $k$ of characteristic $p\neq2$ one obtains an analog of the above if one replaces $W(R)$ by the Zink ring, 
which is the product of $W(k)$ with the points of the group of formal Witt vectors over the maximal ideal of $R$. In this theory one does not have to impose 
a nilpotence condition. Please see \cite[Proposition 7.1.5]{lau2} for details, and \cite{zink3} for further background material on such ``Dieudonn\'e displays''.
\end{rem}

\begin{rem}
The essence of the above proposition was already proved in \cite[Theorem 3.5.4]{pappas} under the slightly more restrictive assumption 
that $\Fa$ is a $p$-adically closed and topologically nilpotent pd-ideal in a $p$-adically separated and complete $W(k_0)$-algebra $R$. Our 
peculiar detour from $G(W(\Fa))$ via $G(I(\Fa))\cong I(\Fa)\otimes\Fg$ to $W(\Fa)\otimes\Fg$ allows to work with a weaker assumption on $\Fa$.\\

The classical method of \cite[proof of Theorem 44]{zink2} also works for all $p$-torsion pd-ideals $\Fa$ without any further 
restrictions. So up to our limit process the above lemma seems to be the ``natural'' generalization to $(G,\mu)$-displays!  
\end{rem}

\begin{rem}
Please see \cite[part (i) of lemma 3.24]{E7} for some variants.
\end{rem}

\subsection{Definition of $\Psi_\mu^\Fa$}
Recall that Zink's logarithmic ghost components induce a decomposition $W(\Fa)=\Fa_{[w_0]}\oplus I(\Fa)$ whenever $\Fa\subset R$ is an ideal 
with pd-structure $\{\gamma_k:\Fa\rightarrow\Fa\}_{k\in\BN}$. If the ambient ring $R$ has the structure of a $W(k_0)$-algebra we may deduce the 
decomposition $L^+G(\Fa)\cong G(\Fa)\times L^{>0}G(\Fa)$, so that $G(\Fa)$ can be regarded as a (normal) subgroup of both $G(R)$ and $L^+G(R)$. 
Also, we note that the absolute Frobenius $F:W(R)\rightarrow W(R)$ induces a shift of the logarithmic ghost components on $W(\Fa)$, so that 
$\Fa_{[w_0]}$ lies in the kernel of $F$. We deduce that the morphism $\Phi_\mu: H_\mu\rightarrow{^FL^+G}$ (cf. proposition \ref{coFrob}) annihilates 
$P_\mu(\Fa)=H_\mu(R)\cap G(\Fa)\subset H_\mu(R)$, regarding $G(\Fa)$ as a subgroup of $G(W(R))$ (cf. \cite[corollary 3.11]{E7}). We need the 
following generalization, in which a pd-morphism from $(R,\Fa,\gamma)$ to $(S,\Fb,\delta)$ is defined to be a ring homomorphism $f:R\rightarrow S$ 
satisfying $f(\Fa)\subset\Fb$ and $f(\gamma_k(x))=\delta_k(f(x))$ for all $k\in\BN$ and $x\in\Fa$, as in \cite[Chapter III, Definition (1.5)]{messing1}.

\begin{lem}
\label{LemExtPhi}
Consider the functors $\Gamma_{\mu}$ and $L^{+}G$ from the category of pd-structures $(R,\Fa,\gamma)$ over $(W(k_0),0_{W(k_0)},0)$ to the 
category of groups, sending a triple $(R,\Fa, \gamma)$ to $\Gamma_{\mu}^{W(R),W(\Fa)+I(R)}$ and respectively $L^{+}G(R)$. There exists a natural 
transformation \[\Psi_{\mu}:\Gamma_{\mu}\to L^{+}G\] denoted on a triple $(R,\Fa,\gamma)$ by $\Psi_{\mu}^{\Fa}$, satisfying the following two properties:
\begin{itemize}
\item
the restriction of $\Psi_\mu^\Fa$ to $H_\mu(R)=\Gamma_\mu^{W(R),I(R)}$ agrees with $\Phi_\mu$.
\item
$\Psi_\mu^\Fa$ vanishes on $G(\Fa)$ (note that this is a subgroup of $G(W(\Fa))$ which is in turn contained in $\Gamma_\mu^{W(R),W(\Fa)+I(R)}$).
\end{itemize}
\end{lem}
\begin{proof}
Consider the three localisations $R_p$, $\tilde R=(1+pR)^{-1}R$, $\tilde R_p$ and observe that 
$R=\tilde R\times_{\tilde R_p}R_p$, as $R\rightarrow\tilde R\times R_p$ is faithfully flat. Consequently we have:
\begin{equation*}
L^+G(R)\cong L^+G(\tilde R)\times_{L^+G(\tilde R_p)}L^+G(R_p)
\end{equation*}
Observe that $\tilde\Fa=\tilde R\otimes_R\Fa$ is still a pd-ideal of $\tilde R$ and so are $\Fa_p\subset R_p$ and $\tilde\Fa_p\subset\tilde R_p$. Furthermore, we have
\begin{equation*}
P_\mu(R/\Fa)\cong P_\mu(\tilde R/\tilde\Fa)\times_{P_\mu(\tilde R_p/\tilde\Fa_p)}P_\mu(R_p/\Fa_p),
\end{equation*}
as $R/\Fa\rightarrow\tilde R/\tilde\Fa\times R_p/\Fa_p$ is still faithfully flat. We deduce
\begin{equation*}
\Gamma_\mu^{W(R),W(\Fa)+I(R)}\cong\Gamma_\mu^{W(\tilde R),W(\tilde\Fa)+I(\tilde R)}
\times_{\Gamma_\mu^{W(\tilde R_p),W(\tilde\Fa_p)+I(\tilde R_p)}}\Gamma_\mu^{W(R_p),W(\Fa_p)+I(R_p)},
\end{equation*}
as $\Gamma_\mu^{W(R),W(\Fa)+I(R)}$ is the inverse image of $P_\mu(R/\Fa)$ in $L^+G(R)$, and likewise for $\Gamma_\mu^{W(R_p),W(\Fa_p)+I(R_p)}$, 
$\Gamma_\mu^{W(\tilde R),W(\tilde\Fa)+I(\tilde R)}$ and $\Gamma_\mu^{W(\tilde R_p),W(\tilde\Fa_p)+I(\tilde R_p)}$. Therefore it is enough to restrict our 
attention to triples $(R,\Fa,\gamma)$ satisfying either $pR=R$ or $pR\subset\rad(R)$. Recall that $pR\subset\rad(R)$ implies the direct generation of 
$\Gamma_\mu^{W(R),W(\Fa)+I(R)}$ by $L^+P_\mu(R)$ and $U_{\mu^{-1}}(W(\Fa)+I(R))$, as $pR\subset\rad(R)$ implies $\Fa\subset\rad(R)$, given that 
$\Fa\subset\sqrt{p\Fa}\subset\sqrt{pR}$ along with corollary \ref{productIII}, i.e.:
\begin{eqnarray}
&&L^+P_\mu(R)\times U_{\mu^{-1}}(W(\Fa)+I(R))\cong\\
&&\Gamma_\mu^{W(R),W(\Fa)+I(R)}\cong H_\mu(R)\times U_{\mu^{-1}}(\Fa)
\end{eqnarray}
Whence it follows $\Gamma_\mu^{W(R),W(\Fa)+I(R)}=G(\Fa)\Gamma_\mu^{W(R),I(R)}$, and the construction of $\Psi_\mu^\Fa$ is evident. In the 
other case $p$ is invertible in $R$ and we appeal to the formula $\Psi_\mu^\Fa(g)={^F(\mu(p)g\mu(\frac1p))}\in L^+G(R)$. It remains to check that the 
assignment $(R,\Fa,\gamma)\leadsto\Psi_\mu^\Fa$ is functorial with respect to $W(k_0)$-linear pd-morphisms $f:(R,\Fa,\gamma)\rightarrow(S,\Fb,\delta)$, 
where we may assume that $pS=S$ and $pR\subset\rad(R)$ (the other three subcases are trivial). To this end it is enough to check that $g\mapsto{^F(\mu(p)g\mu(\frac1p))}$ annihilates $G(\Fb)$, since we already know that $\Phi_\mu$ is a transformation from $H_\mu$ 
to $^FL^+G=L^+{^FG}=L^+G$ (cf. proposition \ref{coFrob}). We are done as $F$ annihilates $G(\Fb)=\mu(p)G(\Fb)\mu(\frac1p)$.
\end{proof}

Observe that $\Phi_\mu(h^{>0})=\Psi_\mu^\Fa(h)$ for any $h\in L^+G(\Fa)$. When using 
this language, proposition \ref{GMZCF} turns into the following more elegant version:

\begin{cor}
\label{GMZCF04}
Fix a $p$-adically separated and complete pd-ideal $\Fa$ in a $W(k_0)$-algebra $R$ (i.e. such that $R$ is separated 
and complete with respect to the topology defined by the descending chain of pd-ideals $\{p^i\Fa\}_{i\in\BN}$ 
by part (ii) of lemma \ref{Lem-pA+a-adic-top}). Let $O$ and $U$ be elements of $G(W(R))$, satisfying
$$O,U\in C_{W(R)}^{nil}(G,\mu)$$ and $$h^{-1}U\Psi_\mu^\Fa(h)O^{-1}\in G(W(\Fa)),$$
for some $h\in\Gamma_\mu^{W(R),W(\Fa)+I(R)}$. Then there exists a unique $\tilde h\in \Gamma_\mu^{W(R),W(\Fa)+I(R)}$ 
satisfying $\tilde hh^{-1}\in G(W(\Fa))$ and $O=\tilde h^{-1}U\Psi_\mu^\Fa(\tilde h)$.
\end{cor}

\subsection{Rigidity of isomorphisms}

This result for automorphisms can be regarded as an analog of a strengthening of \cite[Proposition 40]{zink2}:

\begin{thm}
\label{separ}
Let $k_{0}$ be a finite field and $(G,\mu)$ a display datum over $W(k_0)$. Let $\phi$ be an automorphism of a $(G,\mu)$-display $\D$ over a $W(k_0)$-algebra 
$R$, such that $\D\times_RR/pR$ is adjoint nilpotent. Consider an ideal $\Fa\subset R$ such that one of the following two assertions hold:
\begin{itemize}
\item[(i)]
$\Fa$ is $(pR+\Fa)$-adically separated.
\item[(ii)]
$\Fa$ allows a pd-structure and is $p$-adically separated i.e. $\bigcap_{i\in\BN_0}p^i\Fa$ vanishes.
\end{itemize}
Then $\phi=\id_{\D}$ if and only if $\phi\equiv\id_\D\mod\Fa$.
\end{thm}
\begin{proof}
Recall that by lemma \ref{Lem-pA+a-adic-top}, $\Fa$ is $(pR+\Fa)$-adically separated if and only if $\bigcap_{i\in\BN_0}(p^i\Fa+\Fa^{i+1})$ vanishes. Also, observe that the diagonal of our stack of displays with $(G,\mu)$-structure is separated according to 
\cite[Lemma 3.2.9b)]{pappas}. Towards checking $\phi=\id_\D$ it therefore suffices to prove $\phi\equiv\id_\D\mod(\Fa+pR)^i\Fa$ for all $i$ or 
$\phi\equiv\id_\D\mod p^i\Fa$ for all $i$, provided that $\Fa$ is a pd-ideal. By an easy induction it suffices to check this for $i=1$. This means we 
must check $\phi=\id_\D$ under the assumption $p\Fa=0$ where $\Fa$ satisfies $\Fa^2=0$ or has a pd-structure. Note that the former case follows from 
the latter, as all square-zero ideals can be endowed with a pd-structure. So, let us consider the latter case. After passage to an affine fpqc covering, we can 
assume that $\D$ is banal, hence represented by an element $O\in G(W(R))$. By assumption, the definition of adjoint-nilpotence and lemma \ref{stressful01}, we have 
$O\in C^{nil}_{W(R)}(G,\mu)$. The automorphism $\phi$ is represented by an element $h\in G(W(R))$ satisfying \[O=h^{-1}O\Phi_{\mu}(h^{>0})\] The assumption 
$\phi\equiv\id_{\D}\,\mod\Fa$ means that $h\in G(W(\Fa))$. The uniqueness part of proposition \ref{GMZCF} then implies that $h=1$, i.e. $\phi=\id_{\D}$.
\end{proof}

\subsection{Rigidity of quasi-isogenies}
In this subsection we give a proof of theorem \ref{rigid01}. In order to handle concepts of quasi-isogenies we have to introduce some non-representable fpqc sheaves.

\begin{defn}
\label{neighbor03}
Let $\Fa$ be an ideal in a ring $R$ and let $X$ be a scheme over $R$. We denote by $\xbar^\Fa$ (resp. $\xbar$) the fpqc sheaf on $\Spec R$ which is obtained as the direct image of $X\times_RR/\Fa$ (resp. $X\times\BZ[\frac1p]$) along the closed (resp. open) immersion from $\Spec R/\Fa$ (resp. $\BZ[\frac1p]\otimes R=:R[\frac1p]$) to $\Spec R$, i.e. 
\begin{eqnarray*}
&&\xbar^\Fa(T):=X(T\times_{\Spec R}\Spec R/\Fa)\\
&&\xbar(T):=X(T\times_{\Spec R}\Spec R[\frac1p]),
\end{eqnarray*}
where $T$ is any scheme over $R$.
\end{defn}

\begin{rem}
Let $f:T\rightarrow S$ be a morphism of schemes.  Any fpqc sheaf $\F$ on $S$ (resp. $\G$ on $T$) allows an inverse image 
$f^*\F:\Scheme_T\rightarrow\Ens;\,V\mapsto\F(V_{[f]})$ (resp. direct image $f_*\G:\Scheme_S\rightarrow\Ens;\,U\mapsto\G(T\times_{f,S}U)$). Moreover we have the adjunction isomorphism \[\Hom_T(f^*\F,\G)\cong\Hom_S(\F,f_*\G)\] which induces canonical adjunction morphisms $\F\rightarrow f_*f^*\F$ and $f^*f_*\G\rightarrow\G$. 
In the setting of the above definition \ref{neighbor03}, the $R$-scheme $X$ gives rise to a representable fpqc sheaf $h_{X}$ on $\Spec R$, and one obtains 
$\xbar^\Fa=i_*i^*h_{X}$ and $\xbar=j_*j^*h_{X}$ from the immersions $i:\Spec R/\Fa\hookrightarrow\Spec R$ and $j:\Spec R[\frac1p]\hookrightarrow\Spec R$.
\end{rem}

\begin{lem}
\label{neighbor07}
Let $R$ be a ring and let $X$ be a smooth and quasi-projective $R$-scheme. For any ideal $\Fa\subset R$, the adjunction homomorphism
\begin{equation}
\label{neighbor06}
X\rightarrow\xbar^\Fa
\end{equation}
is an epimorphism of fpqc-sheaves on $\Spec R$.
\end{lem}
\begin{proof}
The lemma boils down to the following claim: For every $R$-algebra $S$ and any $x_0\in X(S/\Fa S)$ there exists a section $y$ of $X$ over a 
faithfully flat $S$-algebra $T$ of which the $\mod\Fa$-reduction in $X(T/\Fa T)$ agrees with the image of $x_0$ therein. Since smooth quasi-projectivity 
is stable under base extensions it clearly suffices to assume $S=R$, moreover we can restrict to the noetherian case, according to Grothendieck's 
``\'Elimination des hypoth\`eses noeth\'eriennes'', just as in the proof of lemma \ref{torsor05}. Appealing to part (ii) of lemma \ref{torsor05} we 
may start out out from a lift $\hat x\in X(\hat R)$ of $x_0$, where $\hat R$ denotes the $\Fa$-adic completion of $R$. Observe that $\hat R$ is flat 
over $R$, but it might not be faithfully flat. So let $f_1,\dots,f_n$ generate the ideal $\Fa$ and choose an \'etale faithfully flat $R$-algebra $Q$ with 
$X(Q)\neq\emptyset$, according to \cite[Corollaire (17.16.3.ii)]{ega4}. The construction is completed 
by taking $T:=\hat R\times\prod_{i=1}^nQ_{f_i}$ 
and by taking $y$ to be any element of 
$X(\hat R)\times\prod_{i=1}^nX(Q_{f_i})=X(T)$ 
whose projection to $X(\hat R)$ agrees with $\hat x$.
\end{proof}

Now let $X$ be a scheme. The assignment $U\mapsto W(\O_X(U))$ defines a Zariski sheaf on $X$ that we denote by $W(\O_X)$. Also, $\BQ\otimes W(\O_X)$ denotes the Zariski sheaf obtained by sheafifying the presheaf $U\mapsto\BQ\otimes W(\O_X(U))$ on $X$. There is a canonical homomorphism from $\BQ\otimes W(\O_X(X))$ 
to $\Gamma(X,\BQ\otimes W(\O_X))$, which is bijective (resp. injective) for quasi-compact and quasi-separated (resp. quasi-compact) $X$. We need the following:

\begin{lem}
For every smooth affine group scheme $H$ over $W(k_0)$ the contravariant group functor\begin{equation*}
LH:\Scheme_{W(k_0)}\rightarrow\Grp:\:X\mapsto H(\Gamma(X,\BQ\otimes W(\O_X)))
\end{equation*}
called the \emph{loop group of $H$}, 
is a sheaf for the fpqc topology.
\end{lem}
\begin{proof}
Observe that for a covering of a $W(k_0)$-scheme 
$X$ by a family of open subsets $\{U_i\}_{i\in I}$ 
one obtains the exactness of
$$H\big(\Gamma(X,\BQ\otimes W(\O_X))\big)\hookrightarrow
\prod_{i\in I}H\big(\Gamma(U_i,\BQ\otimes W(\O_X))\big)\rightrightarrows\prod_{i,j\in I}H\big(\Gamma(U_i\cap U_j,\BQ\otimes W(\O_X))\big)$$
from the mere sheaf property of $\BQ\otimes W(\O_X)$. 
It remains to check that the same holds for
$$LH(X)\hookrightarrow LH(Y)\rightrightarrows LH(Y\times_XY)$$
for any faithfully flat and quasi-compact morphism $\phi:Y\rightarrow X$. By what we just proved, 
this is local with respect to the Zariski topology 
on $X$. Without loss of generality, we may then assume that $X$ is of the form $\Spec A$, in which case $Y$ 
is covered by finitely many open affines, say 
$\Spec B_1,\dots,\Spec B_n$. Since it is allowed 
to refine our covering we may replace $Y$ by the disjoint union of these, so that it becomes affine 
too, say of the form $\Spec B$. We should therefore 
show the exactness of the following:
\[H\big(\Gamma(X,\BQ\otimes W(\O_{X}))\big)\into H\big(\Gamma(Y,\BQ\otimes W(\O_{Y}))\big)
\rightrightarrows H\big(\Gamma(Y\times_{X}Y,\BQ\otimes W(\O_{Y\times_{X}Y}))\big)\]
It is enough to show that the following is exact:
\[\Gamma(X,\BQ\otimes W(\O_{X}))\into \Gamma(Y,\BQ\otimes W(\O_{Y}))\rightrightarrows \Gamma(Y\times_{X}Y,\BQ\otimes W(\O_{Y\times_{X}Y}))\]

Since $X$ and $Y$ are affine, and so are quasi-compact and quasi-separated, it follows from what we said 
just before the lemma, this sequence is identified 
with the following sequence:
\begin{align}
\label{WittSeq}
\BQ\otimes W(A)\into\BQ\otimes W(B)\rightrightarrows\BQ\otimes W(B\otimes_AB)
\end{align}
Since $A\to B$ is faithfully flat, the sequence \[A\hookrightarrow B\rightrightarrows B\otimes_AB\] is exact (cf. 
\cite[Chapter I, Proposition 2.18]{milne1}). By \cite[Lemma 30]{zink2} the sequence \[W(A)\hookrightarrow W(B)\rightrightarrows W(B\otimes_AB)\] 
is exact. This implies the exactness of \eqref{WittSeq}, since localization is an exact functor.
\end{proof}

\begin{rem}
One can think of $LH$ as the fpqc sheafification of the functor $X\mapsto H(\BQ\otimes\Gamma(X,W(\O_X)))$, 
but since sheafification is not always a well-defined operation in the fpqc topology, our detour via the Zariski topology 
seems to be inevitable. Please see \cite{waterhouse} for some presheaves which do not possess fpqc sheafifications.
\end{rem}

The fpqc sheaf $LH$ is a prime example of a non-representable group functor. The following two 
lemmas investigate basic properties of the ``neutral locus'' of an automorphism of a $LH$-torsor:

\begin{lem}
\label{rigid10}
Let $R$ be a $W(k_0)$-algebra. The set 
\begin{equation}
\label{rigid08}
\Fa_n=\{x\in R\mid\;p^nx=p^{n-1}x^p=\dots=px^{p^{n-1}}=x^{p^n}=0\}
\end{equation}
is an ideal for each $n\in\BN_0$. Moreover, one has $pW(\Fa_n)\subset F(W(\Fa_n))\subset W(\Fa_{n-1})$ for all $n>0$.
\end{lem}
\begin{proof}
Note that we have 
\begin{equation}
\label{rigid09}
\Fa_n=\{x\in R\mid\;px,x^p\in\Fa_{n-1}\}
\end{equation}
We show by induction that $\Fa_n$ is an ideal. Observe that $\Fa_0=0$ and thanks to the identity 
\[(x+y)^p=x^p+y^p+pxy(\sum_{i=1}^{p-1}\frac{(p-1)!}{i!(p-i)!}x^{i-1}y^{p-i-1})\] if $\Fa_{n-1}$ is an ideal, $\Fa_n$ is an ideal too.\\

The formula $p=FV$ on $W(R)$ implies $pW(\Fa_n)\subset F(W(\Fa_n))$. In order to show $F(W(\Fa_n))\subset W(\Fa_{n-1})$ it suffices to restrict to $n=1$ 
because of \eqref{rigid09}. Consider the $W(k_0)$-algebra $R_1:=W(k_0)\oplus\Fa_1$, and notice that $\Fa_1$ can be recovered as $\{x\in R_1\mid\;x^p=px=0\}$. Observe that $R'_1:=R_1/pR_1\cong k_0\oplus\Fa_1$, and notice that $\Fa'_1:=\{x\in R'_1\mid\;x^p=px=0\}$ is still canonically isomorphic to $\Fa_1$. Since the 
absolute Frobenius on $W(R'_1)$ is given by the familiar formula $(x_0,x_1,\dots)\mapsto(x_0^p,x_1^p,\dots)$ we deduce $F(W(\Fa'_1))=0$, so that 
$F(W(\Fa_1))\subset W(pR_1)$. However, we have $F(W(\Fa_1))\subset W(\Fa_1)$ too, and $pR_1\cap\Fa_1=0$. It follows that $F(W(\Fa_1))$ vanishes in 
the ring $W(R_1)$. At last we easily deduce the result for $R$, since the canonical map $R_1\rightarrow R$ induces an isomorphism of the respective $\Fa_1$'s. 
\end{proof}

\begin{lem}
\label{rigid07}
Let $H$ be a smooth affine group scheme over $W(k_0)$. Let $\alpha:P\isoto P$ be an automorphism of an $LH$-torsor $P$ over a $W(k_0)$-scheme 
$X$ and consider the subset $X_0=\{x\in X\mid\,\alpha_x=\id_{P_x}\}$, where $\alpha_x:P_x\isoto P_x$ denotes the fiber of $\alpha$ at $x$. 
Let $\alpha_{X_\red}:P_{X_\red}\isoto P_{X_\red}$ stand for the pull-back of $\alpha$ to the reduction $X_\red$. Then the following hold:
\begin{itemize}
\item[(i)]
$X_0$ is closed with respect to the Zariski topology of $X$.
\item[(ii)]
$\alpha_{X_\red}=\id_{P_{X_\red}}$ if and only if $X_0=X$.
\item[(iii)]
Assume that $X$ is locally noetherian and that its generic fiber $X\times_{W(k_0)}W(k_{0})[\frac1p]$ 
is reduced. If $P$ is trivial, then $\alpha_{X_\red}=\id_{P_{X_\red}}$ if and only if $\alpha=\id_P$.
\end{itemize}
\end{lem}

\begin{proof}
We can assume that $X$ is affine, say it is $\Spec R$. Statements (i) and (ii) are not changed if we replace $X$ by $X_{\red}$, and so, in those two statements, we assume further that $R$ is reduced. Since by \cite[Corollaire (2.3.12)]{ega2} (see also by 
\cite[Lemme 1.2]{fga12}) fpqc morphisms are submersive, and $P$ trivializes after base change of $X$ by a faithfully flat morphism, we can assume that $P$ is already trivial over $X$. So, in all three statements, $X=\Spec R$ and $P$ is trivial over $X$ and so, the automorphism $\alpha$ is represented by an element $A\in LH(\Spec R) = H(\BQ\otimes W(R))$ and we want to investigate when it is equal to the identity element. In (i) and (ii) $R$ is reduced, in (iii) $R$ is noetherian.\\

(i): Pick a closed embedding $\rho:H\into\GL_{n}$ over $W(k_{0})$. The element $\rho(A)-1\in \BM_{n}(\BQ\otimes W(R))$ 
can be written as $\frac{1}{p^m}B$ where $m\geq 0$ is a natural number and $B\in \BM_{n}(W(R))$. It then follows 
that $X_{0}$ is the locus of vanishing of all Witt components of all entries of the matrix $B$, which is closed in $X$.\\

(ii): Assume that $X_0=X$ and let $B$ be as in (i). By the construction of $X_0$ in (i) and by assumption, all Witt 
components of all entries of $B$ belong to all prime ideals of $R$. Since $R$ is reduced, $B$ is the zero matrix.\\

(iii): We ought to show that $A$ is the identity element of $H(\BQ\otimes W(R))$ given that its reduction in $H(\BQ\otimes W(R/\sqrt{0_{R}}))$ 
is the identity. We are going to show that the canonical morphism \[\BQ\otimes W(R)\to \BQ\otimes W(R/\sqrt{0_{R}})\] is in fact an 
injection, which would then imply that $A$ is the identity. It is enough to show that $W(\sqrt{0_{R}})$ is a subset of $W(R)_{\tor}$, the 
$p$-power torsion elements of $W(R)$. Consider the ascending chain of ideals $\Fa_n$, as in \eqref{rigid08}. For all $n\geq 0$ we have:
\begin{align}
\label{Wa_nSubWR-tor}
W(\Fa_n)\subset W(R)_{\tor}
\end{align} 
The union of $\Fa_n$ is equal to $R_\tor\cap\sqrt{0_R}$, which is equal to $\sqrt{0_{R}}$, since by assumption $R\otimes_{W(k_0)}W(k_0)[\frac1p]$ 
is reduced. On the other hand, since $R$ is noetherian, this union is equal to $\Fa_n$ for $n\gg 0$. To summarize, for $n\gg 0$, 
we have $\Fa_n=\sqrt{0_R}$. Thus, by \eqref{Wa_nSubWR-tor}, we have $W(\sqrt{0_R})=W(\Fa_n)\subset W(R)_{\tor}$ as desired.
\end{proof}

From now on we work again with the display datum  $(G,\mu)$ over $W(k_{0})$. The relevance of 
$LG$ in this context stems from $^F\mu(p)\in LG(W(k_0))$ and the commutative diagram 
\[\begin{tikzcd} L^+G \ar[r, "F"] & L^+G\ar[r] & LG\ar[d,"\intaut({^F\mu}(p))"]\\ H_{\mu}\ar[u,hook]\ar[r,"\Phi_{\mu}"']& L^+G\ar[r] & LG\end{tikzcd}\]
which implies that for any $H_\mu$-torsor $Q$, we have a canonical isomorphism:
$$^F(Q*LG)\isoto Q*^{\Phi_\mu}LG$$
Also, observe that we have a commutative diagram:
\[\begin{tikzcd}
H_{\mu}\ar[d,hook]\ar[r] & H_{\mu}^{\BQ}\ar[d,hook]\\
L^+G\ar[r]&LG
\end{tikzcd}\]
where $H_\mu^\BQ:=\pbar_\mu\times_{\gbar,w_0^G}LG$. 

\begin{defn}
Let $X$ be a scheme over $\BZ_p$ (resp. over $W(k_0)$). By a \emph{$G$-isodisplay} (resp. \emph{filtered $(G,\mu)$-isodisplay}) 
we mean a pair $(P,b)$ where $P$ is a torsor 
under $LG$ (resp. $H_\mu^\BQ$) over $X$ and 
$b$ is an isomorphism
$$b:{^F(P*LG)}\isoto P*LG$$
of $LG$-torsors over $X$. We call $(P,b)$ \emph{banal} if $P$ is the trivial torsor.\\

A morphism between $G$-isodisplays (resp. filtered $(G,\mu)$-isodisplays) $(P_1,b_1)$ and $(P_2,b_2)$ 
is a morphism $\alpha:P_1\rightarrow P_2$ between 
their underlying torsors, such that $b_2\circ{^F\alpha}\circ b_1^{-1}=\alpha$ (these 
are morphisms between the $LG$-torsors $P_1*LG$ 
and $P_2*LG$ over $X$). 
\end{defn}

\begin{defn}
Let $\D=(Q,u)$ be a $(G,\mu)$-display over $X$. Then, 
as we discussed above, we have canonical isomorphisms:
\[^F(Q*LG)\isoto Q*^{\Phi_\mu}LG\overset{u}{\cong}Q*LG\]Therefore, by change of group from 
$H_{\mu}$ to $LG$ (resp. $H_{\mu}^{\BQ}$), the 
$(G,\mu)$-display $\D$ gives rise to a canonical 
$G$-isodisplay (resp. filtered $(G,\mu)$-isodisplay), which we denote by $\D^0$ (resp. $\D^\BQ$).
Let $\D_{1}$ and $\D_{2}$ be $(G,\mu)$-displays over 
$X$. A \emph{quasi-isogeny} from $\D_1$ to $\D_2$ 
over $X$, denoted by $\D_1\dashrightarrow\D_2 $ 
is an element of $\Hom(\D_1^\BQ,\D_2^\BQ)$. 
\end{defn}

\begin{rem}
The groups $\Hom(\D_1^0,\D_2^0)$, $\Hom(\D_1^\BQ,\D_2^\BQ)$, and $\Hom(\D_1,\D_2)$ 
are closely related to each other, 
observe that we have:
$$\Hom(\D_1^\BQ,\D_2^\BQ)=\Hom_{H_\mu^\BQ}(\D_1^\BQ,\D_2^\BQ)\times_{\Hom_{LG}(\D_1^0,\D_2^0)}\Hom(\D_1^0,\D_2^0)$$
It follows that we have a natural injection $\Hom(\D_1^\BQ,\D_2^\BQ)\into\Hom(\D_1^0,\D_2^0)$, 
since $H_\mu^\BQ$ is a subsheaf of $LG$. In 
particular, the content of theorem \ref{rigid01} 
remains unaltered if each occurrence of ``$\D^0$'' 
is replaced by the arguably more natural ``$\D^\BQ$''.
\end{rem}

\begin{prop}
\label{rigid06}
Let $A$ be a $W(k_0)$-algebra and fix an element $U\in C_{W(A)}^{nil}(G,\mu)$. Then the equation 
$$g^{-1}(U{^F\mu(p)}){^Fg}=U{^F\mu(p)}$$ 
has no non-trivial solutions in elements $g\in G(\BQ\otimes W(\Fa))$, in each of the following two cases:
\begin{itemize}
\item[(i)]
$\Fa=\sqrt{pA}$ and $A$ is noetherian and $p$-adically separated and complete.
\item[(ii)]
$A$ is a complete discrete valuation ring with maximal ideal $\Fa$ and its residue field contains $k_0$.
\end{itemize}
\end{prop}

\begin{proof}
Let $b:=U{^F\mu(p)}$, so that $g^{-1}b{^Fg}=b$, where $g$ is as in the proposition. Lemma \ref{rigid10} shows that if $\Fb$ is an ideal, 
and there exists an $n$ such that for all $x\in \Fb$, $x^{n}=0=p^{n}x$, then $W(\Fb)$ is $p$-power torsion, and so $\BQ\otimes W(\Fb)=0$. 
So, in order to deal with (i) we may replace $\Fa=\sqrt{pA}$ by $pA$, since $\BQ\otimes W(\sqrt{pA}/pA)=0$ and therefore we have:
$$\BQ\otimes W(pA)\cong \BQ\otimes W(\sqrt{pA})$$ 
Moreover, we can assume that $A$ has no $p$-torsion, which is due to
$$\BQ\otimes W(pA)\cong\BQ\otimes W(pA_0),$$
where $A_0$ is the quotient of $A$ by its $p$-torsion (again, we have $\BQ\otimes W\big(\ker(pA\onto pA_{0}) \big)=0$). Now observe that $F$ maps $W(pA)$ to $pW(pA)$, so there exists some $s\in\BN$ with $^{F^s}g=:\tilde g\in G(W(pA))$. Let us put $^{F^s}U=:\tilde U$ and $^{F^s}\mu=:\tilde\mu$ and 
$^{F^s}b=\tilde U{^F\tilde\mu(p)}=:\tilde b$. Observe that $\tilde U$ is an element of $C_{W(A)}^{nil}(G,\tilde\mu)$ and that $\tilde g^{-1}\tilde b{^F\tilde g}=\tilde b$ holds.
It follows that \[\tilde g^{-1}\tilde U^{F}\big(\tilde\mu(p){\tilde g} {\tilde\mu(\frac{1}{p})}\big)=\tilde U\] and so, $\tilde g$ is trivial by proposition \ref{GMZCF}. By iterating 
the identity $g=b{^{F}g}b^{-1}$, one concludes that $g$ and $^{F^{s}}g=\tilde g$ are conjugates and so $g$ is trivial as well.\\

In case (ii) we must focus on the equal characteristic case, since we could appeal to part (i) otherwise. 
Let $k$ denote the residue field of $A$. Our proof starts out from the Cohen structure theorem, which implies that $A$ is isomorphic to the ring $k[[t]]$ 
of power series with coefficients in $k$ in an indeterminate $t$. It does no harm to assume that $k$ 
is algebraically closed. Observe that the equation $g=b{^Fg}b^{-1}$ implies that $g$ is a unipotent 
element (of say $G\big(\BQ\otimes W(k((t)))\big)$), indeed look at some faithful representation $\rho:G_{\BQ_p}\rightarrow\GL(n)_{\BQ_p}$ and 
consider the characteristic polynomial $\chi\in\big(\BQ\otimes W(k[[t]])\big)[x]$ of $\rho(g)$. Notice that since $g$ is conjugated to $^Fg$, 
the coefficients of $\chi$ must be fixed by $F$, and 
so, they are elements of $\BQ\otimes W(\{a\in k((t))\mid\:a=a^p\})=\BQ_p$. Also, notice that 
$$\chi\equiv(x-1)^n\mod{\BQ\otimes W(tk[[t]])},$$
as $\rho(g)$ specializes to the identity matrix at 
$t=0$. We deduce $\chi=(x-1)^n$, which is due to $\BQ_p\cap\BQ\otimes W(tk[[t]])=0$. Since $g$ is unipotent it is allowed to consider its logarithm $\log(g)=:X$, which satisfies the identity 
$X=\Ad(b){^FX}$ and the rest of the proof is the 
lemma \ref{rigid02} below.
\end{proof}

\begin{lem}
\label{rigid02}
Let $k$ be an algebraically closed field of charachteristic $p$ and let $B$ be a 
$d\times d$-Matrix with entries in $\BQ\otimes W(k[[t]])$. If there exists a number $s\in\BN$ with 
\begin{equation}
\label{rigid03}
B\cdot{^FB}\dots{^{F^{s-1}}B}\in \frac{1}{p^{s-1}}\BM_{d}\big(W(k[[t]])\big),
\end{equation}
then the equation 
\begin{equation}
\label{rigid05}
x=B\cdot{^Fx}
\end{equation}
has no non-trivial solutions in column vectors $x$ with entries in $\BQ\otimes W(tk[[t]])$.
\end{lem} 
\begin{proof}
Observe that \cite[Definition 5.15; Proposition 5.16]{effeff} introduced, for any $r\in \BR_{\geq 0}$, valuations: 
\begin{equation}
\label{rigid04}
v_r:\BQ\otimes W(\widehat{k[[t]]^{perf}})\rightarrow\BR\cup\{\infty\};\:\sum_{\BZ\ni i\gg-\infty}p^i[a_i]\mapsto\inf\{ri+v(a_i)\mid i\in\BZ\}
\end{equation}
where $v:\widehat{k((t))^{perf}}\rightarrow\BZ[\frac1p]\cup\{\infty\}$ is the usual $t$-adic valuation (see also \cite{messing3}). 
When regarding $\sum_{i\in\BZ}p^i[a_i]$ as a bi-Witt vector its $i$th component is given by $a_i^{p^i}$, when using the notation 
of \eqref{rigid04}. Also, notice that we have $v_r([t])=1$ and $v_r(p)=r$, and observe that every $a\in W(tk[[t]])$ satisfies 
$$v_r(a)\geq\frac r{\log p}(1+\log\log p-\log r),$$
simply because the $i$th Witt component $a_i^{p^i}$ of $a$ is an element of $tk[[t]]$, while
$$ri+\frac1{p^i}\geq\frac r{\log p}(1+\log\log p-\log r)$$
by a little bit of high school calculus (look at $i=\frac{\log\log p-\log r}{\log p}$, which is the 
only zero of $r-\frac{\log p}{p^i}$ being the derivative with respect to $i$ of the left hand side).\\

Now assume that \eqref{rigid05} holds for a column vector $x=(x_1,\dots,x_d)^{T}$, and let $s$ be as in \eqref{rigid03}, so that we have 

\begin{align}
\label{X=CFX}
x=C\cdot{^{F^s}x}
\end{align}

where $C:=B\cdot{^FB}\dots{^{F^{s-1}}B}\in \frac{1}{p^{s-1}}\BM_{d}\big(W(k[[t]])\big)$. Notice that $v_r({^Fx_i})=pv_\frac rp(x_i)$ holds for $i\in\{1,\dots,d\}$, 
so that looking at the function 
$$\nu(r):=\min\{v_r(x_1),\dots,v_r(x_d)\}\geq\frac r{\log p}(1+\log\log p-\log r),$$
allows us to deduce $$\nu(r)\geq(1-s)nr+p^{sn}\nu(\frac r{p^{sn}})$$ because of the identity
$$x=C\cdot{^{F^s}C}\dots{^{F^{s(n-1)}}C}\cdot{^{F^{sn}}x}$$
which follows by iterating \eqref{X=CFX}. It follows that we have 
$$\nu(r)\geq rn+\frac r{\log p}(1+\log\log p-\log r)$$ for every $n$, which shows that the valuations of $x_{i}$ are not bounded and therefore $x_{i}$ are zero as desired.
\end{proof}

We are now all set to prove the following, where $\D_x$ denotes the pull-back of a $(G,\mu)$-display $\D$ along the morphism $\Spec\kappa(x)\rightarrow X$ for any $x\in X$
of which the stalk $\O_{X,x}$ gives rise to its residue field $\kappa(x)$.

\begin{thm}[Rigidity of Quasi-Isogenies]
\label{rigid01}
Let $k_{0}$ be a finite field and $(G,\mu)$ a display datum over $W(k_0)$. Let $\D$ be a $(G,\mu)$-display over a $W(k_0)$-scheme $X$ whose restriction to the special fiber $X\times_{W(k_0)}k_0$ is adjoint nilpotent.
\begin{itemize}
\item[(i)]
If $X$ is connected and locally noetherian and its structural morphism to $\Spec W(k_0)$ is assumed to be \emph{closed} (in the sense of \cite[(2.2.6)]{egai} or \cite[\href{https://stacks.math.columbia.edu/tag/004E}{Tag 004E}]{stacks-project}), then
a self quasi-isogeny $\alpha\in\Aut(\D^0)$ agrees with the identity if and only if this holds for at least one of its specializations $\alpha_x\in\Aut(\D_x^0)$.
\item[(ii)] 
If $p$ is assumed to be locally nilpotent on $\O_X$, then the natural map $\Aut(\D)\rightarrow\Aut(\D^0)$ is 
injective (i.e. $\D\mapsto\D^0$ is faithful).
\end{itemize}
\end{thm}

\begin{proof}
We begin with the non-trivial direction of (i), so let the specialization of $\alpha$ to some point of $X$ be the identity. Let $X_0\subset X$ be the set of all $x$ such 
that $\alpha_x$ agrees with the identity element of $\Aut(\D_x^0)$. By assumption $X_0$ is not empty and by part (i) of lemma \ref{rigid07} $X_0$ is closed. 
In order to prove that $X_0$ is open, it suffices to check that it is stable under generation, in view of \cite[Th\'eor\`eme (1.10.1)]{ega1} or \cite[\href{https://stacks.math.columbia.edu/tag/0542}{Tag 0542}]{stacks-project} as closed sets are certainly constructible. So assume that some $x\in X$ specializes to some 
$x_0\in X_0$. It does no harm to assume that the characteristic 
of $\kappa(x_0)$ is $p$, because $X\rightarrow\Spec W(k_0)$ 
is a closed map and so the special point of $\Spec W(k_{0})$ belongs to the image of the non-empty closed subset 
$\ol{\{x\}}\cap X_{0}$. According to a classical specialization result (e.g. \cite[Proposition (7.1.9)]{egaii}, or \cite[\href{https://stacks.math.columbia.edu/tag/01J8}{Tag 01J8}]{stacks-project}) there exists a discrete valuation ring 
$R$ and a morphism $\Spec R\rightarrow X$ that sends the generic point to $x$ and the special point to $x_0$, furthermore we may assume that $R$ is complete and has an algebraically closed residue field (of characteristic $p$). It comes in handy that a $(G,\mu)$-display over such a ring $R$ is banal. 
Therefore we are allowed to apply part (ii) of proposition \ref{rigid06} and deduce $x\in X_0$, as required. The 
connectedness of $X$ enforces $X_0=X$.\\

Now that we have established that all specializations of 
$\alpha$ are identity, we will use a descent argument to show 
that $\alpha$ is the identity (globally). So choose an affine \'etale covering $\{U_i\}_{i\in I}$ trivializing $q_0(\D)$. 
Let us write $A_i$ for the noetherian $W(k_0)$-algebra with 
$\Spec A_i=U_i$ and pick $f_i\in1+pA_i$ with $f_i(\bigcap_{n\in\BN_0}p^nA_i)=0$ (use Artin-Rees and 
Nakayama). Observe that $U_i$ and $D(f_i)=\Spec A_i[\frac1{f_i}]$ 
have the same fiber over $\Spec k_0$, as $f_i\equiv1\pmod p$. 
It follows that the special fiber of $X$ is contained in the 
image of the map
$$\coprod_{i\in I}D(f_i)\rightarrow X$$ 
which is an open subset of $X$ (the composition 
$D(f_{i})\to U_{i}\to X$ is flat, and so is an open map). Since $X\to\Spec W(k_{0})$ is a closed map, the image of the above map must be the whole of $X$ and therefore, $\{D(f_i)\}_{i\in I}$ is still an \'etale covering of $X$. We are going to show that the pull-back of $\alpha$ along this covering is the identity and therefore, it must be the identity as well. So, we need to 
show that the pull-back of $\alpha$ to each $A_{i}[\frac1{f_{i}}]$ is the identity. Using parts (i) and (ii) of lemma \ref{banana01}, we conclude that the pull-back of $\D$ to the $p$-adic completion $\hat A_i$ of $A_i$ is banal. Moreover, the reducedness of $A_i/\sqrt{pA_i}$ allows to apply part (ii) of lemma \ref{rigid07} and we deduce that the pull-backs of $\alpha$ to each of 
those $A_i/\sqrt{pA_i}$'s are identities. Furthermore, over 
each $\hat A_i$ we can take advantage of part (i) of proposition \ref{rigid06} and deduce that the pull-backs of $\alpha$ to 
each of those $\hat A_i$'s are still identities. At last the 
pull-backs of $\alpha$ to each $A_i[\frac1{f_i}]$ must agree 
with the identity too, because the map 
$A_i[\frac1{f_i}]\rightarrow\hat A_i[\frac1{f_i}]$ is injective.\\

In order to proceed to assertion (ii), we may assume $X=\Spec R$ 
for some $W(k_{0})$-algebra $R$ satisfying $pR\subset\sqrt{0_R}$ 
and that $\D$ is banal, say with representative $U\in G(W(R))$. 
Part (i) of theorem \ref{separ} tells us that some 
$H_\mu(R)\ni h=U\Phi_\mu(h)U^{-1}$ is neutral if and only 
if $h\in H^\mu(\Fa)$ for some nilpotent ideal $\Fa\subset R$. 
This allows us to assume $pR=0$. Furthermore, $h$ and $1$ 
give rise to the same self quasi-isogeny if and only if 
$h\in G(W(R)_{\tor})$. In view of the finite presentation of $G$, 
the latter condition is equivalent to $h\in\bigcup_lH_\mu(\Fa_l)$ where $\Fa_l$ is as in lemma \ref{rigid10}. By working with a 
minimal counterexample, it remains to check that
$$H_\mu(\Fa_1)\ni h=U\Phi_\mu(h)U^{-1}$$
does not have non-trivial solutions (N.B. recall 
$\Fa_l=\{x\in R\mid\,x^p\in\Fa_{1-1}\}$). Observe that one has $I(R)W(\Fa_1)=0$ and $W(\Fa_1)^2\subset\hat W(\Fa_1)$, where
$$\hat W(\Fa):=\{(x_0,x_1,\dots)\in W(\Fa)\mid\,\sharp\{i\mid\,x_i\neq0\}<\infty\}$$
denotes the ideal of formal Witt vectors over 
$\Fa\subset\sqrt{0_R}$. Indeed, the first claim follows 
from $^Vxy={^V(x{^Fy})}$ and the second one from 
$W(\Fa_1)=\hat W(\Fa_1)+I(\Fa_1)$, so that 
$(\hat W(\Fa_1)+I(\Fa_1))^2=\hat W(\Fa_1)^2$. The extension
$$S:=W(\Fa_1)/\hat W(\Fa_1)\hookrightarrow 
W(R)/\hat W(\Fa_1)\twoheadrightarrow W(R/\Fa_1)$$
leads to a commutative diagram
$$\begin{CD}
H_\mu(\Fa_1)@>>>\Fg\otimes_{\BZ_p}S\\
@V{\Phi_\mu}VV@V{\varphi^\mu}VV\\
U_{^F\mu}(W(\Fa_1))@>>>\Fg\otimes_{\BZ_p}S
\end{CD},$$
where $\varphi^\mu:=\pi\otimes V^{-1}$ and 
$\pi$ is the projector that we defined in definition \ref{stressful02}. Here notice $V^{-1}$ is well-defined because 
$V$ preserves $\hat W(\Fa_1)$ and induces a bijection on $S$. 
Our element $h$ maps to an element 
$X\in\Fg\otimes_{\BZ_p}S$ satisfying: 
$$\Ad^G(w_0(U))\varphi^\mu(X)=X$$
Since this must vanish we deduce 
$h\in G(\hat W(\Fa_1))$, so that 
$h\in H_\mu(\Fa)$ for some nilpotent $\Fa\subset\Fa_1$. 
This is what we wanted to show and it achieves the proof.
\end{proof}

The part (ii) of the above result can be seen as 
a $(G,\mu)$-analog of the fact that homomorphism 
groups between nilpotent Zink-displays are torsion 
free, which is proved in \cite[Corollary 41]{zink2} 
under the additional assumption that $\sqrt{pR}$ 
is a nilpotent ideal. A straightforward manipulation 
of our proof could be used to remove the latter 
assumption, as asserted in the following lemma. 
We leave the details to the reader, as this can 
also be seen as a consequence of \cite{lau1}, 
because homomorphism groups between $p$-divisible 
groups are obviously torsion free:

\begin{lem}
Let $P$ and $P'$ be Zink-displays over a $p$-adically separated 
and complete ring $R$. If the $\pmod p$-reduction of at 
least one of $P'$ or $P^*$ (the Cartier dual of $P$) is a 
nilpotent display over $R/pR$, then $\Hom(P',P)_{\tor}=0$.
\end{lem}

We would like to conclude this circle of ideas 
with another consequence, namely a strengthening 
of \cite[Proposition 5.2.3]{pappas}:

\begin{lem}
Let $(G,\mu)$ and $(G',\mu')$ be display data over $W(k_0)$ in 
the sense of definition \ref{concept02}. Let $i:G\rightarrow G'$ 
be a closed immersion of group schemes over $\BZ_p$ such that $\mu'=i\circ\mu$ (i.e. $(G,\mu)\stackrel i{\hookrightarrow}(G',\mu')$
in the sense of \cite[subsubsection 3.7.1]{pappas}). Let $R$ be a 
$W(k_0)$-algebra satisfying $pR\subset\sqrt{0_R}$ and let $\D_1$ 
and $\D_2$ be $(G,\mu)$-displays over $R$ such that their images $i(\D_1)$ and $i(\D_2)$ are adjoint nilpotent $(G',\mu')$-displays. 
A pair consisting of an isomorphism $\psi:i(\D_1)\rightarrow i(\D_2)$ and a quasi-isogeny $\phi:\D_1\dashrightarrow\D_2$ such that 
$\psi$ induces the quasi-isogeny 
$i(\phi):i(\D_1)\dashrightarrow i(\D_2)$, comes from a 
unique isomorphism $\epsilon:\D_1\rightarrow\D_2$.
\end{lem}
\begin{proof}
It suffices to consider the banal case, so let $\psi$ 
(resp. $\D_1$ and $\D_2$) be represented by 
$H_{\mu'}(R)\ni h'=U_2\Phi_{\mu'}(h')U_1^{-1}$ (resp. 
$U_1,U_2\in G(W(R))$). Again, part (i) of theorem \ref{separ} 
allows us to ignore any nilpotent ideals, in particular we 
will assume $pR=0$ from now on. The existence of $\phi$ proves 
that $h'$ lies in the inverse image of $G(W(R)/W(R)_{\tor})$ in $G'(W(R))$. Since the closed immersion $G\hookrightarrow G'$ is defined by a finitely generated ideal, there exists an integer 
$l$ such that $h'$ lies already in the inverse image of $G(W(R/\Fa_l))$ in $G'(W(R))$, due to 
$W(R)_{\tor}=\bigcup_l W(\Fa_l)$ using the notation of lemma \ref{rigid10}. 

Towards checking that $h'$ does come from the subgroup 
$G$ we can assume $l=1$, by reduction to a minimal counterexample. Given that $0=S^2=0$ we can apply an appropriate version of the unique lifting lemma in this situation (cf. \cite[Proposition 7.4]{lau2}, 
\cite[Theorem 3.2/Remark 3.3 ]{lau3}). One obtains a lift 
$h\in G(W(R)/\hat W(\Fa_1))$ of the $\mod\Fa_1$-reduction of $h'$, which is an element of $H_\mu(R/\Fa_1)\subset G(W(R/\Fa_1))$. In fact, $h$ must 
agree with the $\mod{\hat W(\Fa_1)}$-reduction of $h'$, by the same 
unique lifting applied in the larger group $G'$. We will use again that 
the Witt coordinates of any element of $\hat W(\Fa_1)$ are contained in 
a suitable nilpotent subideal $\Fa$ of $\Fa_1$. Since we proved 
that $h'$ lies in the inverse image of 
$G(W(R)/\hat W(\Fa_1))$ in $G'(W(R))$ we deduce that 
it lies already in the inverse image of 
$G(W(R)/\hat W(\Fa))$ in $G'(W(R))$ and we are done.
\end{proof}

\subsection{$(G,\mu)$-triples}
\label{rigidgrowth}

For any ideal $\Fa$ in any $W(k_0)$-algebra 
$R$ we wish to consider the following fpqc 
sheaves of groups on $\Spec R$:
\begin{eqnarray}
\label{neighbor04}
&&\hat P_\mu^\Fa:=\pbar_\mu^\Fa\times_{\gbar^\Fa}G\\
\label{neighbor01}
&&\hat H_\mu^\Fa:=\hat P_\mu^\Fa\times_{G,w_0^G}L^+G
\end{eqnarray}
So, for any $R$-algebra $S$, in the 
previous notations, we have 
$\hat P_\mu^\Fa(S)=\Gamma_{\mu}^{S,\Fa S}$. 
We also have a canonical monomorphism 
\begin{align}
\label{CanMapHtoHatH}
H_{\mu}\times_{W(k_0)}R\into \hat H_{\mu}^{\Fa}
\end{align} 
of fpqc sheaves of groups over $\Spec R$. Suppose 
that $\Fa$ is equipped with divided powers 
$\{\gamma_k:\Fa\rightarrow\Fa\}_{k\in\BN}$. 
Observe that for any flat $R$-algebra $S$, a fundamental observation of Berthelot (see e.g. \cite[Chapter III, Lemma (1.8)]{messing1}) grants 
a unique extension of $\gamma$ to a pd-structure $\gamma_S$ on the ideal $\Fa\otimes_RS=\Fa S$. Note that in the notations of lemma 
\ref{LemExtPhi}, $\hat H_{\mu}^{\Fa S}(S)$ is nothing but the group $\Gamma_{\mu}(S,\Fa S,\gamma_{S})$, and so, we have a group homomorphism  
$$\hat H_\mu^\Fa(S)=\hat H_\mu^{\Fa S}(S)\stackrel{\Psi_\mu^{\Fa S}}{\rightarrow}L^+G(S)$$ 
Whence we obtain a mini-morphism from $\hat H_\mu^\Fa$ to $L^+G\times_{W(k_0)}R$, in the 
sense of appendix \ref{mini} (i.e., a morphism 
of functors on the full subcategory of flat $R$-algebras). Following the arguments therein (cf. \eqref{awkward01}) we construct an awkward functor: 
\begin{equation}
\label{awkward02}
\Tor_R(\hat H_\mu^\Fa)\rightarrow\Tor_R(L^+G):\,\hat Q\mapsto\hat Q*^{\Psi_\mu^\Fa}L^+G
\end{equation}
With these preparations we are now in the position to define the notion of $(G,\mu)$-triple:

\begin{defn}
\label{triple01}
Let $R$ be a $W(k_0)$-algebra and let $\Fa\subset R$ be an ideal with divided powers $\{\gamma_k:\Fa\rightarrow\Fa\}_{k\in\BN}$. 
\begin{itemize}
\item[(1)]
A $(G,\mu)$-triple over $(R,\Fa,\gamma)$ is a pair $(\hat Q,u)$, where $\hat Q$ is a $\hat H_\mu^\Fa$-torsor equipped with a $L^+G$-equivariant isomorphism 
$$u:\hat Q*^{\Psi_\mu^{\Fa}}L^+G\isoto\hat Q*L^+G$$ 
\item[(2)]
We write $\Tri_{R,\Fa,\gamma}(G,\mu)$ for the groupoid of $(G,\mu)$-triples over $(R,\Fa,\gamma)$. 
Equivalently $\Tri_{R,\Fa,\gamma}(G,\mu)$ is the groupoid rendering the diagram:

\begin{equation}
\label{poshIII}
\begin{tikzcd}
\Tri_{R,\Fa,\gamma}(G,\mu)\ar[rrr]\ar[dd]&&&\Tor_R(L^+G)\arrow[dd,"\Delta_{\Tor_R(L^+G)}"]\\\\
\Tor_R(\hat H_\mu^\Fa)\arrow[rrr, "({\Psi_\mu^\Fa}\times\id)\circ\Delta_{\Tor(\hat H_\mu^\Fa)}" below]&&&\Tor_R(L^+G)^2
\end{tikzcd}
\end{equation}

$2$-Cartesian. 
\item[(3)]
The passage from a $(G,\mu)$-triple $\hat\D$ over $(R,\Fa,\gamma)$ to its underlying 
$\hat H_\mu^\Fa$-torsor (i.e. the vertical projection in \eqref{poshIII}) will be denoted by 
$$q:\Tri_{R,\Fa,\gamma}(G,\mu)\rightarrow\Tor_R(\hat H_\mu^\Fa),\,(\hat Q,u)\mapsto\hat Q$$ 
and $\hat\D$ is called banal if $q(\hat\D)$ possesses a global section, i.e. is trivial.
\item[(4)]
The precomposition of the natural $1$-morphism $\ton^{\hat H_\mu^\Fa}\rightarrow\ton^{\hat P_\mu^\Fa}$, induced by $w_0^G:L^+G\to G$ 
with $q$ is called the \emph{lowest truncation} and will be denoted by $q_0:\Tri_{R,\Fa,\gamma}(G,\mu)\rightarrow\ton^{\hat P_\mu^\Fa}$.
\end{itemize}
\end{defn}

\begin{rem}
\label{triple02}
\begin{itemize}
\item[(a)] 
In principle, it would have been possible to define a notion of $(G,\mu)$-triple over an arbitrary $W(k_0)$-scheme which is equipped with a 
pd-ideal sheaf. We refrained from doing that, because we did not see any applications, very much unlike $(G,\mu)$-displays over non-affine bases.
\item[(b)] 
Notice that a $(G,\mu)$-triple over $(R,0_R,0)$ is nothing but a $(G,\mu)$-display over $R$. Furthermore, a $(G,\mu)$-triple over 
$(R,\Fa,\gamma)$ allows a pull-back along any pd-morphism say $(R,\Fa,\gamma)\rightarrow(S,\Fb,\delta)$, in particular any 
$(G,\mu)$-triple over $(R,\Fa,\gamma)$ determines a $(G,\mu)$-display over $R/\Fa$, which we will refer to as its \emph{canonical reduction}.
\end{itemize}
\end{rem}

\begin{defn}
\label{DefAdjNilTrip}
We call a $(G,\mu)$-triple over $(R,\Fa,\gamma)$ to be \emph{adjoint nilpotent} if the $\pmod p$ reduction of its 
canonical reduction is an adjoint nilpotent $(G,\mu)$-display over $R/(pR+\Fa)$. The class of adjoint nilpotent $(G,\mu)$-triples 
over $(R,\Fa,\gamma)$ forms a natural full subcategory, which will be denoted by $\Tri_{R,\Fa,\gamma}^{nil}(G,\mu)$
\end{defn}

\begin{rem}
\label{triple03}
Notice that the full subcategory of banal $(G,\mu)$-triples over $(R,\Fa,\gamma)$ is equivalent to the groupoid $[L^+G(R)/_{\Psi_\mu^\Fa}\hat H_\mu^\Fa(R)]$ 
(cf. appendices \ref{app1st} and \ref{app2nd}). Its full subcategory of adjoint nilpotent banal triples over $(R,\Fa,\gamma)$ is equivalent 
to $[C_{W(R)}^{nil}(G,\mu)/_{\Psi_\mu^\Fa}\hat H_\mu^\Fa(R)]$. In the sequel we will quite often assume that $\Fa$ is $p$-adically 
separated and complete. In this case a $(G,\mu)$-triple is banal if and only if its canonical reduction is banal, by lemma \ref{torsor03} (i).
\end{rem}

Unfortunately, we need another preparation for the final result of this section, namely an extension of corollary 
\ref{torsor06} to the categories of torsors under our non-representable sheaves of groups $\hat H_\mu^?$:

\begin{lem}
\label{triple04}
Let $\Fa$ be an ideal in a $W(k_0)$-algebra $R$.
\begin{itemize}
\item[(i)]
The natural $2$-commutative diagram 

\[\begin{tikzcd}
\Tor_R(\hat H_\mu^\Fa)\ar[r]\ar[d]&\Tor_R(L^+G)\ar[d]\\
\Tor_{R/\Fa}(P_\mu)\ar[r]&\Tor_{R/\Fa}(G)
\end{tikzcd}\]

is $2$-Cartesian.
\item[(ii)]
Assume that $\Fa$ is closed with respect to some admissible topology $\{\Fa_i\}_{i\in\BN_0}$ on $R$. Then we have:
$$\Tor_R(\hat H_\mu^\Fa)\isoto 2-\invlim_i\Tor(\hat H_\mu^{(\Fa_i+\Fa)/\Fa_i})$$
\end{itemize}
\end{lem}
\begin{proof}
Observe, that there are equivalences of categories: 
$$\Tor_R(\gbar^\Fa)\isoto\Tor_{R/\Fa}(G)$$ 
$$\Tor_R(\pbar_\mu^\Fa)\isoto\Tor_{R/\Fa}(P_\mu)$$ 
So the first assertion follows from lemma \ref{FuncChangGroup5} and the faithful flatness of $w_0^G:L^+G\rightarrow G$ together with $G\onto\gbar^\Fa$ being 
an epimorphism, which in turn follows from the lemmas \ref{answer04} and \ref{neighbor07}.\\

The second assertion is a straightforward consequence 
of the first one together with corollary \ref{torsor06}, 
here notice that the descending chain of ideals 
$\{(\Fa+\Fa_i)/\Fa\}_{i\in\BN_0}$ defines an admissible 
topology on $R/\Fa$.
\end{proof}

Let $R$ be a $W(k_{0})$-algebra and $(R,\Fa,\gamma)$ a 
pd-structure. Let $S$ be a flat $R$-algebra. As we said 
before, $(S,\Fa S,\gamma_S)$ is then a pd-structure 
over $(R,\Fa,\gamma)$. The pull-back of triples 
along the pd-morphism
$$(R,\Fa,\gamma)\rightarrow(S,\Fa S,\gamma_S)$$ 
deserves special attention and we will write
\begin{align*}
\Tri_{R,\Fa,\gamma}(G,\mu)&\rightarrow\Tri_{S,\Fa S,\gamma_S}(G,\mu)\\
\hat\D&\mapsto\hat\D_S
\end{align*}

where $\hat\D$ is a $(G,\mu)$-triple over $(R,\Fa,\gamma)$. 
Notice that this is well defined only for flat $R$-algebras $S$! 

\begin{thm}
\label{triple05}
Let $k_{0}$ be a finite field and $(G,\mu)$ a display datum over $W(k_0)$. Let $R$ be a $W(k_0)$-algebra, and let $\Fa\subset R$ be a $p$-adically separated and complete ideal with pd-structure $\{\gamma_k:\Fa\rightarrow\Fa\}_{k\in\BN}$, 
so that $R$ is separated and complete with respect to the topology defined by the descending chain of pd-ideals $\{p^i\Fa\}_{i\in\BN_0}$. The canonical reduction establishes an equivalence of categories from the category of adjoint nilpotent 
$(G,\mu)$-triples over $(R,\Fa,\gamma)$ to the category of $(G,\mu)$-displays over $R/\Fa$ satisfying the adjoint nilpotence condition over $R/(pR+\Fa)$.
\end{thm}
\begin{proof}
Let us first settle the banal case. Observe that $W(\Fa)$ is a $p$-adically separated and complete pd-ideal. Part (ii) of 
lemma \ref{torsor05} yields the surjectivity of $G(W(R))\rightarrow G(W(R/\Fa))$ (N.B.: The requested topological nilpotence of 
$W(\Fa)$-elements follows from the pd-property as $x^{p^n}\in p^n!\Fa$ for $x\in W(\Fa)$). By the same token we obtain the surjectivity of 
$\hat H_\mu^\Fa(R)\rightarrow H_\mu(R/\Fa)$. Moreover, even $C_{W(R)}^{nil}(G,\mu)\to C_{W(R/\Fa)}^{nil}(G,\mu)$ is surjective, as 
any lift of an $C_{W(R/\Fa)}^{nil}(G,\mu)$-element to $G(W(R))$ lies automatically in $C_{W(R)}^{nil}(G,\mu)$, according to lemma 
\ref{stressful01} together with $\Fa\subset\sqrt{p\Fa}\subset\sqrt{pR}$. Together with corollary \ref{GMZCF04} this implies the equivalence:
$$[C_{W(R)}^{nil}(G,\mu)/_{\Psi_\mu^\Fa}\hat H_\mu^\Fa(R)]\isoto[C_{W(R/\Fa)}^{nil}(G,\mu)/_{\Phi_\mu}H_\mu(R/\Fa)],$$
which is nothing but the banal case of the requested result, in view of remark \ref{triple03}.\\ 

For the proof in the non-banal case, we use descent theory to reduce to the banal case: First and foremost, notice that part (ii) of lemma \ref{triple04} 
allow us to restrict to the $p$-adically discrete case i.e. $p^i\Fa=0$ for $i\gg0$, here notice that the condition ``$p^i\Fa=0$'' is local for the fpqc topology, 
while ``$\Fa$ is $p$-adically separated and complete'' is not. The full faithfulness of the canonical reduction is straightforward: If $\D$ (resp. $\D'$) are 
the canonical reductions of adjoint nilpotent $(G,\mu)$-triples $\hat\D=(\hat Q,u)$ (resp, $\hat\D'=(\hat Q',u')$), then we choose a faithfully flat $R$-algebra 
$S$, such that $(S,\Fa\otimes_RS,\gamma_S)$ banalizes $\hat\D$ and $\hat\D'$. The result follows from the full faithfulness of the canonical reduction 
over $S$ together with descent theory. The determination of the essential image is a little trickier, because it is not so clear how to choose our trivialisation 
$S$: Starting out from a $(G,\mu)$-display $\D=(Q,u)$ over $R/\Fa$ such that $\D\times_{R/\Fa}R/pR+\Fa$ is adjoint nilpotent, we need to check the 
following auxiliary assertion: There exists a faithfully flat extension $R\rightarrow S$ such that $S/\Fa S$ banalizes $\D$. Indeed, part (ii) of lemma 
\ref{proetale01} yields a section of $Q$ over, say $\dirlim_nA_n$, where $R/\Fa\rightarrow A_1\rightarrow A_2\dots$ is a direct system with \'etale faithfully 
flat transition maps. By \cite[Th\'eor\`eme (18.1.2)]{ega4} each of the \'etale $R/\Fa$-algebras $A_n$ has a unique lift, which gives rise to a direct system 
$R\rightarrow S_1\rightarrow S_2\dots$ again with \'etale faithfully flat transition maps. The limit $S=\dirlim_nS_n$ is still faithfully flat and so it solves our 
auxiliary assertion (although $S$ may not even be finitely generated let alone \'etale). Now it is clear that there exists an adjoint nilpotent $(G,\mu)$-triple 
$\hat\D$ over $(S,\Fa S,\gamma_S)$ of which the canonical reduction agrees with $\D\times_{R/\Fa}S/\Fa S$. It remains to construct an effective descent 
datum on $\hat\D$, which may be deduced from the full faithfulness of the canonical reduction.
\end{proof}

\subsection{Lifts of $(G,\mu)$-displays}
Now we can complete our description of the set of lifts of a $(G,\mu)$-display. As before, let $R$ be a $W(k_0)$-algebra and 
$(R,\Fa,\gamma)$ a pd-structure. Our point of departure is the canonical monomorphism \eqref{CanMapHtoHatH} together 
with its truncated pendant $P_\mu\times_{W(k_0)}R\into\hat P_\mu^\Fa$. These lead to a commutative diagram
\[
\begin{tikzcd}
\B(G,\mu)(R)\ar[r]\ar[d]&\Tor_R(H_\mu)\ar[r]\ar[d]&\Tor_R(P_\mu)\ar[d]\\
\Tri_{R,\Fa,\gamma}(G,\mu)\ar[r]&\Tor_R(\hat H_\mu^\Fa)\ar[r]&\Tor_R(\hat P_\mu^\Fa)
\end{tikzcd}
\]

within which both squares are $2$-Cartesian (cf. part (i) of lemma \ref{triple04}). Fix a $\hat P_\mu^\Fa$-torsor $\hat Q_0$. Denote its 
extension of structure group via $\hat P_\Fa^\mu\rightarrow G\times_{\BZ_p}R$ by $P$. Moreover, let $Q_0$ be the $P_\mu$-torsor over 
$R/\Fa$ obtained by pulling back $\hat Q_0$ (N.B.: The pull-back of $\hat P_\mu^\Fa$ to $R/\Fa$ agrees with $P_\mu\times_{W(k_0)}R/\Fa$). 
Lemma \ref{descent06} justifies to consider the $R$-scheme $X_{\hat Q_0}:=P/P_\mu$, of which the $\mod\Fa$-reduction is pointed, as 
$X_{\hat Q_0}\times_RR/\Fa\cong Q_0*G_{W(k_0)}/P_\mu$ (N.B.: The $P_\mu$-action on $G_{W(k_0)}/P_\mu$ preserves the neutral element 
therein). It follows that the inverse image of the aforementioned base point, say $x_{\hat Q_0}\in X_{\hat Q_0}(R/\Fa)$, under the $\mod\Fa$-reduction: 
\begin{equation}
\label{poshVI}
X_{\hat Q_0}(R)\rightarrow X_{\hat Q_0}(R/\Fa)
\end{equation}
agrees with the set of $P_\mu$-structures on the $\hat P_\mu^\Fa$-torsor $\hat Q_0$. Summing up, we have proven the following:

\begin{cor}
\label{triple06}
Let $(\Fa,\gamma)$ be a $p$-adically separated and complete pd-ideal in some $W(k_0)$-algebra $R$. Let $\D$ be a $(G,\mu)$-display over $R/\Fa$ of which the $\pmod p$-reduction is 
assumed to be adjoint nilpotent. Then there exists a 
canonical bijection between the set of isomorphism classes 
of lifts of $\D$ to a $(G,\mu)$-display over $R$ and the set 
of lifts of $x_{\hat Q_0}\in X_{\hat Q_0}(R/\Fa)$ to an 
$R$-valued point of $X_{\hat Q_0}$, where $\hat Q_0$ is the 
lowest truncation of the $(G,\mu)$-triple $\hat\D$ over $(R,\Fa,\gamma)$ which is associated to $\D$ according 
to theorem \ref{triple05}.
\end{cor}

\begin{rem}
Notice that there exists at least one such lift of $\D$ as a consequence of lemma \ref{torsor03} while \eqref{poshVI} is surjective by part (ii) of lemma \ref{torsor05}.
\end{rem}

Combining theorem \ref{separ} and corollary \ref{triple06}, 
we obtain the following:

\begin{thm}
\label{compl}
Let $k_{0}$ be a finite field and $(G,\mu)$ a display datum over $W(k_0)$. Let $R$ be a $W(k_{0})$-algebra and $\Fa\subset R$ an ideal. Finally, let $\D$ a $(G,\mu)$-display over $R/\Fa$, where $\Fa$ is assumed to satisfy one of the following:
\begin{itemize}
\item[(i)]
$\Fa$ is $(pR+\Fa)$-adically separated and complete.
\item[(ii)]
$\Fa$ allows a divided power structure and is $p$-adically separated and complete.
\end{itemize}

Assume that $\D\times_{R/\Fa}R/(pR+\Fa)$ is adjoint nilpotent. Then, for each lift $\tilde\D$ of $\D$ the reduction of automorphisms gives an injection $\Aut(\tilde\D)\into\Aut(\D)$. Moreover there exists at least one such lift.
\end{thm}

\section{$(G,\mu)$-Windows}

Fix a finite field $k_0$, a smooth affine group scheme $G$ over $W(k_0)$, and a 
cocharacter $\mu:\BG_{m,W(k_0)}\rightarrow G$, as in subsection \ref{simplerversion}.

\begin{defn}
\label{wind03}
Let $P$ be a $G$-torsor over a flat $W(k_0)$-algebra $A$.
\begin{itemize}
\item[(1)]
By a \emph{modification} of $P$ we mean a pair $(P_1,\eta_1)$, where $P_1$ is a $G$-torsor over 
$A$ and \[\eta_{1}:P\times_AA[\frac1p]\to P_1\times_AA[\frac1p]\] is a $G$-equivariant isomorphism.
\item[(2)]
Let $\eta_0:Q_0*G\isoto P\times_AA/pA$ be a $G$-equivariant isomorphism, where $Q_0$ is a $P_\mu$-torsor over $A/pA$. We say that a modification 
$(P_1,\eta_1)$ is \emph{adapted} to $(Q_0,\eta_0)$ if there exist a faithfully flat $A$-algebra $B$, and elements $x_0\in Q_0(B/pB)$, $x_1\in P_1(B)$ and 
$x\in P(B)$ such that the $\pmod p$ reduction of $x$ agrees with $\eta_0(x_0)\in P(B/pB)$ while $\eta_1(x)=x_1\mu(p)$ holds in $P_1(B[\frac1p])$, i.e.:
$$P_1(B[\frac1p])\supset P_1(B)\mu(p)\ni x_1\mu(p)\stackrel{\eta_1}{\longmapsfrom}x\rightsquigarrow\eta_0(x_0)\in P(B/pB)$$
\end{itemize}
\end{defn}

\begin{defn}
\label{wind06}
If $(P_2,\eta_2)$ is another modification of $P$, we say that $(P_1,\eta_1)$ and $(P_2,\eta_2)$ are \emph{equivalent} and write $(P_1,\eta_1)\sim(P_2,\eta_2)$ 
if the $G$-equivariant isomorphism $$\eta_2\circ\eta_1^{-1}:P_1\times_AA[\frac1p]\rightarrow P_2\times_AA[\frac1p]$$ is already defined over $A$.
\end{defn}

\begin{lem}
\label{wind04} 
Let $G/W(k_0)$ and $\mu:\BG_{m,W(k_0)}\rightarrow G$ be as above, but assume that $-1$ is the 
only negative $\mu$-weight of $\Lie G$. Let $P$ be a $G$-torsor over a flat $W(k_0)$-algebra $A$. 
\begin{itemize}
\item[(i)]
Up to equivalence, there exists at most one modification of $P$ which is adapted to a given $P_\mu$-structure $(Q_0,\eta_0)$ on $P\times_AA/pA$.
\item[(ii)]
Let $\eta_1:P\times_AA[\frac1p]\stackrel{\cong}{\rightarrow}P_1\times_AA[\frac1p]$ be a modification 
adapted to a $P_\mu$-structure $(Q_0,\eta_0)$ on $P\times_AA/pA$. If $A$ is $p$-adically separated and 
complete and $x_0$ is a global section of $Q_0$, then there exist $x_1\in P_1(A)$ and $x\in P(A)$ satisfying
$$P_1(A[\frac1p])\ni x_1\mu(p)\stackrel{\eta_1}{\longmapsfrom}x\rightsquigarrow\eta_0(x_0)\in P(A/pA),$$
in particular $P_1$ and $P$ are trivial too.
\end{itemize}
\end{lem}
\begin{proof}
We start with assertion (i), so let $(P_1,\eta_1)$ and $(P'_1,\eta'_1)$ be modifications of $P$ that are adapted to the same 
$P_\mu$-structure $(Q_0,\eta_0)$ on the $\pmod p$ reduction $P\times_AA/pA$. Over a sufficiently large faithfully flat extension 
$A\rightarrow B$ there exist elements $x_1\in P_1(B)$, $x_1'\in P_1'(B)$, and $x,x'\in P(B)$ and $x_0,x_0'\in Q_0(B/pB)$ satisfying
\begin{eqnarray*}&&P_1(B[\frac1p])\ni x_1\mu(p)\stackrel{\eta_1}{\longmapsfrom}x\rightsquigarrow\eta_0(x_0)\in P(B/pB)\\
&&P_1'(B[\frac1p])\ni x_1'\mu(p)\stackrel{\eta_1'}{\longmapsfrom}x'\rightsquigarrow\eta_0(x_0')\in P(B/pB)\end{eqnarray*}
according to the definition. There is a unique element $g\in G(B)$ with $x'g=x$ and likewise one and only one $g_0\in P_\mu(B/pB)$ with $x_0'g_0=x_0$.
Observe that $g_0$ agrees with the reduction of $g$ modulo $p$, as $\eta_0$ sends $x_0$ (resp. $x'_0$) to the reduction of $x$ (resp. $x'$) modulo $p$. Also, 
observe that the isomorphism $\eta:=\eta_1'\circ\eta_1^{-1}$ from the $G$-torsor $P_1\times_AA[\frac1p]$ to $P_1'\times_AA[\frac1p]$ sends $x_1\in P_1(B)$ 
to $x_1'\mu(p)g\mu(\frac1p)\in P_1'(B[\frac1p])$. According to corollary \ref{factd} we have $\mu(p)g\mu(\frac1p)\in G(B)$, as 
$$G(B)\ni g\rightsquigarrow g_0\in P_\mu(B/pB)$$
and recall that $-1$ is the only negative $\mu$-weight of $\Lie G$. Whence it follows $x_1'\mu(p)g\mu(\frac1p)\in P_1'(B)$, so that 
$\eta$ is defined over $B$, and afortiori over $A$ too. We deduce $(P_1,\eta_1)\sim(P'_1,\eta'_1)$, in the sense of definition \ref{wind06}.\\

Towards the assertion (ii), we observe that part (ii) of lemma \ref{torsor05} yields a section $x\in P(A)$ lifting $\eta_0(x_0)\in P(A/pA)$, so that $P$ is trivial. 
Consider the $G$-equivariant isomorphism $\eta_1'$ from $P\times_AA[\frac1p]$ to $G\times_AA[\frac1p]$ which sends the global section $x$ to $\mu(p)$. 
The modification $(G,\eta_1')$ is adapted to the $P_\mu$-structure $(Q_0,\eta_0)$ as can be seen by taking $B=A$ and $x_1'=1$ in definition \ref{wind03}. 
By part (i), we have $(P_1,\eta_1)\sim(G,\eta_1')$ and therefore the integral section $1=x_{1}'\in G(A)$ would yield the desired element $x_{1}\in P_{1}(A)$.
\end{proof}

\subsection{Frobenius $G$-torsors}

Suppose that $\sigma$ is a Frobenius lift on a torsion-free $p$-adically separated and complete ring $A$. For any $n\in\BN$, the pull-back along 
$\sigma^n:A\rightarrow A$ will be denoted with a superscripted ``$\sigma^n$'', e.g. $^{\sigma^n}M:=M\otimes_{\sigma^n,A}A$ for any $A$-module $M$ and 
$^{\sigma^n}X:=X\times_{A,\sigma^n}A$ for any $A$-scheme $X$. It is folklore to call a pair $(M,\phi)$, where $M$ is a finitely generated projective $A$-module and 
$$\phi:{^\sigma M}\otimes_AA[\frac1p]\cong M\otimes_AA[\frac1p]$$ 
is an isomorphism, a \emph{Frobenius module}. We would like to spell out the analog in the setting of $G$-torsors:

\begin{defn}
\label{wind01}
Let $\sigma:A\rightarrow A$ be as above and let $G$ be a smooth affine $\BZ_p$-group. By a \emph{Frobenius $G$-torsor} 
over $(A,\sigma)$ we mean a pair $(P,b)$, where $P$ is a $G$-torsor over $A$ and $b$ is a $G$-equivariant isomorphism 
from $^\sigma P\times_AA[\frac1p]$ to $P\times_AA[\frac1p]$. If $(P,b)$ is a Frobenius $G$-torsor over $(A,\sigma)$ we write 
$P^{\Ad}$ for the finitely generated projective $A$-module correspoding to the $\GL(\Fg)$-torsor $P*^{\Ad^G}\GL(\Fg)$ and 
$$b^{\Ad}:{^\sigma P}^{\Ad}\otimes_{A}A[\frac1p]\rightarrow P^{\Ad}\otimes_{A}A[\frac1p]$$ 
for the isomorphism induced by $b$. A Frobenius $G$-torsor $(P,b)$ over $(A,\sigma)$ is called 
\emph{adjoint-nilpotent} if one (and hence both) of the following equivalent properties is fulfilled:
\begin{itemize}
\item
There exists $n\in\BN$ such that the composition $p^n\cdot b^{\Ad}\circ\sigma(b^{\Ad})\circ\dots\circ\sigma^n(b^{\Ad})$ 
is integral (i.e. sends $^{\sigma^{n+1}}P^{\Ad}$ to $P^{\Ad}$)
\item
The $\sigma$-linear endomorphism 
$p\cdot b^{\Ad}\circ(\id_{P^{\Ad}}\otimes\sigma)$ is $p$-adically topologically nilpotent on $P^{\Ad}\otimes_{A}A[\frac1p]$
\end{itemize}
We call $(P,b)$ \emph{banal} if $P$ is a trivial $G$-torsor, i.e. $P(A)\neq\emptyset$. A morphism between Frobenius $G$-torsors $(P_1,b_1)$ and $(P_2,b_2)$ 
is a morphism $\alpha:P_1\rightarrow P_2$ of their underlying torsors, such that $b_2\circ\sigma(\alpha)\circ b_1^{-1}=\alpha$ holds over $A[\frac1p]$. 
\end{defn}

It is clear that Frobenius $G$-torsors over $(A,\sigma)$ form a groupoid which we denote by $\B_{A,\sigma}(G)$. Its full subcategory consisting of adjoint 
nilpotent Frobenius $G$-torsors over $(A,\sigma)$ will be denoted by $\B_{A,\sigma}^{nil}(G)$. In either case the full subcategories of banal objects read:
\begin{eqnarray}
&&[G(A[\frac1p])/_\sigma G(A)]\into \B_{A,\sigma}(G)\\
&&[C_A^{nil}(G)/_\sigma G(A)]\into \B_{A,\sigma}^{nil}(G)
\end{eqnarray}
where $$C_A^{nil}(G):=\{b\in G(A[\frac1p])\,|\,\exists n\in\BN:\,p^n\Ad^G(b\sigma(b)\dots\sigma^n(b))(\Fg)\subset A\otimes_{\BZ_p}\Fg\}$$ 
equivalently, $C_A^{nil}(G)$ is the subset of $G(A[\frac1p])$ consisting of elements $b$ such that $p^n\big(\Ad^G(b)\circ(\sigma\otimes\id_\Fg)\big)^n$ 
converges $p$-adically towards $0$. Observe that we have an equality (cf. definition \ref{stressful02}):
\[G(A)\cap C_A^{nil}(G){^F\mu(\frac1p)}=C_A^{nil}(G,\mu)\]
In order to see this, consider $b=U{^F\mu(p)}$ for some $U\in G(A)$, and observe that the endomorphism 
$p\Ad^G(b)\circ(\sigma\otimes\id_\Fg)$ preserves $A\otimes_{\BZ_p}\Fg$ and that its $\pmod p$-reduction is given by 
the formula \eqref{stressful04}. Hence $U$ satisfies the adjoint nilpotence condition if and only if $b$ is contained in $C_A^{nil}(G)$.

\subsection{$(G,\mu)$-Windows over Frames}

From now on we fix a finite field $k_0\supset\BF_p$ and a display datum $(G,\mu)$ over $W(k_0)$ (cf. definition \ref{concept02}).

\begin{defn}
\label{wind07}
By a $W(k_0)$-frame we mean a triple $(A,\Fa,\sigma)$ where $\Fa$ is a $p$-adically closed pd-ideal in a 
torsion-free $p$-adically separated and complete $W(k_0)$-algebra $A$ on which $\sigma$ is a Frobenius lift.
\end{defn}

Recall, that pd-structures on ideals in torsion-free rings are unique. Furthermore, given the existence of a lift $\sigma$ 
of the absolute Frobenius, we observe that the existence of a pd-structure on $\Fa$ implies a crucial inclusion
\begin{equation}
\label{wind09}
\Fa\subset\sigma^{-1}(pA)
\end{equation}
as $\sigma^{-1}(pA)$ is in fact the largest pd-ideal of $A$. When working with one $W(k_0)$-frame $(A,\Fa,\sigma)$ at a time, we would write
\begin{equation}
\label{wind05}
\sigma_0:A/\Fa\rightarrow A/pA.
\end{equation}
for the homomorphism induced by $\sigma$. Pull-back along this morphism will be indicated by superscripted ``$\sigma_0$'', in particular we let 
$^{\sigma_0}X:=X\times_{A/\Fa,\sigma_0}A/pA$ for any $A/\Fa$-scheme $X$. In the same vein, observe that a $P_\mu$-structure $(Q_0,\eta_0)$ on the 
$\mod\Fa$ reduction of some $G$-torsor $P$ over $A$ gives rise (by pulling-back along $\sigma_{0}$) to a $P_{^F\mu}$-structure on the $\pmod p$ reduction of $^\sigma P$, which we denote by $({^{\sigma_0}Q_0},\sigma_0(\eta_0))$.\\

In the sequel we write $\hat\Phi_\mu^{A,\sigma}:\Gamma_\mu^{A,\sigma^{-1}(pA)}\rightarrow G(A)$ for the composition of the following group homomorphisms:

$$\Gamma_\mu^{A,\sigma^{-1}(pA)}\stackrel{\sigma}{\longrightarrow}\Gamma_{^F\mu}^{A,pA}\stackrel{\intaut({^F\mu}(p))}{\longrightarrow}G(A)$$

Notice that the second map is well-defined, because of corollary \ref{factd}. Clearly one has 
$\hat\Phi_\mu^{A,\sigma}(g)=\sigma(\mu(p)g\mu(\frac1p))={^F\mu(p)}\sigma(g){^F\mu(\frac1p)}$. Observe that 
$[C_A^{nil}(G,\mu)/_{\hat\Phi_\mu^{A,\sigma}}\Gamma_\mu^{A,\Fa}]$ agrees with $[C_A^{nil}(G,\mu){^F\mu(p)}/_\sigma\Gamma^{A,\Fa}_\mu]$. 
We would like to add some general comments on morphisms between frames:

\begin{defn}
\label{wind08}
Fix a finite field $k_0$ and let $(A,\Fa,\sigma)$ and $(B,\Fb,\tau)$ be $W(k_0)$-frames. By a morphism $\varkappa$ from 
the former to the latter we mean a $W(k_0)$-algebra homomorphism $A\rightarrow B$ such that $\varkappa(\Fa)\subset\Fb$ 
and $\varkappa\circ\sigma=\tau\circ\varkappa$. At last we call $\varkappa$ a $p$-adic covering if and only if its special fiber 
$A/pA\rightarrow B/pB$ is faithfully flat while $\Fb$ is topologically generated by $\varkappa(\Fa)$ with respect to the $p$-adic topology.
\end{defn}

For any $p$-adic covering $\varkappa:(A,\Fa,\sigma)\rightarrow(B,\Fb,\tau)$ we claim that 
the maps $A/p^nA\rightarrow B/p^nB$ are flat for every $n$, in fact we have the following:

\begin{lem}
\label{descent03}
Let $\Lambda$ be a ring, $M\subset\Lambda$ an ideal and $\varkappa:A\rightarrow B$  a morphism of flat $\Lambda$-algebras. If the 
induced morphism $A/MA\rightarrow B/MB$, is flat then the same holds for $A/M^nA\rightarrow B/M^nB$, for any $n\in\BN$.
\end{lem}
\begin{proof}
Let us prove this by induction on $n$, so given an ideal $I\subset A/M^nA$, we need to show the injectivity of 
$I\otimes_{A/M^nA}B/M^nB\rightarrow B/M^nB$. Consider the ideal $I\cap MA/M^nA=:J\subset MA/M^nA$ and the commutative diagram

\[
\begin{tikzcd}
J\otimes_{A/M^nA}B/M^nB\ar[r]\ar[d]&I\otimes_{A/M^nA}B/M^nB\arrow[r]\ar[d]&I/J\otimes_{A/MA}B/MB\ar[d]\\
MB/M^nB\ar[r]&B/M^nB\ar[r]&B/MB
\end{tikzcd}
\]

Notice the identifications
\begin{eqnarray*}
&&MA/M^nA\otimes_{A/M^nA}B/M^nB\cong MB/M^nB\\
&&J\otimes_{A/M^{n-1}A}B/M^{n-1}B\cong J\otimes_{A/M^nA}B/M^nB\\
&&I/J\otimes_{A/M^nA}B/M^nB\cong I/J\otimes_{A/MA}B/MB
\end{eqnarray*}
The map $$J\otimes_{A/M^{n-1}A}B/M^{n-1}B\rightarrow B/M^{n-1}B$$ is injective by the induction hypothesis. The map 
$$I/J\otimes_{A/MA}B/MB\rightarrow B/MB$$ 
is injective by the assumed flatness of $B/MB$ over $A/MA$. At last observe the exactness of the sequence 
$$J\otimes_{A/M^nA}B/M^nB\rightarrow I\otimes_{A/M^nA}B/M^nB\rightarrow I/J\otimes_{A/M^nA}B/M^nB$$
so that a straightforward diagram chase yields the result.
\end{proof}

\begin{defn}
\label{wind02}
Let $(A,\Fa,\sigma)$ be a $W(k_0)$-frame. Consider a Frobenius $G$-torsor $(P,b)$ over $(A,\sigma)$ 
together with a $P_\mu$-structure $(Q_0,\eta_0)$ on $P_0:=P\times_AA/\Fa$, in other words, we have:
\begin{itemize}
\item
$P$ is a $G$-torsor over $A$.
\item
$b$ is a $G$-equivariant isomorphism from $^\sigma P\times_AA[\frac1p]$ to $P\times_AA[\frac1p]$.
\item
$Q_0$ is a $P_\mu$-torsor over $A/\Fa$. 
\item
$\eta_0$ is a $G$-equivariant isomorphism from $Q_0*G$ to $P_0$. 
\end{itemize}
We call the quadruple $(P,Q_0,\eta_0,b)$ a \emph{$(G,\mu)$-window} over $(A,\Fa,\sigma)$ if and only if the following axioms hold: 
\begin{itemize}
\item[(i)]
The modification of $^\sigma P$ given by $b:{^\sigma P}\times_AA[\frac1p]\isoto P\times_AA[\frac1p]$ is adapted to the $P_{^F\mu}$-structure 
$({^{\sigma_0}Q_0},\sigma_0(\eta_0))$ on $^\sigma P\times_AA/pA$, which is the pull-back of $(Q_0,\eta_0)$ along the morphism $\sigma_0:A/\Fa\rightarrow A/pA$.
\item[(ii)]
The Frobenius $G$-torsor $(P,b)$ over $(A,\sigma)$ is adjoint nilpotent, in the sense of definition \ref{wind01}.
\end{itemize}
A $(G,\mu)$-window $(P,Q_0,\eta_0,b)$ over $(A,\Fa,\sigma)$ is called \emph{banal} if and only if $Q_0$ is the trivial $P_\mu$-torsor over $A/\Fa$.
\end{defn}

It is clear that $(G,\mu)$-windows over $(A,\Fa,\sigma)$ form a groupoid, which we denote by $\B_{A,\Fa,\sigma}^{nil}(G,\mu)$. Note that this is a full subcategory 
of the groupoid of quadruples $(P,Q_0,\eta_0,b)$ satisfying (i) only, which shall be denoted by $\B_{A,\Fa,\sigma}(G,\mu)$. Let us elucidate their full 
subcategories of banal quadruples $\W=(P,Q_0,\eta_0,b)$: Given the triviality of $Q_0$ we may choose isomorphisms $P_\mu\times_{W(k_0)}A/\Fa\cong Q_0$ 
and $G\times_{\BZ_p}A\cong P$, of whose the latter is a lift of $G\times_{\BZ_p}A/\Fa\cong Q_0*G\stackrel{\eta_0}{\rightarrow}P\times_AA/pA$ 
(cf. part (i) of lemma \ref{torsor03}). This allows to rewrite $\W$ in the form $(G\times_{\BZ_p}A,P_{\mu}\times_{W(k_0)}A/pA,\id,b)$ where 
$b\in\Hom_G({^\sigma P}\times_A[\frac1p],P\times_AA[\frac1p])=G(A[\frac1p])$, as $^\sigma P=G\times_{W(k_0)}A=P$. It remains to 
analyze which elements of $G(A[\frac1p])$ satisfy the adaptedness condition (i) of definition \ref{wind02}. So $P$ is an adapted modification of $^\sigma P$ if and 
only if there exists $x_0\in P_{^F\mu}(A/pA)$ and $x,\;x_1\in G(A)$ such that the $\pmod p$ reduction of $x$ agrees with $x_0$ while $bx=x_1{^F\mu}(p)$. In view of 
corollary \ref{factd} this means that $b$ can be written as $U{^F\mu}(p)$ for some $U\in G(A)$. At last, it is easy to see that the set of
$\B_{A,\Fa,\sigma}(G,\mu)$-morphisms from another $\B_{A,\Fa,\sigma}(G,\mu)$-object of the form 
$\W'=(G\times_{\BZ_p}A,P_{\mu}\times_{W(k_0)}A/pA,\id,b')$ to $\W$ is given by 
$$\Hom(\W',\W)=\{k\in\Gamma_\mu^{A,\sigma}\mid\,b'=k^{-1}b\sigma(k)\}=\{k\in\Gamma_\mu^{A,\sigma}\mid\,U'=k^{-1}U\hat\Phi_\mu^{A,\sigma}(k)\},$$ 
where $b'=U'{^F\mu}(p)$. In the language of part \ref{app1st} of the appendix this means that the full subcategory 
of banal $(G,\mu)$-windows over $(A,\Fa,\sigma)$ (resp. banal $\B_{A,\Fa,\sigma}(G,\mu)$-objects) is equivalent to 
$[C_A^{nil}(G,\mu){^F\mu(p)}/_\sigma\Gamma_\mu^{A,\Fa}]\cong[C_A^{nil}(G,\mu)/_{\hat\Phi_\mu^{A,\sigma}}\Gamma_\mu^{A,\Fa}]$
(resp. $[G(A){^F\mu(p)}/_\sigma\Gamma_\mu^{A,\Fa}]\cong[G(A)/_{\hat\Phi_\mu^{A,\sigma}}\Gamma_\mu^{A,\Fa}]$). 
Note also that there is a natural commutative diagram of groupoids:
\[
\begin{tikzcd}
\B_{A,\Fa,\sigma}^{nil}(G,\mu)\arrow[r]\ar[d]&\B_{A,\Fa,\sigma}(G,\mu)\ar[d]\\
\B_{A,\sigma}^{nil}(G)\arrow[r]&\B_{A,\sigma}(G)
\end{tikzcd}
\]
of which the banal pendant is: 
\[
\begin{tikzcd}
{[}C_A^{nil}(G,\mu){^F\mu(p)}/_\sigma\Gamma_\mu^{A,\Fa}{]}\ar[r]\ar[d] & {[}G(A){^F\mu(p)}/_\sigma\Gamma_\mu^{A,\Fa}{]}\ar[d]\\
{[}C_A^{nil}(G)/_\sigma G(A){]}\ar[r]&{[}G(A[\frac1p])/_\sigma G(A){]}
\end{tikzcd}
\]

\begin{rem}
We note some elementary properties of these diagrams:
\begin{itemize}
\item[(a)]
All horizontal arrows in these diagrams are fully faithful.
\item[(b)]
Both of these diagrams are $2$-Cartesian.
\item[(c)]
All vertical arrows in these diagrams are faithful.
\end{itemize}
\end{rem}

\subsection{Descent}
\label{subsectionDescent}

Let $k_0$ and $(G,\mu)$ be as above and let $\varkappa:(A,\Fa,\sigma)\rightarrow(B,\Fb,\tau)$ be $W(k_0)$-linear $p$-adic 
covering. Consider the $p$-adically completed $m$-fold tensor product of $B$, i.e. the $A$-algebra
$$B_{[m]}:=B^{\hat\otimes_Am}=\widehat{B^{\otimes_Am}},$$
and let us write $\Fb_{[m]}$ for the $p$-adically closed $B_{[m]}$-ideal which is topologically generated by the set 
$\bigcup_{k+l=m-1}B^{\otimes_Ak}\otimes_A\Fb\otimes_A B^{\otimes_Al}$, so that $B_{[m]}/\Fb_m$ is isomorphic to the 
$p$-adically completed $m$-fold tensor product of the $A/\Fa$-algebra $B/\Fb$. Let us write $\tau_{[m]}:B_{[m]}\rightarrow B_{[m]}$ for the 
endomorphism induced by $x_1\otimes\dots\otimes x_m\mapsto\tau(x_1)\otimes\dots\otimes\tau(x_{m})$. With the usual notations 
for $\varkappa_{12},\varkappa_{23},\varkappa_{13}:B_{[2]}\rightarrow B_{[3]}$ (resp. $\varkappa_1,\varkappa_2:B\rightarrow B_{[2]}$) 
we obtain further morphism $(B_{[2]},\Fb_{[2]},\tau_{[2]})\rightarrow(B_{[3]},\Fb_{[3]},\tau_{[3]})$ (resp. $(B,\Fb,\tau)\rightarrow(B_{[2]},\Fb_{[2]},\tau_{[2]})$). 
In order to give a sober formulation of descent we need the following notions: 

\begin{itemize}
\item
Let $\W$ be a $(G,\mu)$-window over $(B,\Fb,\tau)$. By a $\varkappa$-descent datum on $\W$ we mean a $\B_{B_{[2]},\Fb_{[2]},\tau_{[2]}}^{nil}(G,\mu)$-morphism 
$\alpha:\varkappa_1(\W)\rightarrow\varkappa_2(\W)$. The descent datum is said to satisfy the cocycle condition if 
$\varkappa_{23}(\alpha)\circ\varkappa_{12}(\alpha)=\varkappa_{13}(\alpha)$ holds in $\B_{B_{[3]},\Fb_{[3]},\tau_{[3]}}^{nil}(G,\mu)$.
\item 
Let $\alpha'$ be another $\varkappa$-descent datum on another $\B_{B,\Fb,\tau}^{nil}(G,\mu)$-object $\W'$. By an isomorphism from $(\W,\alpha)$ to $(\W',\alpha')$,
 we mean a $\B_{B,\Fb,\tau}^{nil}(G,\mu)$-morphism $\eta:\W\rightarrow\W'$ satisfying $\alpha'\circ\varkappa_1(\eta)=\varkappa_2(\eta)\circ\alpha$.    
\end{itemize}

\begin{prop}
\label{descent01}
Let $k_0$ and $(G,\mu)$ be as above and let $\varkappa:(A,\Fa,\sigma)\rightarrow(B,\Fb,\tau)$ be a $W(k_0)$-linear $p$-adic covering. 
Then $\varkappa$ induces a fully faithful functor from $\B_{A,\Fa,\sigma}^{nil}(G,\mu)$ to the groupoid of $(G,\mu)$-windows over $(B,\Fb,\tau)$ 
with $\varkappa$-descent datum. The essential image of this functor consists of all pairs $(\W,\alpha)$, where $\alpha$ satisfies the cocycle condition.
\end{prop}
\begin{proof}
In order to descend the underlying $G$-torsor $P$ (resp. the $P_\mu$-structure $(Q_0,\eta_0)$ on its $\mod\Fa$-reduction) one simply combines corollary \ref{torsor06} with the usual descent along the maps $A/p^nA\rightarrow B/p^nB$ (resp. $A/(p^nA+\Fa)\rightarrow B/(p^nB+\Fb)$). Due to lemma \ref{descent03} these maps are faithfully flat, so that one has effective descent. In order to descend the Frobenius of the window i.e. the datum ``$b$'' in the quadruple of definiton \ref{wind02} we use the reformulation in terms of the ``$u$'' given in corollary \ref{poshI}, together with the aforementioned faithful flatness of $A/p^nA\rightarrow B/p^nB$. Descending the applicable adjoint nilpotence condition is straightforward.
\end{proof}

\subsection{A reformulation}

In this subsection we give a technical reformulation of the concept of $(G,\mu)$-window over a $W(k_0)$-frame $(A,\Fa,\sigma)$. This 
is based on the formalism of pull-back along so-called mini-morphism, as developed in part \ref{mini} of the appendix: Let us write 
$\hat\intaut_\mu$ for the mini-morphism from $\hat P_\mu^{pW(k_0)}$ to $G_{W(k_0)}$ whose effect on any flat $W(k_0)$-algebra $R$ is given by 
$$\hat P_\mu^{pW(k_0)}(R)=\Gamma_\mu^{R,pR}\rightarrow G(R);\,g\mapsto \mu(p)g\mu(\frac1p)$$
(see corollary \ref{factd}). Following the formalism \eqref{awkward01} this induces a functor
\begin{equation}
\label{wind10}
\text{mo}_\mu:\Tor_R(\hat P_\mu^{pW(k_0)})\rightarrow\Tor_R(G);\,\hat P\mapsto\hat P*^{\hat\intaut_\mu}G
\end{equation} 
for any such $R$. Observe also that for any $\hat P_\mu^{pW(k_0)}$-torsor $\hat P$ over $R$ there is a canonical isomorphism between the generic fibers of $\hat P*G$ 
and $\hat P*^{\hat\intaut_\mu}G$, i.e. the latter is a canonical modification of the former in the sense of definition \ref{wind03}. We have to start with another result:

\begin{lem}
\label{torsor04}
Let $A$ be a flat $W(k_0)$-algebra.
\begin{itemize}
\item[(i)] 
The natural $2$-commutative diagram
\[\begin{tikzcd}
\Tor_A(\hat P_\mu^{pW(k_0)})\ar[r]\ar[d]&\Tor_A(G)\ar[d]\\\Tor_{A/pA}(P_\mu)\ar[r]&\Tor_{A/pA}(G)
\end{tikzcd}\]
is 2-Cartesian, where the top horizontal map is induced by change of group along the natural inclusion $\hat P_\mu^{pW(k_0)}\into G$.
\item[(ii)]
Suppose that $\hat P$ is the $\hat P_\mu^{pW(k_0)}$-torsor that corresponds to a $P_\mu$-structure $(Q_0,\eta_0)$ on the $\pmod p$-reduction 
of some $G$-torsor $P$ over $A$, as in the diagram above. Then the modification $P_1:=\hat P*^{\hat\intaut_\mu}G$ is adapted to $(Q_0,\eta_0)$. 
\end{itemize}
In particular, adapted modifications always exist (and they are unique by part (i) of lemma \ref{wind04}).
\end{lem}
\begin{proof}
Part (i) is a canonical consequence of the criterion (i) of lemma \ref{FuncChangGroup5}. In order to deal with 
part (ii) we choose a section $\hat x\in\hat P(B)$, where $B$ is a faithfully flat $A$-algebra. This gives rise to sections $x\in P(B)$ and $x_0\in Q_0(B/pB)$ satisfying $x\rightsquigarrow\eta_0(x_0)\in P(B/pB)$. Moreover $P_1$ has a section over $B$ too, namely $x_1:=\hat x*1$. At last the canonical homomorphism from $P\times_AA[\frac1p]$ to $P_1\times_AA[\frac1p]$ sends $x$ to $x_1\mu(p)$, as can be seen immediately from the formula \eqref{FuncChangGroup6}. 
\end{proof}

We are now in the position to give the reformulation of definition \ref{wind02}, as promised at the beginning of this subsection. Let $(A,\Fa,\sigma)$ be a 
$W(k_0)$-frame and let us write $^FX$ (resp. $^\sigma X$) for the pull-back of an fpqc sheaf $X$ on $\Spec W(k_0)$ (resp. $\Spec A$) along the morphism 
$\Spec W(k_0)\stackrel{F}{\rightarrow}\Spec W(k_0)$ (resp. $\Spec A\stackrel{\sigma}{\rightarrow}\Spec A$). In particular, observe that the pull-back 
$^\sigma\hat P_\mu^\Fa$ is naturally isomorphic to $\hat P_{^F\mu}^{\sigma(\Fa)A}$, as $^FP_\mu\cong P_{^F\mu}$. Likewise, pulling back some 
$\hat P_\mu^\Fa$-torsor $\hat Q$ on $\Spec A$ yields a $\hat P_{^F\mu}^{\sigma(\Fa)A}$-torsor $^\sigma\hat Q$ on $\Spec A$. The crucial inclusion 
\eqref{wind09} allows to precompose with an extension of structure groups $\hat P_{^F\mu}^{\sigma(\Fa)A}\subset\hat P_{^F\mu}^{pA}$, thus arriving at a functor:
\begin{equation}
\label{wind14}
\Tor_A(\hat P_\mu^\Fa)\rightarrow\Tor_A(\hat P_{^F\mu}^{pA});\,\hat Q\mapsto{^\sigma\hat Q}*\hat P_{^F\mu}^{pA}
\end{equation}
Composing the functor $\text{mo}_{^F\mu}$, as defined in \eqref{wind10}, with \eqref{wind14} yields:
\begin{equation}
\label{wind11}
\text{mo}_\mu^\sigma:\Tor_A(\hat P_\mu^\Fa)\rightarrow\Tor_A(G);\,\hat Q\mapsto({^\sigma\hat Q}*\hat P_{^F\mu}^{pA})*^{\hat\intaut_{^F\mu}}G
\end{equation}

\begin{cor}
\label{poshI}
Let $(A,\Fa,\sigma)$ be a $W(k_0)$-frame and $(G,\mu)$ a $W(k_0)$-display datum, where $k_0$ is a finite field.
\begin{itemize}
\item[(i)]
The natural category structure on pairs $(\hat Q,u)$ consisting of a $\hat P_\mu^\Fa$-torsor $\hat Q$ together with a $G$-equivariant 
isomorphism $u:\text{mo}_\mu^\sigma(\hat Q)\isoto\hat Q*G$ is naturally equivalent to the groupoid $\B_{A,\Fa,\sigma}(G,\mu)$.
\item[(ii)]
Its full subcategory consisting of those pairs $(\hat Q,u)$ such that its scalar extension to $A[\frac1p]$ satisfies the adjoint nilpotence condition (cf. definition 
\ref{wind01}) is equivalent to the groupoid $\B_{A,\Fa,\sigma}^{nil}(G,\mu)$ of $(G,\mu)$-windows over $(A,\Fa,\sigma)$.
\end{itemize}
\end{cor}
\begin{proof}
Our point of departure is part (i) of lemma \ref{torsor04}, which says that giving a $\hat P_\mu^\Fa$-torsor $\hat Q$ is equivalent to giving a $G$-torsor $P$ 
over $A$ allowing a $P_\mu$-structure $(Q_0,\eta_0)$ on its $\mod\Fa$-reduction, i.e. a $G$-isomorphism $\eta_0:Q_0*^{P_\mu}G\isoto P\times_AA/\Fa$.
In order to turn $(P,Q_0,\eta_0)$ into a $\B_{A,\Fa,\sigma}(G,\mu)$-object $(P,Q_0,\eta_0,b)$ we need to give a $G$-isomorphism 
$b:^\sigma P\times_AA[\frac1p]\isoto P\times_AA[\frac1p]$ with respect to which $P$ becomes a $({^{\sigma_0}Q_0},\sigma_0(\eta_0))$-adapted modification of 
$^\sigma P$. Applying part (ii) of lemma \ref{torsor04} to the $\hat P_{^F\mu}^{pA}$-torsor ${^\sigma\hat Q}*\hat P_{^F\mu}^{pA}$ reveals that there is a unique such 
$G$-torsor, namely $\text{mo}_\mu^\sigma(\hat Q)$. We deduce that the datum ``$b$'' is equivalent to giving an isomorphism $u$ from $\text{mo}_\mu^\sigma(\hat Q)$ 
to $P$, which is $\hat Q*G$.
\end{proof}

\begin{rem}
One might want to think of this as yet another $2$-Cartesian diagram, namely:

\begin{equation}
\begin{tikzcd}
\B_{A,\Fa,\sigma}(G,\mu)\ar[rr]\ar[d]&&\Tor_A(G)\arrow[d,"\Delta_{\Tor_A(G)}"]\\
\Tor_A(\hat P_\mu^\Fa)\arrow[rr,"(\text{mo}_{\mu}^{\sigma}{,} \id)" below]&&\Tor_A(G)^2
\end{tikzcd}
\end{equation}

\end{rem}

\section{The Main Theorem}
\label{main03}
The strategy of the proof of theorem \ref{main01} consists of passing from a given $W(k_0)$-frame $(A,\Fa,\sigma)$ to the so-called Witt-frame, i.e. 
$(W(A),W(\Fa)+I(A),F)$. In the banal case this is an intricate application of corollary \ref{GMZCF04} and in the non-banal case one also has to worry 
about ``banalization'' via passage to other frames, supplemented by a study of descent properties. A last and final ingredient is the equivalence
\begin{equation}
\label{wind12}
\B_{W(A),W(\Fa)+I(A),F}^{nil}(G,\mu)\isoto\Tri_{A,\Fa,\gamma}^{nil}(G,\mu),
\end{equation}
of which the banal pendant is a tautology, namely:
\begin{equation}
\label{wind13}
[C_{W(A)}^{nil}(G,\mu)/_{\hat\Phi_\mu^{W(A),F}}\Gamma_\mu^{W(A),W(\Fa)+I(A)}]\isoto[C_{W(A)}^{nil}(G,\mu)/_{\Psi_\mu^\Fa}\hat H_\mu^\Fa(A)]
\end{equation}

Every morphism of $W(k_0)$-frames $\varkappa:(A,\Fa,\sigma)\rightarrow(B,\Fb,\tau)$ induces ``compatible'' maps

\begin{eqnarray}
\label{independence04}
&&\Gamma_\mu^{A,\Fa}\rightarrow\Gamma_\mu^{B,\Fb}\\
\label{independence03}
&&C_A^{nil}(G,\mu)\rightarrow C_B^{nil}(G,\mu)
\end{eqnarray}
along with functors
\begin{eqnarray}
\label{independence02}
&&[C_A^{nil}(G,\mu){^F\mu(p)}/_\sigma\Gamma_\mu^{A,\Fa}]
\rightarrow[C_B^{nil}(G,\mu){^F\mu(p)}/_\tau\Gamma_\mu^{B,\Fb}]\\
\label{independence01}
&&\B_{A,\Fa,\sigma}^{nil}(G,\mu)
\rightarrow \B_{B,\Fb,\tau}^{nil}(G,\mu)
\end{eqnarray}
as well as $\B_{A,\sigma}^{nil}(G)\rightarrow \B_{B,\tau}^{nil}(G)$ and $\B_{A,\sigma}(G)\rightarrow \B_{B,\tau}(G)$. 
If $\varkappa$ is injective, it is easy to see that all of the aforementioned functors are faithful.\\

Observe that the triple $(W(A),W(\Fa)+I(A),F)$ is a $W(k_0)$-frame whenever $\Fa$ is a $p$-adically closed pd-ideal in a $p$-adically separated 
and complete torsion-free $W(k_0)$-algebra $A$. We write $\varkappa_\sigma:(A,\Fa,\sigma)\rightarrow(W(A),W(\Fa)+I(A),F)$ for Cartier's diagonal 
map, which is associated to our Frobenius lift $\sigma$. This Frobenius lift induces by functoriality of the Witt vectors, a morphism on $W(A)$, which we will still denote by $\sigma$. Note that $\varkappa_{\sigma}$ identifies $A$ with the equalizer of the two endomorphisms $F$ and $\sigma$ on $W(A)$. Concerning a change of frame via $\varkappa_\sigma$ we have the following result, which is one of the main ingredients of for the proof of the main theorem:

\begin{prop}
\label{independence07}
Fix a $W(k_0)$-frame $(A,\Fa,\sigma)$. The frame change
\begin{equation}
\label{independence05}
\B_{A,\Fa,\sigma}^{nil}(G,\mu)\rightarrow \B_{W(A),W(\Fa)+I(A),F}^{nil}(G,\mu),
\end{equation}
which is induced from Cartier's diagonal morphism
$$\varkappa_\sigma:(A,\Fa,\sigma)\rightarrow(W(A),W(\Fa)+I(A),F)$$
is an equivalences of categories.
\end{prop}
\begin{proof}
Observe that we have $W(A)/(W(\Fa)+I(A))\cong A/\Fa$, so that the $P_\mu$-torsor $Q_0$ which is implicit in giving a window (being the 2nd datum in the quadruple of definition 
\eqref{wind02}) is not altered by the functor \eqref{independence05}. In particular a $(G,\mu)$-window over $(A,\Fa,\sigma)$ is banal if and only if its base change under 
$\varkappa_\sigma$ is a banal $(G,\mu)$-window over $(W(A),W(\Fa)+I(A),F)$. We are going to first prove the equivalence of the frame change functor restricted to the banal objects and then prove the general case by descent to this case. So let us consider the banal pendant of the functor \eqref{independence05}, which reads:\begin{equation}
\label{independence06}
[C_A^{nil}(G,\mu){^F\mu(p)}/_\sigma\Gamma_\mu^{A,\Fa}]\rightarrow[C_{W(A)}^{nil}(G,\mu){^F\mu(p)}/_F\Gamma_\mu^{W(A),W(\Fa)+I(A)}]
\end{equation}

Recall that for any frame $(B,\Fb,\tau)$ we have a canonical identification 
\[[C_B^{nil}(G,\mu){^F\mu(p)}/_\sigma\Gamma_\mu^{B,\Fb}]\cong[C_B^{nil}(G,\mu)/_{\hat\Phi_\mu^{B,\tau}}\Gamma_\mu^{B,\Fb}]\]We are going to work with the description on the right hand side.
Our proof starts with the full faithfulness of \eqref{independence06}. Its faithfulness is clear, from $\varkappa_\sigma$ inducing 
an inclusion $\Gamma_\mu^{A,\Fa}\into\Gamma_\mu^{W(A),W(\Fa)+I(A)}$. Let us consider $O,U\in C_A^{nil}(G,\mu)$, so that 
$\varkappa_\sigma(O),\varkappa_\sigma(U)\in C_{W(A)}^{nil}(G,\mu)$. Suppose that $k\in\Gamma_\mu^{W(A),W(\Fa)+I(A)}$ 
satisfies $$\varkappa_\sigma(O)=k^{-1}\varkappa_\sigma(U)\hat\Phi_\mu^{W(A),F}(k)$$ Applying $F$ and $\sigma$ yields
$${^F\varkappa_\sigma(O)}={^Fk}^{-1}{^F\varkappa_\sigma(U)}\hat\Phi_{^F\mu}^{W(A),F}({^Fk})$$
and
$${^F\varkappa_\sigma(O)}=\sigma(k)^{-1}{^F\varkappa_\sigma(U)}\hat\Phi_{^F\mu}^{W(A),F}(\sigma(k))$$
as $^F\varkappa_\sigma(O)=\sigma(\varkappa_\sigma(O))$ and $^F\varkappa_\sigma(U)=\sigma(\varkappa_\sigma(U))$. 
In view of proposition \ref{GMZCF} we are allowed to deduce $^Fk=\sigma(k)$ from $^Fk\equiv\sigma(k)\mod W(pA)$ together 
with the fact that $pA$ is a pd-ideal of $A$. It follows that $k$ lies in the image of $\varkappa_\sigma$. At last, notice that
\[\begin{tikzcd}
\Gamma_\mu^{A,\Fa}\arrow{r}{\varkappa_{\sigma}}\arrow[d,hook]&\Gamma_\mu^{W(A),W(\Fa)+I(A)}\ar[d,hook]\\
G(A)\arrow{r}[below]{\varkappa_{\sigma}}&G(W(A))
\end{tikzcd}
\]

is a Cartesian diagram, so that we have $k\in\varkappa_\sigma(\Gamma_\mu^{A,\Fa})$, which is what we wanted.\\

We continue with the essential surjectivity of \eqref{independence06}. Let us consider some $U\in C_{W(A)}^{nil}(G,\mu)$. Observe 
that we have $\tau(U)\equiv F(U)\mod W(pA)$, so that by proposition \ref{GMZCF} and the fact that $pA$ is a pd-ideal implies there exists some $k\in G(W(pA))$ with 
$$\tau(U)=k^{-1}F(U)\hat\Phi_\mu^{W(A),F}(k)$$

Consider, for each $i$, the ghost coordinate 
$k_i:=w_i(k)$ and the element $h_i:=k_{i-1}\tau(k_{i-2})\dots\tau^{i-1}(k_0)$ of $G(A)$. Since $h_{i+1}\tau(h_i)^{-1}=k_i\equiv1\pmod{p^{i+1}}$, there exists an 
element $h\in G(W(A))$ whose $i$th ghost coordinate agrees with $h_i$ for any $i$, in fact $h$ lies in $G(I(A))\subset\Gamma_\mu^{W(A),W(\Fa)+I(A)}$, since 
$h_0=1$. Now it follows that $U':=h^{-1}U\hat\Phi_\mu^{W(A),F}(h)$ lies in the image of $\varkappa_\sigma$, as $^FU'=\tau(U')$ in view of $F(h)\tau(h)^{-1}=k$. This achieves the proof in the banal case.\\

We are now in the position to return to the functor \eqref{independence05}, and again we start with the full faithfulness. Consider a 
$\B_{W(A),W(\Fa)+I(A),F}^{nil}(G,\mu)$-isomorphism $h$ between $\B_{A,\Fa,\sigma}^{nil}(G,\mu)$-objects $(P',Q_0',\eta_0',b')$ and $(P,Q_0,\eta_0,b)$. According to \cite[Corollaire (17.16.3.ii)]{ega4} we may choose a \'etale faithfully flat $A/(pA+\Fa)$-algebra $\tilde A_0$ satisfying $Q_0'(\tilde A_0)\neq\emptyset\neq Q_0(\tilde A_0)$. 
Let us write $\tilde A_i$ for the \'etale faithfully flat $A/p^iA$-algebra satisfying $\tilde A_i/(p\tilde A_i+\Fa\tilde A_i)\cong\tilde A_0$, which exists up to a unique 
isomorphism (cf. \cite[Th\'eor\`eme (18.1.2)]{ega4}). The uniqueness implies that there are natural morphisms $\tilde{A}_{i+1}\to\tilde{A}_{i}$ inducing an isomorphism $\tilde{A}_{i+1}/p^{i}\cong \tilde{A}_{i}$. By lemma \ref{descent02} below it is possible to extend $\sigma$ to a Frobenius lift $\tilde\sigma$ on 
the $p$-adically formally \'etale, separated and complete $A$-algebra $\tilde A:=\invlim_i\tilde A_i$. At last writing $\tilde\Fa$ for the $p$-adic closure of $\Fa\tilde A$, 
yields a $W(k_0)$-linear $p$-adic covering $\varkappa:(A,\Fa,\sigma)\rightarrow(\tilde A,\tilde\Fa,\tilde\sigma)$. We write $(\tilde P',\tilde Q_0',\tilde\eta_0',\tilde b')$ 
(resp. $(\tilde P,\tilde Q_0,\tilde\eta_0,\tilde b)$) for the $\B_{\tilde A,\tilde\Fa,\tilde\sigma}^{nil}(G,\mu)$-objects that arise from $(P',Q_0',\eta_0',b')$ (resp. 
$(P,Q_0,\eta_0,b)$) by base change along $\varkappa$. From the full faithfulness of \eqref{independence06} we deduce an isomorphism, say 
$\tilde h:(\tilde P',\tilde Q_0',\tilde\eta_0',\tilde b')\rightarrow(\tilde P,\tilde Q_0,\tilde\eta_0,\tilde b)$ whose Cartier diagonal agrees with the base change of 
$h$ to the frame $(W(\tilde A),W(\tilde\Fa)+I(\tilde A),F)$. We have to show that $\tilde h$ descends to $(A,\Fa,\sigma)$. We therefore have to check that its images in 
$(\tilde A_{[2]},\tilde\Fa_{[2]},\tilde\sigma_{[2]})$ 
agree (see notations at the beginning of subsection \ref{subsectionDescent}), which can be done  using the previously established faithfulness of $\varkappa_{\tilde\sigma_{[2]}}$.\\

It remains to check the essential surjectivity of \eqref{independence05}, so consider a $(G,\mu)$-window $(P',Q_0',\eta_0',b')$ over $(W(A),W(\Fa)+I(A),F)$. 
Let us choose a $W(k_0)$-linear $p$-adic covering $\varkappa:(A,\Fa,\sigma)\rightarrow(\tilde A,\tilde\Fa,\tilde\sigma)$ such that $(P',Q_0',\eta_0',b')$ 
becomes banal over $(W(\tilde A),W(\tilde\Fa)+I(\tilde A),F)$, here notice that $W(\tilde A)/W(\tilde\Fa)+I(\tilde A)\cong\tilde A/\tilde\Fa$, so that we can construct 
$\tilde A$ by means of \cite[Corollaire (17.16.3.ii)]{ega4} and \cite[Th\'eor\`eme (18.1.2)]{ega4}. The self-explanatory descent argument is completed by means of the faithfulness of $\varkappa_{\tilde\sigma_{[3]}}$ and the full faithfulness of $\varkappa_{\tilde\sigma_{[2]}}$.
\end{proof}

\begin{lem}
\label{descent02}
Let $A$ be a ring and $B$ be a $p$-adically formally smooth, separated and complete $A$-algebra. 
Let $\sigma:A\rightarrow A$ be a Frobenius lift. Then there exists a Frobenius lift $\tau:B\to B$ that restricts to $\sigma$ on $A$.
\end{lem}

\begin{proof}
In view of $B$ being $p$-adically separated and complete, we should construct a compatible system of Frobenius lifts $\tau_{n}:B/p^{n}\to B/p^{n}$, extending $\sigma_{n}:A/p^{n}\to A/p^{n}$. We can therefore assume that 
$p^nA$ vanishes for some $n$. The case $n=1$ is trivial. By induction we may assume the result for all formally 
smooth $A/p^{n-1}A$-algebras, so that we have a Frobenius lift extending $\sigma_{n-1}$
$$\tau_{n-1}:B/p^{n-1}\rightarrow B/p^{n-1}$$
The homomorphism $\tilde\tau_{n-1}:B\onto B/p^{n-1}\arrover{\tau_{n-1}} B/p^{n-1}$ fits into a commutative diagram

\[
\begin{tikzcd}
A\ar[r]\arrow{d}[left]{\sigma}&B\arrow[dr, dashed, "\tau_{n}"]\arrow{r}{\tilde\tau_{n-1}}&B/p^{n-1}\\
A\ar[rr]&&B\arrow[u, two heads]
\end{tikzcd}
\]

in which the formal smoothness yields the desired oblique morphism $\tau_n:B\rightarrow B$ filling the diagram.
\end{proof}

We are now in a position to prove the following:

\begin{thm}
\label{main01}
Let $(A,\Fa,\sigma)$ be a $W(k_0)$-frame and $(G,\mu)$ a $W(k_0)$-display datum, where $k_0$ 
is a finite field. The natural functor from $\B_{A,\Fa,\sigma}(G,\mu)$ to the $A/\Fa$-valued objects of $\B(G,\mu)$ restricts to an equivalence from the groupoid of $(G,\mu)$-windows over $(A,\Fa,\sigma)$ to the groupoid of $(G,\mu)$-displays 
over $A/\Fa$ of which the $\pmod p$-reduction is adjoint nilpotent. If $\Fa$ is $p$-adically open, one has in particular:
$$\B_{A,\Fa,\sigma}^{nil}(G,\mu)\isoto\B^{nil}(G,\mu)(A/\Fa)$$ 
\end{thm}
\begin{proof}
In view of theorem \ref{triple05} and $\varkappa_\sigma:\B_{A,\Fa,\sigma}^{nil}(G,\mu)\isoto \B_{W(A),W(\Fa)+I(A),F}^{nil}(G,\mu)$ (i.e. proposition 
\ref{independence07}) it only remains to check \eqref{wind12}. However, this follows from our reformulation i.e. corollary \ref{poshI}. We merely have to observe that
\begin{itemize}
\item
$\Tor_{W(A)}(G)\isoto\Tor_A(L^+G)$ by lemma \ref{torsor02}.
\item
$\Tor_{W(A)}(\hat P_\mu^{W(\Fa)+I(A)})\isoto\Tor_A(\hat H_\mu^\Fa)$ again by lemma \ref{torsor02} in conjunction with part (i) of lemma \ref{triple04}
\item
The $2$-commutativity of the diagram:
\end{itemize}
\[
\begin{tikzcd}
\Tor_{W(A)}(\hat P_\mu^{W(\Fa)+I(A)})\ar[r]\arrow[d,"F" left]&\Tor_A(\hat H_\mu^\Fa)\arrow[r, "\Psi_\mu^\Fa"]&\Tor_A(L^+G) \\
\Tor_{W(A)}(\hat P_{^F\mu}^{pW(k_0)})\arrow[rr,"\text{mo}_\mu^\sigma" below]&&\Tor_{W(A)}(G)\arrow[u,"\cong" right]
\end{tikzcd}\]
\end{proof}

\appendix

\section{Groupoids}
\label{app1st}
Let $\Gamma$ be an abstract group acting from the right on a set $X$. We denote by $[X/\Gamma]$ the following category: 
objects are elements of $X$ and if $x,y\in X$, we set $\Hom_{[X/\Gamma]}(x,y):=\{\gamma\in\Gamma\,\mid\, x=y\cdot\gamma\}$. 
This defines a groupoid, as every morphism has an inverse.\\

Suppose that $\Xi,\Phi:\Gamma\rightarrow M$ are homomorphisms to another abstract group $M$. By $(\Xi,\Phi)$-conjugation we mean the action
\begin{equation}
\label{Frobconj}
M\times\Gamma\rightarrow M;\quad (m,\gamma)\mapsto\Xi(\gamma)^{-1}m\Phi(\gamma),
\end{equation}
and we let $[M/_{\Xi,\Phi}\Gamma]$ be the resulting quotient groupoid, where $M$ is regarded as a right $\Gamma$-set, according to the 
$(\Xi,\Phi)$-conjugation \eqref{Frobconj}. If $M$ contains $\Gamma$ as a subgroup and $\Xi$ is the inclusion, we drop $\Xi$ from the notation. 
Note that this construction is functorial in quadruples $(\Xi,\Phi,\Gamma,M)$ in the obvious way. Whenever some subset $X\subset M$ is invariant 
under $(\Xi,\Phi)$-conjugation, we also write $[X/_{\Xi,\Phi}\Gamma]$ for the full subcategory of $[M/_{\Xi,\Phi}\Gamma]$ whose class of objects is $X$.

\section{Torsors}
\label{app2nd}
\subsection{Fpqc sheaves}

For every scheme $S$ we let $\Scheme_S$ be the category of schemes over the base $S$. Unless otherwise said we endow $\Scheme_S$ 
with the so-called fpqc topology, which is the Grothendieck topology generated by all faithfully flat and quasi-compact coverings together with the 
usual Zariski ones (cf. \cite[Expos\'e IV 6.3]{sga1}). Recall that a covering $\V=\{V_j\}_{j\in J}$ of some $S$-scheme $U$ is said to refine another 
$\U=\{U_i\}_{i\in I}$ one, if for each $j\in J$ there exists $i\in I$ such that $V_j\rightarrow U$ factors through $U_i\rightarrow U$. Recall also that 
a presheaf on $S$ is nothing but a contravariant set-valued functor on $\Scheme_S$. Let us recall the sheaf concepts for this particular site:

\begin{notation}
\label{NotationsTorsors01}
\begin{itemize}
\item
$\U\times\V:=\{U_i\times_UV_j\}_{i\in I,\,j\in J}$ is the ``canonical common refinement'' of coverings $\U=\{U_i\}_{i\in I}$ and $\V=\{V_j\}_{j\in J}$ of some $S$-scheme $U$. 
\item
For every presheaf $X:\Scheme_S^{\rm op}\rightarrow\Ens$ and every covering $\U=\{U_i\}_{i\in I}$ of $U$, we let $X(\U):=\prod_{i\in I}X(U_i)$. Moreover, 
\begin{equation}
\label{NotationsTorsors02}
\proj_1,\,\proj_2:X(\U)\rightrightarrows X(\U\times\U)
\end{equation}
denote the projections to the two factors. Observe that $U_{ij}:=U_i\times_UU_j$ are the ``intersections'' and that 
$X(\U\times\U)=\prod_{i,\,j\in I}X(U_{ij})$.
\item
A presheaf $X:\Scheme_S^{\rm op}\rightarrow\Ens$ is called \emph{separated} if the natural map 
\begin{equation}
\label{NotationsTorsors03}
X(U)\rightarrow X(\U)
\end{equation} 
is injective for all coverings $\{U_i\}_{i\in I}$ of $U$.
\item
A separated presheaf $X:\Scheme_S^{\rm op}\rightarrow\Ens$ is called a \emph{sheaf} if \eqref{NotationsTorsors03} agrees 
with the equalizer of \eqref{NotationsTorsors02}, whenever $\{U_i\}_{i\in I}$ is a covering of $U$.
\end{itemize}
\end{notation}

Notice that the sheaf property implies the convenient formula $X(\dot\bigcup_{i\in I}U_i)=X(\U)$, since $\U=\{U_i\}_{i\in I}$ is a Zariski covering of 
$\dot\bigcup_{i\in I}U_i$. A homomorphism between sheaves or (separated) presheaves is defined to be a transformation of contravariant functors.

\begin{rem}
\label{glueing01}
Let $\{U_i\}_{i\in I}$ be a covering of $U$ and let $X_i$ be a sheaf on $U_i$. By a glueing datum on the family $\{X_i\}_{i\in I}$ 
one means a family of isomorphisms $\psi_{ij}$ between sheaves $X_j\vert_{U_{ij}}$ and $X_i\vert_{U_{ij}}$ on $U_{ij}$ satisfying 
$\psi_{ij}\circ\psi_{jk}=\psi_{ik}$. The family of restrictions to $U_i$ of a sheaf $X$ on $U$ possesses a tautological glueing datum. One knows that:
\begin{itemize}
\item[(a)]
The set of homomorphisms $\alpha$ between sheaves $X'$ and $X$ on $U$ can be recovered as the family of homomorphisms 
$\alpha_i$ between the respective restrictions $X'_i=X'\vert_{U_i}$ and $X_i=X\vert_{U_i}$ which preserve their glueing data 
$\psi_{ij}'$ and $\psi_{ij}$ respectively in the sense that $\alpha_i\vert_{U_{ij}}\circ\psi_{ij}'=\psi_{ij}\circ\alpha_j\vert_{U_{ij}}$ holds.
\item[(b)]
To every glueing datum $\{\psi_{ij}\}_{i,\,j\in I}$ on $\{X_i\}_{i\in I}$ there is a sheaf 
$X$ on $U$ together with a family $e_i:X_i\stackrel{\cong}{\rightarrow}X\vert_{U_i}$ satisfying 
$e_i\vert_{U_{ij}}\circ\psi_{ij}=e_j\vert_{U_{ij}}$. The pair $(X,\{e_i\}_{i\in I})$ is unique up to unique isomorphism.
\end{itemize}
\end{rem}

\subsection{Torsors}
\label{ExampleCocyclGroup}
In this subsection we recall some basic definitions and facts on torsors. The main references are \cite[Chapitre III]{giraud} and \cite[Chapter III]{milne1}. So let $G$ be a 
fpqc sheaf in groups over $S$ and let $\U=\{U_i\}_{i\in I}$ be an fpqc covering of some $S$-scheme $U$. Let us endow the set $G(\U\times\U)$ with a right $G(\U)$-action via:
\[(\bold g,\vec h)\mapsto\proj_1(\vec h)^{-1}\,\bold g\,\proj_2(\vec h)\]whenever $\bold g\in G(\U\times\U)$ and $\vec h\in G(\U)$. In the 
language of appendix \ref{app1st} this is just our $(\proj_1,\proj_2)$-conjugation. Recall that a family $\bold g:=\{g_{ij}\}_{i,\,j\in I}\in G(\U\times\U)$ 
satisfying the cocycle relation
\[\forall i,j,k:\quad (g_{ij}\vert_{U_{ijk}})(\,g_{jk}\vert_{U_{ijk}})=g_{ik}\vert_{U_{ijk}}\]
is called a \emph{1-cocycle} for $\U$, according to \cite[Chapitre III, Definition 3.6.1.(i)]{giraud}. The set of all 1-cocycles for $\U$ forms a $G(\U)$-invariant 
subset, which we denote by $\check{Z}^1(\U/S,G)\subset G(\U\times\U)$. Two 1-cocycles $\bold g=\{g_{ij}\}_{i,\,j\in I}$ and $\bold f=\{f_{ij}\}_{i,\,j\in I}$ 
for $\U$ are \emph{cohomologuous} if they lie in the same $G(\U)$-orbit (i.e. such that there exists some $\vec h=\{h_i\}_{i\in I}\in G(\U)$ satisfying
\begin{align}
\label{1-CocyclCoh}
(h_i\vert_{U_{ij}})g_{ij}=f_{ij}(h_j\vert_{U_{ij}})
\end{align}
in accordance with \cite[Chapitre III, Definition 3.6.1.(ii)]{giraud}).\\

The set of \emph{cohomology classes} of 1-cocycles modulo this equivalence relation is denoted by $\check{H}^1(\U/S,G)$. Notice that 
$\check{H}^1(\U/S,G)$ has a distinguished element, namely, the class of trivial 1-cocycles: Here, a cocycle is called \emph{trivial} if it lies in the $G(\U)$-orbit 
of the neutral element (i.e. if it can be written as $\bold g=\{(e_i^{-1}\vert_{U_{ij}})(e_j\vert_{U_{ij}})\}_{i,\,j\in I}$, where $\{e_i\}_{i\in I}=\vec e\in G(\U)$).

\begin{defn}
Let $G$ be an fpqc sheaf of groups over $S$. A \emph{pseudo $G$-torsor} on an $S$-scheme $U$ is an fpqc sheaf $X$ on $U$ on which $G$ acts freely and transitively from the right, i.e. such that the morphism 
\[X\times G\vert_U\to X\times X,\quad (x,g)\mapsto (x, x\cdot g)\] 
is an isomorphism. For any two pseudo $G$-torsors $X$ and $Y$ over $U$, we write $\Hom_G(X,Y)$ for the $G$-equivariant morphism between these two 
fpqc-sheaves on $U$. 
\end{defn}

\begin{rem}
\label{glueing02}
Recall that for every pseudo $G$-torsor $X$ on $U$ any global 
section $s\in X(U)$ gives rise to a canonical isomorphism:
\[G\vert_U\isoto X,\quad g\mapsto s\cdot g\] 
Therefore we have canonical identifications: 
$$\Isom_{G\vert_U}(X,G\vert_U)\cong\Isom_{G\vert_U}(G\vert_U,X)\cong\Hom_{G\vert_U}(G\vert_U,X)\cong X(U)$$
We deduce that a pseudo $G$-torsor whose set of global sections is 
not empty must be isomorphic to $G\vert_U$ regarded as a pseudo 
torsor under itself. If this happens the pseudo $G$-torsor $X$ 
is called \emph{(globally) trivial}. Also, notice that we have 
a canonical identification $\Aut_{G\vert_U}(G\vert_U)=G(U)$, but 
do notice that the isomorphism $\Aut_{G\vert_U}(X)=G(U)$ for 
an arbitrary trivial $G$-torsor $X$ may depend on the choice 
of $G\vert_U\cong X$, i.e. on the base point $s\in X(U)$.
\end{rem}

\begin{defn}
Let $G$ be an fpqc sheaf of groups over $S$. A pseudo $G$-torsor 
$X$ over an $S$-scheme $U$ is called a \emph{$G$-torsor} if and 
only if there exists a covering $\U=\{U_{i}\}_{i\in I}$ of $U$ 
such that each $X\vert_{U_i}$ is trivial. We denote by $\Tor_U(G)$ 
the category of $G$-torsors over $U$ and by $\Tor_{\U/U}(G)$ the 
full subcategory of $G$-torsors that are trivialized by a given covering $\U$. We let $\phs_U(G)$ and $\phs_{\U/U}(G)$ 
be the respective classes of isomorphism classes of objects in 
these categories.\\

At last, we write $\vec{\Tor}_{\U/U}(G)$ for the category 
of $G$-torsors over $U$ that are equipped with a specific 
trivialization over $\U$.
\end{defn}

It is folklore that the set of isomorphism classes of torsors are in bijection with the first cohomology with non-abelian coefficients. 
We would like to make this statement precise, the bijection explicit and present a variant of it. To this end we fix a covering 
$\{U_i\}_{i\in I}=\U$ of $U$. An element $\{g_{ij}\}_{i,\,j\in I}=\bold g\in\check{Z}^1(\U/U,G)$ may be interpreted as a glueing 
datum on the family $X_i:=G\vert_{U_i}$ of trivial $G$-torsors on $U_i$. Whence we obtain a $G$-torsor $X$ over $U$ together 
with a family $\{x_i\}_{i\in I}=\vec x\in X(\U)$. We denote the assignment thus obtained by:
$$\Fc:\check{Z}^1(\U/U,G)\rightarrow\vec{\Tor}_{\U/U}(G);\,\bold g\mapsto(X,\vec x)$$
Notice that we have $x_i\vert_{U_{ij}}g_{ij}=x_j\vert_{U_{ij}}$ in this setting.

\begin{prop}
\label{PropTorCoho}
Let $G$ be an fpqc sheaf of groups on a scheme $S$ and fix a covering $\U=\{U_i\}_{i\in I}$ of an $S$-scheme $U$:
\begin{itemize}
\item[(i)]
The assignment $\Fc$ induces a bijection between the set $\check{H}^1(\U/U,G)$ and the set of isomorphism classes of the category $\Tor_{\U/U}(G)$.
\item[(ii)]
The assignment $\Fc$ induces a bijection between $\check{Z}^1(\U/U,G)$ and the set of isomorphism classes of the (discrete!) category $\vec{\Tor}_{\U/U}(G)$.
\item[(iii)]
The assignment $\Fc$ induces an equivalence $$[\check{Z}^1(\U/U,G)/_{p_1,\,p_2}G(\U)]\cong\Tor_{\U/U}(G)$$ in the sense that there exists a commutative diagram

\[
\begin{tikzcd}
\check{Z}^1(\U/U,G)\ar[d]\arrow[r,"\Fc"]& \vec{\Tor}_{\U/U}(G)\ar[d]\\
{[}\check{Z}^1(\U/U,G)/_{\proj_1,\,\proj_2}G(\U){]}\arrow[r,"\cong" below] &\Tor_{\U/U}(G)
\end{tikzcd}
\]

both of whose horizontal arrows are category equivalences.
\end{itemize}  
\end{prop}
\begin{proof}
We could just refer to \cite[Chapitre III, Proposition 3.6.3.(i)]{giraud}, but for the convenience of the reader we explain in our own words what is going on. Part (ii) follows from the remarks \ref{glueing01} and \ref{glueing02}. We continue with the construction of a 
fully faithful functor from 
$[\check Z^1(\U/U,G)/_{\proj_1,\,\proj_2}G(\U)]$ to $\Tor_{\U/U}(G)$: 
All we have to do is exhibit a bijection
\begin{align}
\label{glueing04}
\{\vec h\in G(\U)\mid\,\proj_1(\vec h)\, \bold g=\bold f\,\proj_2(\vec h)\}\stackrel{\cong}{\rightarrow}\Hom_{\Tor_{\U/U}(G)}(X,Y)
\end{align}
for any two instances 
\begin{eqnarray*}
\check Z^1(\U/U,G)\ni\bold g=\{g_{ij}\}_{i,\,j\in I}&\stackrel{\Fc}{\mapsto}(X,\vec x)\\
\check Z^1(\U/U,G)\ni\bold f=\{f_{ij}\}_{i,\,j\in I}&\stackrel{\Fc}{\mapsto}(Y,\vec y)
\end{eqnarray*} 
of our glueing construction $\Fc$. Here, $\vec x\in X(\U)$ and $\vec y\in Y(\U)$ are respective trivalisations of 
$G$-torsors $X$ and $Y$ over $U$, as objects of the category $\vec{\Tor}_{\U/U}(G)$ are ``pointed'' $G$-torsors over $U$.
These elements satisfy:
\begin{eqnarray*}
&&x_i\vert_{U_{ij}}g_{ij}=x_j\vert_{U_{ij}}\\
&&y_i\vert_{U_{ij}}f_{ij}=y_j\vert_{U_{ij}}
\end{eqnarray*} 
For any $i\in I$ and $h\in G(U_i)$ we will now write $\alpha_i(h):X_i\rightarrow Y_i$ for the unique $G$-morphism which sends $x_i$ to $y_ih$. Let us see for which 
families $\{h_i\}_{i\in I}=\vec h\in G(\U)$ one obtains a global morphism $\alpha:X\rightarrow Y$ from those $\alpha_i(h_i)$'s: The restriction of $\alpha_i(h_i)$ to $U_{ij}$ 
sends $x_j\vert_{U_{ij}}$ to: 
$$y_i\vert_{U_{ij}}h_i\vert_{U_{ij}}g_{ij}$$ 
The restriction of $\alpha_j(h_j)$ to $U_{ij}$ sends $x_j\vert_{U_{ij}}$ to: 
$$y_i\vert_{U_{ij}}f_{ij}h_j\vert_{U_{ij}}$$
It follows that $\alpha_i(h_i)$ and $\alpha_j(h_j)$ agree on the overlap $U_{ij}$ if and only if $h_i\vert_{U_{ij}}g_{ij}=f_{ij}h_j\vert_{U_{ij}}$, 
in which case we do obtain our global homomorphism. This yields the bijection \eqref{glueing04}, and we 
leave it to the reader to check that it is compatible with compositions. At last (i) is a mere corollary to (iii). 
\end{proof}

\begin{rem}
\label{glueing05}
Let $\U$ be a covering of an $S$-scheme $U$ and let $\V$ be a refinement of $\U$. There are canonical arrows
\begin{eqnarray*}
&&\Tor_{\U/U}(G)\rightarrow\Tor_{\V/U}(G)\\
&&\phs_{\U/U}(G)\rightarrow\phs_{\V/U}(G),
\end{eqnarray*}
which do not depend on the specific factorisations $V_j\rightarrow U_i$, which are implicitly 
present when saying that $\V$ refines $\U$! Taking the direct limit over all coverings, we obtain 
\begin{eqnarray*}
&&\dirlim_{\U}\Tor_{\U/U}(G)\stackrel{\cong}{\rightarrow}\Tor_U(G)\\
&&\dirlim_{\U}\phs_{\U/U}(G)\stackrel{\cong}{\rightarrow}\phs_U(G),
\end{eqnarray*} 
and all of these maps are ``pointed'' in the sense that they preserve the distinguished element (cf. \cite[Remarque 3.6.5]{giraud}). It is somewhat more fiddly to describe the 
whole groupoid $\Tor_U(G)$ in terms of ``non-abelian Cech cohomology $1$-classes''. Rather than taking the direct limit we consider the disjoint union over all coverings
$$\check{Z}^1(U,G):=\dot\bigcup_\U\check{Z}^1(\U/U,G),$$ 
which is the proper class consisting of pairs $(\U,\bold g)$ where $\U$ is a covering of $U$ and $\bold g\in\check{Z}^1(\U/U,G)$. Next we define a groupoid structure 
on $\check{Z}^1(U,G)$ by decreeing that a homomorphism from $(\U,\{g_{ij}\}_{i,\,j\in I})$ to $(\V,\{f_{kl}\}_{k,\,l\in J})$ shall consist of $\vec h\in G(\U\times\V)$ satisfying 
$$(h_{ik}\vert_{U_{ij}\times_UV_{kl}})(g_{ij}\vert_{U_{ij}\times_UV_{kl}})=(f_{kl}\vert_{U_{ij}\times_UV_{kl}})(h_{jl}\vert_{U_{ij}\times_UV_{kl}})$$
This groupoid is equivalent to  $\Tor_U(G)$, but it does not fit into the framework of appendix \ref{app1st}.
\end{rem}

\subsection{Twists and change of group}
We will say that $X$ is a $G$-sheaf if each $X(U)$ is endowed with the structure of a left $G(U)$-set, functorially in 
$U\in\Ob_{\Scheme_S}$. In several settings we need to twist $G$-sheaves with $G$-torsors. The ultimate reference 
for this is \cite[Chapitre III, Proposition 2.3.1]{giraud}, but again we would like to spell out a version that fits the style of our paper: 

\begin{prop}
\label{FuncChangGroup1}
Let $G$ be an fpqc sheaf of groups on $S$, let $X$ be a $G$-sheaf and let $P$ be a $G$-torsor. For any $S$-scheme $U$ consider the equivalence relation
$$(p',x')\sim_G(p,x)\Leftrightarrow\exists g\in G(U):\;(p',x')=(pg^{-1},gx)$$
on the set $P(U)\times X(U)$. Then the presheaf 
\begin{align*}
P\times X/\sim_G:\Scheme_S\rightarrow\Ens:\,U\mapsto P(U)\times X(U)/\sim_G
\end{align*}
possesses a sheafification denoted by $P*^GX$.
\end{prop}
\begin{proof}
Let us choose a trivialization of $P$ over some covering $\U=\{U_i\}_{i\in I}$ of $S$, and consider the $1$-cocycle 
$\{g_{ij}\}_{i,\,j\in I}=\bold g\in\check{Z}^1(\U/U,G)$ that goes with it. Let $X_i$ be the pull-back of $X$ to $U_i$. The family of maps
$$\psi_{ij}:X_j\vert_{U_{ij}}\isoto X_i\vert_{U_{ij}};\,x\mapsto g_{ij}x$$
may be interpreted as a glueing datum on the family $\{X_i\}_{i\in I}$ in the sense of remark \ref{glueing01}. This yields a sheaf $P*^GX$, 
which is equipped with a canonical morphism $P\times X\rightarrow P*^GX$, as $P$ is obtained by an analogous glueing construction.
\end{proof}

Let $G$ and $H$ be a sheaves of groups over $S$ and let $\phi:G\rightarrow H$ be a homomorphism between them. 
Since $H$ may be regarded as a $G$-sheaf proposition \ref{FuncChangGroup1} can be used to construct a functor
\begin{equation}
\label{FuncChangGroup4}
\Tor_S(G)\rightarrow\Tor_S(H);\,P\mapsto P*^{G,\phi}H
\end{equation}
It is noteworthy, that the isomorphism class of \eqref{FuncChangGroup4} depends only on the conjugacy class of 
$\phi$. Indeed if $Q$ is any $H$-torsor and $h$ any global section of $H$, then there is a canonical isomorphism
\begin{equation}
\label{FuncChangGroup6}
Q\isoto Q*^{H,\intaut(h)}H;\,q\mapsto q*h
\end{equation}
as $\intaut(h):H\rightarrow H$ is the inner automorphism 
$g\mapsto h\vert_U\cdot g\cdot h\vert_U^{-1}$ for varying $U\in\Ob_{\Scheme_S}$ and $g\in H(U)$.

\subsection{Cartesianism of $\Tor_S$}
\label{FuncChangGroup3}

It turns out useful to have a criterion for the Cartesianism of $\Tor_S$:

\begin{lem}
\label{FuncChangGroup5}
Consider fpqc sheaves of groups $H$, $G$ and $G'$ on $S$ together with homomorphisms 
$\phi:G\rightarrow H$ and $\phi':G'\rightarrow H$. The following statements are equivalent:
\begin{itemize}
\item[(i)]
The composition of the group law $m_H:H^2\rightarrow H$ with the outer product 
$\phi\times\phi':G\times G'\rightarrow H^2$ induces an epimorphism of fpqc sheaves $G\times G'\onto H$.
\item[(ii)]
The natural $2$-commutative diagram

\[
\begin{tikzcd}
\Tor_S(G\times_{\phi,H,\phi'}G')\arrow[r, "\proj_2"]\arrow[d,"\proj_1"] & \Tor_S(G')\arrow[d,"\phi'"]\\
\Tor_S(G)\arrow[r,"\phi"] & \Tor_S(H)
\end{tikzcd}
\]

is $2$-Cartesian.
\end{itemize}
In particular (ii) holds if one of $\phi$ or $\phi'$ is already an epimorphism.
\end{lem}
\begin{proof}
Suppose that (i) holds and that we are given torsors $P$ (resp. $P'$) under $G$ (resp. $G'$) together with $\alpha:P'*^{G',\phi'}H\isoto P*^{G,\phi}H=:Q$. 
We have to construct a $G'':=G\times_{\phi,H,\phi'}G'$-torsor $P''$ with compatible isomorphisms $P''*^{G'',\proj_1}G\cong P$ and $P''*^{G'',\proj_2}G'\cong P'$. 
The sheaf $P'':=P\times_{\phi,Q,\alpha\circ\phi'}P'$ is a pseudo-torsor under $G''$. In order to check that $P''$ is a torsor we may assume that all of 
$P'$, $P$ and $P*^{G,\phi}H=Q$ are trivial. In view of (i) we may also assume that $\alpha$ is the identity, as it does no harm to replace it by 
$\phi(\gamma)\circ\alpha\circ\phi'(\gamma')$ for any $\gamma\in\Aut_G(P)$ and $\gamma'\in\Aut_{G'}(P')$. However $P=G$, $P'=G'$ and $\alpha=\id_G$ 
gives rise to $P''=G''$ and we are done. We leave the implication (ii)$\Rightarrow$(i) to the reader because it is not used in the text. 
\end{proof}
 
\subsection{Functoriality properties of $\Tor_S$}
\label{mini}
In fact we need a slight generalization of the obvious functoriality \eqref{FuncChangGroup4}: Let us write $G_{small}$ (resp. $H_{small}$) for the restrictions 
of fpqc groups $G$ (resp. $H$) on $S$ to the full subcategory $\Flat_S\subset\Scheme_S$ consisting of flat $S$-schemes. By a \emph{mini-morphism} from $G$ to $H$ 
we understand a transformation of functors $G_{small}\rightarrow H_{small}$. This is a family of homomorphisms
$$\phi_U:G(U)\rightarrow H(U)$$
where $U$ runs through the class of flat $S$-schemes, and which is functorial in the sense that the diagrams

\[
\begin{tikzcd}
G(V)\ar[d]\arrow[r,"\phi_V"]&H(V)\ar[d]\\
G(U)\arrow[r,"\phi_U" below]&H(U)
\end{tikzcd}
\]

commute for any (not necessarily flat) $S$-morphism between any flat $S$-schemes $V$ and 
$U$. Let $\U$ be a covering of $S$. Observe that our mini-morphism $\phi$ induces maps
\begin{eqnarray*}
&&\check{H}^1(\U/S,G)\rightarrow\check{H}^1(\U/S,H)\\
&&\check{Z}^1(\U/S,G)\rightarrow\check{Z}^1(\U/S,H),
\end{eqnarray*}
and a functor:
$$[\check{Z}^1(\U/S,G)/_{\proj_1,\,\proj_2}G(\U)]\rightarrow[\check{Z}^1(\U/S,H)/_{\proj_1,\,\proj_2}H(\U)]$$
Using transport of structure and proposition \ref{PropTorCoho} we obtain maps
\begin{eqnarray*}
&&\phs_{\U/S}(G)\rightarrow\phs_{\U/S}(H)\\
&&\vec{\Tor}_{\U/S}(G)\rightarrow\vec{\Tor}_{\U/S}(H),
\end{eqnarray*}
and a functor:
$$\Tor_{\U/S}(G)\rightarrow\Tor_{\U/S}(H)$$
At last, passage to the limit (cf. remark \ref{glueing05}) yields another map
$$\phs_S(G)\rightarrow\phs_S(H)$$
and yet another functor:
\begin{equation}
\label{awkward01}
\Tor_S(G)\rightarrow\Tor_S(H)
\end{equation}
By slight abuse of notation we continue to denote its effect on $G$-torsors $P$ by: $P\mapsto P*^{H,\phi}G$

\subsection{Representability of twists}

In the last part of this appendix we discuss representability issues. Whenever $\L$ is a line bundle on a $S$-scheme $X$, we write $\Aut_S(X,\L)$ for the set 
of pairs $(\alpha,m)$ consisting of an $S$-automorphism $\alpha:X\isoto X$ and an isomorphism $m:\alpha^*(\L)\isoto\L$. This is a group, the group law being: 
$(\alpha,m)(\beta,n):=(\alpha\circ\beta,n\circ\beta^*(m))$. Furthermore, we let $\uAut_S(X,\L)$ be the fpqc sheaf of groups on $S$ whose sections on an $S$-scheme $U$ are given by
$$\Gamma(U,\uAut_S(X,\L))=\Aut_U(X_U,\L_U),$$
where $\L_U$ stands for the pull-back of $\L$ along the projection of $X_U=X\times_SU$ to the first factor. Moreover, $\uAut_S(X)$ shall denote the 
fpqc sheaf of groups given by $U\mapsto\Aut_U(X_U)$. It goes without saying that there is a natural morphism $\uAut_S(X,\L)\rightarrow\uAut_S(X)$.

\begin{prop}
\label{FuncChangGroup2}
Consider an $S$-scheme $X$.
\begin{itemize}
\item[(i)]
Assume that $X\rightarrow S$ is affine and consider $G:=\uAut_S(X)$. Then, for any $G$-torsor $P$ on $S$, the twist 
$P*^GX$ is representable by a relatively affine $S$-scheme.
\item[(ii)]
Assume that $X$ is of finite type over $S$, let $\L$ be a relatively ample line bundle on the $S$-scheme $X$ and consider 
$G:=\uAut_S(X,\L)$. Then, for any $G$-torsor $P$, the twist $P*^GX$ is representable by a quasi-projective $S$-scheme.
\end{itemize}
\end{prop}
\begin{proof}
Part (i) can be proved using faithfully flat descent \`a la \cite[Th\'eor\`eme 2]{fga12}, which owes its effectivity to \cite[Th\'eor\`eme 1]{fga12} 
together with the mere fact that relatively affine $S$-schemes are relative spectra of Zariski sheaves of $\O_S$-algebras 
on $S$. Part (ii) is a straightforward consequence of \cite[Chapter 6.1, Theorem 7]{raynaud03}, see also \cite{fga}.
\end{proof}

\begin{cor}
\label{descent04}
Let $G$ be a relatively affine group scheme over $S$.
\begin{itemize}
\item[(i)]
Every $G$-torsor is representable by a relatively affine $S$-scheme. 
\item[(ii)]
For any two $G$-torsors $P$ and $Q$, the sheaf 
\begin{equation}
\label{descent07}
U\mapsto\Hom_{G\vert_U}(P\vert_U,Q\vert_U)
\end{equation}
 is representable by a relatively affine $S$-scheme.
\end{itemize}
\end{cor}
\begin{proof}
Observe that $P*^GG\cong P$, so that (i) is a special case of part (i) of proposition \ref{FuncChangGroup2}. Towards (ii) notice that the sheaf in question is 
isomorphic to $Q*^GP^{op}$, where $P^{op}$ stands for $P$ equipped with the left $G$-action: \[G\times P\to P;\quad (g,x)\mapsto x\cdot g^{-1}\] Again we are 
done by part (i) of proposition \ref{FuncChangGroup2}. Alternatively one could consult \cite[Chapter III, Theorem 4.3]{milne1}.
\end{proof}

We write $\innHom_G(P,Q)$ for the sheaf \eqref{descent07}.

\begin{lem}
\label{descent05}
Let $G$ be a flat group scheme over $S$ and let $P$ be a representable pseudo $G$-torsor over 
$S$. Then $P$ is a $G$-torsor if and only if the structural morphism $P\rightarrow S$ is faithfully flat. 
\end{lem}
\begin{proof}
This is because $P\rightarrow S$ is an fpqc covering in this case while $P$ becomes trivial over $P$, cf. \cite[Chapter III, Proposition 4.1]{milne1}.
\end{proof}

\begin{lem}
\label{descent06}
Let $i:H\hookrightarrow G$ be a closed immersion of smooth affine group schemes over a discrete valuation ring $\Lambda$. For every 
$G$-torsor $P$ over some $\Lambda$-scheme $S$, the fpqc-quotient sheaf $P/H_S$ is representable by a quasi-projective scheme over $S$.
\end{lem}
\begin{proof}
By Chevalley's trick $H$ can be recovered as the subgroup stabilizing some line $L\subset V$, where 
$$\rho:G\rightarrow\GL(V/\Lambda)$$  
is a suitable representation of $G$ on a free $\Lambda$-module $V$ of finite rank. Let $\Gamma$ be the stabilizer of the same line in the larger group $\GL(V/\Lambda)$, so that $\BP_\Lambda(V)$ represents the quotient of fppf sheaves $\GL(V/\Lambda)/\Gamma$. Our representation $\rho$ induces a natural morphism
\begin{equation}
\label{torsor07}
G/H\rightarrow\GL(V/\Lambda)/\Gamma,
\end{equation}
of which its source and target are the quotients in the category of fppf-sheaves. Here observe that \cite[Th\'eor\`eme 4.C]{onedimensional} (see also 
\cite[Expos\'e XVII, Appendice III, Th\'eor\`eme 1.3]{sga2} or \cite[Th\'eor\`eme 1, Example ii)]{raynaud01}) are granting both of these quotients to be representable 
by separated schemes of finite type over $\Lambda$. Two more remarks on this scenario are in order: First, it turns out that \eqref{torsor07} is a radicial morphism, 
as one can show by inspecting the injectivity of all geometric fibers. In particular it is separated and quasi-finite too. Second, recall that separated and quasi-finite morphisms are quasi-affine, as a consequence of Grothendieck's generalization of Zariski's main theorem, cf. \cite[Proposition (18.12.12)]{ega4} resp. \cite[Th\'eor\`eme (8.12.6)]{ega3} (see also \cite[Corollaire (18.12.13)]{ega4}). Let us write $R$ for the $\GL(V/\Lambda)$-torsor $P*^{G,\rho}\GL(V/\Lambda)$, observe that 
$$R*^{\GL(V/\Lambda)}\BP_\Lambda(V)\cong R/\Gamma_S$$
is (represented by) a familiar scheme, namely $\BP_S(\F)$, where $\F$ is the locally free $\O_S$-module corresponding to the $\GL(V/\Lambda)$-torsor $R$ via 
Hilbert 90. In fact, the whole of \eqref{torsor07} allows a twist with the $G$-torsor $P$ as can be seen by means of part (ii) of proposition \ref{FuncChangGroup2}. 
In this context it is relevant that there exists a $G$-equivariant ample line bundle $\L$ on $G/H$, namely the pull-back along \eqref{torsor07} of the canonical line 
bundle $\O_{\BP_\Lambda(V)}(1)$, which is a $\GL(V/\Lambda)$-equivariant ample line bundle on $\BP_\Lambda(V)$ (N.B.: pull-back along quasi-affine morphisms preserves ampleness).
\end{proof}

\begin{rem}
In the setting of the above lemma, one can identify the sections of $P/H_S$ over some $S$-scheme $T$ with the set of ``$H$-structures'' 
on the $G$-torsor $P\times_ST$, i.e. $H$-torsors $Q$ over $T$ equipped with an isomorphism $Q*^HG\cong P\times_ST$. 
\end{rem}

\end{document}